\documentclass[a4paper,11pt]{amsart}
\usepackage{amsfonts,amssymb,amsthm,cite,amsmath,amstext}
\usepackage{color}
\usepackage[colorlinks,linkcolor=black,citecolor=black]{hyperref}
\usepackage{comment}
\usepackage{framed}
\usepackage{calligra} % for dedication
\usepackage{tikz-cd}
\usepackage{tikz}
\newcommand{\black}{\color{black}}
\newcommand{\red}{\color{red}}

\definecolor{shadecolor}{gray}{0.875}

\newtheorem{thrm}{Theorem}[section]
\newtheorem{thrmx}{Theorem}
\newtheorem{conjx}{Conjecture}

\newtheorem{lem}[thrm]{Lemma}
\newtheorem{cor}[thrm]{Corollary}
\newtheorem{prop}[thrm]{Proposition}
\newtheorem{conj}[thrm]{Conjecture}

\theoremstyle{definition}
\newtheorem{defn}[thrm]{Definition}

\newtheorem{rmk}[thrm]{Remark}

\newtheorem{ques}[thrm]{Question}

\DeclareMathOperator{\Chow}{Chow}
\DeclareMathOperator{\Ann}{Ann}

\DeclareMathOperator{\cone}{cone}

\DeclareMathOperator{\Amp}{Amp}
\DeclareMathOperator{\Nef}{Nef}
\DeclareMathOperator{\Eff}{Eff}
\DeclareMathOperator{\Psef}{Psef}

\DeclareMathOperator{\vol}{vol}
\DeclareMathOperator{\Mov}{Mov}

\DeclareMathOperator{\conv}{Conv}

\DeclareMathOperator{\Image}{Image}

\DeclareMathOperator{\nd}{nd}

\DeclareMathOperator{\codim}{codim}

\DeclareMathOperator{\ord}{ord}
\DeclareMathOperator{\bg}{Big}
\DeclareMathOperator{\sing}{sing}

\title{Positivity in the shadow of Hodge index theorem}
\author{Jiajun Hu and Jian Xiao}
\date{}

\begin{document}
\maketitle

%\begin{dedication}
%In memory of ...
%\end{dedication}

\begin{abstract}
    Taking a compact K\"ahler manifold as playground, we explore the powerfulness of Hodge index theorem. A main object is the Lorentzian classes on a compact K\"ahler manifold, behind which the characterization via Lorentzian polynomials over the K\"ahler cone and hence the validity of Hodge index theorem. Along the exploration, we discover several applications in complex geometry that may be unexpected before. 
    (1) For a Lefschetz type operator given by the complete intersection of nef classes, we give a complete characterization of its kernel face against the pseudo-effective cone.    
    (2) We provide a new approach to Teissier's proportionality problem from the validity of hard Lefschetz property. This perspective enables us to establish the extremals for the Brunn-Minkowski inequality on a strictly Lorentzian class, and thus also characterize the most extremal case for a log-concavity sequence given by the intersection numbers of two nef classes. These Lorentzian classes include the fundamental classes of smooth projective varieties or compact K\"ahler manifolds as typical examples, hence our result extends Boucksom-Favre-Jonsson's and Fu-Xiao's results in respective settings to broader contexts, e.g. certain algebraic cycle classes given by reducible subvarieties.   
    (3) Furthermore, we also strengthen the proportionality characterization by comparing various quantitative deficits and establishing stability estimates.
    Two quantitative sharper stability estimates with close relation with complex Monge-Amp\`ere equations and Newton-Okounkov bodies are also discussed.    
\end{abstract}

\tableofcontents

\section{Introduction}

\subsection{Hodge index theorem}
Let $X$ be a compact K\"ahler manifold of dimension $n$, and denote by $\Amp^1(X)\subset H^{1,1}(X,\mathbb{R})$ the open convex cone of all K\"ahler classes.
The starting point of this article is the Hodge index theorem:
\begin{quote}
(\textbf{HIT})
    \emph{For any K\"ahler class $\omega\in\Amp^1(X)$, the bilinear form
    $$H^{1,1}(X,\mathbb{R})\times H^{1,1}(X,\mathbb{R}) \rightarrow \mathbb{R},  \ (\alpha,\beta)\mapsto \alpha\cdot \beta \cdot \omega^{n-2}$$
    has Lorentzian signature $(+,-,...,-)$.}
\end{quote}
As first noted by Gromov \cite{gromov1990convex} in the setting of K\"ahler geometry (see also \cite{cattanimixedHRR, DN06}), the same conclusion remains valid when the class $\omega^{n-2}$ is replaced a product of different K\"ahler classes.
We refer to this generalization the mixed Hodge index theorem:
\begin{quote}
(\textbf{MHIT})
\emph{For any K\"ahler classes $\omega_1,...,\omega_{n-2}\in \Amp^1(X)$, the bilinear form
$$H^{1,1}(X,\mathbb{R})\times H^{1,1}(X,\mathbb{R}) \rightarrow \mathbb{R},  \ (\alpha,\beta)\mapsto \alpha\cdot \beta \cdot \omega_1\cdots \omega_{n-2}$$
    has Lorentzian signature $(+,-,...,-)$.}
\end{quote}
From $(\textbf{MHIT})$ one can derive numerous inequalities originally due to Khovanskii and Teissier, analogous to the Alexandrov-Fenchel inequalities in convex geometry. If we focus on the inequalities, we can work within the closure of $\Amp^1(X)$, known as the nef cone and denoted by $\Nef^1(X)$.
For a nef class $\alpha \in \Nef^1(X)$, its volume is defined as $\vol(\alpha):=\alpha^n$, also denoted by $ |\alpha|$.
Then for any nef classes $\alpha,\beta,\omega_1,...,\omega_{n-2} \in \Nef^1(X)$, we have the Alexandrov-Fenchel (or Khovanskii-Teissier) inequality:
\begin{equation}\tag{\textbf{AF}}
    (\alpha\cdot \beta \cdot \omega_1\cdot...\cdot\omega_{n-2})^2 \geq (\alpha^2\cdot \omega_1\cdot...\cdot\omega_{n-2})(\beta^2\cdot \omega_1\cdot...\cdot\omega_{n-2}),
\end{equation}
and as a consequence, the Brunn-Minkowski inequality
\begin{equation}\tag{\textbf{BM}}
    |\alpha+\beta|^{\frac{1}{n}}\geq |\alpha|^{\frac{1}{n}}+|\beta|^{\frac{1}{n}}.
\end{equation}
Another consequence of (\textbf{AF}) is the log-concavity of the sequence of the intersection numbers given by two nef classes. Beyond geometry, motivated by different problems in other areas, the same algebraic structure behind (\textbf{HIT}) or (\textbf{MHIT}) has been established in much broader contexts. For example, the Hodge index theorem has its counterparts in combinatorial Hodge theory or in a more basic framework of Lorentzian polynomials (see e.g. \cite{huhHRR, branhuhlorentz, brandenleake23coneLorentzian, AnariMasonconj}), which have found profound applications -- resolving the Heron-Rota-Welsh conjecture and {M}ason's ultra-log-concavity conjecture.  Inside geometry itself, in our previous work \cite{hxhardlef} we prove the Hall-Rado relation in K\"ahler geometry via Hodge index theorem (see also \cite{HuXiaoLefChar} for its extension to Lorentzian polynomials), and in \cite{hxinterineq} we applied the perspective of Lorentzian polynomials to study intersection numbers on a compact K\"ahler manifold and obtained several unexpected results, including rKT property, submodularity of the numerical dimension (or equivalently, polymatroid structure on nef classes), see Sections \ref{sec collect nef}, \ref{sec rkt} for some discussions.

\begin{comment}
Although not that obvious, $(\textbf{MHIT})$ also implies the following reverse Khovanskii-Teissier inequality as noted by \cite{hxinterineq}
\begin{equation}\tag{\textbf{rKT}}
    (\alpha_1\cdots\alpha_k\cdot \omega^{n-k})(\beta^{k}\cdot \omega^{n-k})\leq 2^{l(k-l)}(\alpha_1\cdots\alpha_l\cdot \beta^{k-l}\cdot \omega^{n-k})(\alpha_{l+1}\cdots\alpha_{k}\cdot \beta^{l}\cdot\omega^{n-k}).
\end{equation}
Using the (\textbf{rKT}), the authors of \cite{hxinterineq} also proved
\begin{equation}\tag{\textbf{Submodularity}}
    \nd(\alpha+\beta+\gamma)+\nd(\gamma)\leq \nd(\alpha+\gamma)+\nd(\beta+\gamma).
\end{equation}
\end{comment}

\subsection{Motivation}
Taking a compact K\"ahler manifold as playground, we aim to explore further powerfulness of the Hodge index theorem.

One motivation is the extremals in $(\textbf{AF})$ and $(\textbf{BM})$, that is, 
\begin{quote}
\begin{center}
  \emph{When does the equality hold in (\textbf{AF}) and (\textbf{BM})?}
\end{center}
\end{quote}
When $\omega_1,...,\omega_{n-2}$ are all K\"ahler, equality in $(\textbf{AF})$ holds if and only if $\alpha$ and $\beta$ are proportional by (\textbf{MHIT}). However, when the classes $\omega_1,...,\omega_{n-2}$ are only nef, characterizing the equality case of (\textbf{AF}) becomes extremely challenging and is still open in full generality. In different or specific contexts, the extremal problem for (\textbf{AF}) has been resolved, such as compact complex tori \cite{Panov1987ONSP}, convex polytopes \cite{SvH23Acta} (hence smooth projective toric varieties). In \cite{HuXiaoLefChar}, we solved the extremal problem under certain positivity assumption -- ``supercriticality with a rearrangement'' -- on the collection $\{\omega_1,...,\omega_{n-2}\}$, in particular, this includes the case of big nef collections (i.e., all $|\omega_i|>0$). For the general case, see \cite{HuXiaoLefChar} for the conjectural picture in K\"ahler geometry.
The problem of characterizing the equality in (\textbf{BM}) is originally proposed by Teissier \cite{teissier82} who also conjectured that
for $\alpha,\beta$ big and nef, the equality 
$$|\alpha+\beta|^{\frac{1}{n}}= |\alpha|^{\frac{1}{n}}+|\beta|^{\frac{1}{n}}$$
holds if and only if $\alpha$ and $\beta$ are proportional. This conjecture, known as Teissier's proportionality problem, and has been solved in various settings \cite{BFJ09, fx19, lehmxiaoPosCurveANT} (see also \cite{cutkosky13, Cutkosky24} for extensions to a general field). However, on a cycle class that is given be reducible subvarieties, one may still have (\textbf{BM}) for the volume on the cycle class (see Section \ref{sec exmple lorent}). The extremal problem for this kind of (\textbf{BM}) corresponds to the characterization of the most extremal case (i.e., flat case) of a log-concave sequence, which is given by the intersection of two nef classes on a cycle class. We shall solve the extremal problem for (\textbf{BM}) over a strictly Lorentzian class and also establish quantitative stability estimates.
The same question also arises for  (\textbf{AF}).
As the current tools seems unapplicable, less is known in this setting.

Another basic question is the kernel of a Lefschetz type operator. More precisely, fix a collection of nef classes $\mathcal{L}=\{L_1, ...,L_m\}$, and denote the complete intersection class by $\mathbb{L}=L_1 \cdot ...\cdot L_m$, which can be viewed as a linear operator acting on some cohomology group or cycle classes group. A natural question is, 
\begin{quote}
\begin{center}
  \emph{What is the kernel space of $\mathbb{L}$?}
\end{center}
\end{quote}
A complete resolution to this question will be very challenging, as even in some special cases it could imply a strong form of the generalized Hodge conjecture for smooth projective varieties \cite{voisinPushforw}.

Nevertheless, the problem becomes accessible if we know well the space that the operator $\mathbb{L}$ acts on. Even more, we can first consider the kernel space $\ker \mathbb{L}$ against some convex cone given by certain positivity.
It is closely related to the hard
Lefschetz property of $\mathbb{L}$ and plays an important role in the algebro-geometric formalism of the extremals in log-concavity sequences \cite{HuXiaoLefChar}. 
In this paper, we study its action on the closed convex cone consisting of the classes of all $d$-closed positive $(1,1)$-currents, known as the pseudo-effective (psef) cone $\Psef^1(X)\subset H^{1,1}(X,\mathbb{R})$. It is also closely related the following extremal problem of intersection numbers.
If $\beta-\alpha \in \Psef^1(X)$ is psef, then
\begin{equation*}
    \mathbb{L}\cdot \alpha \leq \mathbb{L}\cdot \beta,
\end{equation*}
where $\leq$ means that $\mathbb{L}\cdot \beta - \mathbb{L}\cdot \alpha$ is a psef class in $H^{m+1, m+1} (X, \mathbb{R})$.
To understand when equality holds in the above monotonicity inequality, it is equivalent to solve the following problem:

\begin{quote}
(\textbf{KF})
\emph{Characterize the kernel face $\Psef^1 (X) \cap \ker \mathbb{L} = \{\alpha\in \Psef^1(X):\mathbb{L}\cdot \alpha =0\}$.}
\end{quote}

It is easy to see that for any K\"ahler class $\omega$, $$\alpha\in \Psef^1 (X) \cap \ker \mathbb{L} \Leftrightarrow \alpha\in \Psef^1 (X) \cap \ker (\mathbb{L} \cdot \omega^{n-m-1}).$$
Therefore, in the study of (\textbf{KF}) problem, we may only need to consider the case when $\mathbb{L}$ is given by a product of $(n-1)$ nef classes. We will give a complete answer to this problem.

\subsection{The main results}

\subsubsection{The kernel face against the pseudo-effective cone}

Our first main result solves the (\textbf{KF}) problem completely. To move on, we introduce the notion of numerical dimension for a collection of nef classes $\mathcal{L}=\{L_1,...,L_m\}$:
    $$\nd(\mathcal{L}):=\min_{\emptyset \neq I\subset [m]}\{\nd(L_I)-|I|+m\},$$
where $L_I = \sum_{i\in I} L_i$ and $\nd(L_I)$ is the usual numerical dimension for a single class. For more properties on $\nd(\mathcal{L})$, see Section \ref{sec collect nef}.

\begin{thrmx}\label{intro thrm kerf}
    Let $X$ be a compact K\"ahler manifold of dimension $n$ and let $\mathcal{L}=\{L_1,...,L_{n-1}\} $ be a collection of nef classes. 
     Then we have:
\begin{enumerate}
    \item If $\nd(\mathcal{L})\leq n-2$, then $\mathbb{L}=0$ and thus $\ker \mathbb{L}\cap \Psef^1(X) =\Psef^1 (X)$.
    \item If $\nd(\mathcal{L})=n-1$, then $\ker \mathbb{L}\cap \Psef^1(X)$ is generated by the classes of hypersurfaces $D$ such that $$\nd_D (L_{I}) < |I|$$ 
        for some $I\subset [n-1]$ and the movable classes $M$ such that 
        $$\nd(L_I+M)=\nd(L_I)$$ 
        for some $I\subset [n-1]$ with $\nd(L_I)=|I|$. Here, $\nd_D (-)$ is the numerical dimension on $D$.
    %$$\ker \mathbb{L}\cap \Mov^1(X)=\sum_{I\subset [n-1]:\nd(L_I)=|I|}\{M \in \Mov^1(X):\nd(L_I+M)=\nd(L_I)\}$$.
    \item If $\nd(\mathcal{L})=n$, then $\ker \mathbb{L}\cap \Psef^1(X)$ is generated by the classes of hypersurfaces $D$ such that $$\nd_D (L_{I}) < |I|$$ 
        for some $I\subset [n-1]$. That is, the movable part in (2) vanishes and only the divisorial part remains.
\end{enumerate}
\end{thrmx}

From the proof, we shall see that the hypersurfaces $D$ have some special properties -- they must be exceptional and have strong rigidity.

\subsubsection{Extremals for (\textbf{BM}) on a strictly Lorentzian class}
One way to introduce Lorentzian classes is using Lorentzian polynomials -- this is one main reason that we give this name. For other equivalent reformulations, see Section \ref{sec Lorenclass}.
Given $\Omega \in H^{n-d, n-d}(X, \mathbb{R})$, we can define the volume polynomial on $\Omega$ as 
$$f_{\Omega}: H^{1,1}(X, \mathbb{R})\rightarrow \mathbb{R}, \ \alpha \mapsto \alpha^d\cdot \Omega.$$

\begin{defn}
The class $\Omega$ is (strictly) Lorentzian iff the polynomial $f_\Omega$ is (strictly) $\Amp^1(X)$-Lorentzian.
\end{defn}

Fix a Lorentzian class $\Omega\in H^{n-d,n-d}(X,\mathbb{R})$ (we call it is of dimension $d$). 
For any nef class $\alpha$, we denote its volume on $\Omega$ by $$|\alpha|_{\Omega}=\alpha^{d}\cdot\Omega$$ and we call that $\alpha$ is big on $\Omega$ if $|\alpha|_{\Omega}>0$. 
Then, by the Lorentzian property of $f_\Omega$, the Brunn-Minkowski inequality holds on $\Omega$: for any nef classes $\alpha $ and $ \beta$,
\begin{equation}\tag{$\textbf{BM}_{\Omega}$}
    |\alpha+\beta|_{\Omega}^{\frac{1}{d}}\geq |\alpha|_{\Omega}^{\frac{1}{d}}+|\beta|_{\Omega}^{\frac{1}{d}}.
\end{equation}

Our second main result solves the extremal problem for ($\textbf{BM}_{\Omega}$) on a strictly Lorentzian class, which include certain cycle classes given by reducible subvarieties.  In this relative setting the previous methods seems to not apply. Therefore our proof provides a novel approach\footnote{While we independently find this approach, when we were finishing the manuscript and checking the references, we realized that a prototype in the convexity analog had already appeared in \cite[Theorem 7.6.9]{schneiderBrunnMbook} where $\Omega$ reads as a complete intersection of K\"ahler classes, which is strictly Lorentzian.
The proof utilizes an inductive procedure established by Fenchel \cite{Fenchel1936}, along with the characterizations of the equality cases of the Alexandrov-Fenchel inequalities for smooth bodies, as developed by Schneider \cite{Schneider90AF}. Our approach applies an inductive procedure which is essentially the same as that in the convexity setting.} even in the absolute case, i.e., when $\Omega =[X]$.

\begin{thrmx}\label{thrmx proportionality}
Let $X$ be a compact K\"ahler manifold of dimension $n$ and let $\Omega$ be a strictly Lorentzian class of dimension $d$ ($2\leq d \leq n-2$), and let $\alpha,\beta$ be big nef classes on $\Omega$. Then the followings are equivalent:
    \begin{enumerate}
        \item $\alpha,\beta$ are proportional;
        \item For some $1\leq k\leq  d-1$, the classes $\alpha^k$ and $\beta^k$ are proportional;
        \item For some $1\leq k\leq d-1$, the classes $\alpha^{k}\cdot \Omega$ and $\beta^{k}\cdot \Omega$ are proportional;
        \item For any $1\leq k\leq d-1$,
         $$(\alpha^k \cdot \beta^{d-k} \cdot \Omega)^2 = (\alpha^{k-1} \cdot \beta^{d-k+1} \cdot \Omega)(\alpha^{k+1} \cdot \beta^{d-k-1} \cdot \Omega);$$
        \item The equality holds in ($\textbf{BM}_{\Omega}$):
        $$|\alpha+\beta|_{\Omega}^{\frac{1}{d}}= |\alpha|_{\Omega}^{\frac{1}{d}}+|\beta|_{\Omega}^{\frac{1}{d}}. $$
    \end{enumerate}
\end{thrmx}

On the necessity of $\alpha,\beta$ being big nef classes on $\Omega$ to conclude proportionality, see the discussion at the beginning of Section \ref{sec teisprop} and Proposition \ref{BMeqnotbig}.

The third statement just says that the log-concave sequence $\{(\alpha^k \cdot \beta^{d-k} \cdot \Omega)\}_{k=0} ^d$ lies in a line -- the most extremal case for a log-concave sequence. For more equivalent descriptions on the extremals for ($\textbf{BM}_{\Omega}$), see Lemma \ref{logconc} and Theorem \ref{TeissierLorentzian}. An interesting consequence is the injectivity property on taking powers of big nef classes on a strictly Lorentzian class.

\subsubsection{Quantitative stability estimates}
Our third main result is the quantitative refinement of Theorem \ref{thrmx proportionality}. Roughly speaking, it is about the question: how to measure the difference of the rays given by $\alpha,\beta$  by using the difference of volumes
$$|\alpha+\beta|_{\Omega}^{\frac{1}{d}} - \left(|\alpha|_{\Omega}^{\frac{1}{d}}+|\beta|_{\Omega}^{\frac{1}{d}}\right)? $$ 

To analyze the stability problem, we introduce several deficit functions.

\begin{defn}\label{defndeficit}
Let $\alpha,\beta$ be big and nef classes on a Lorentzian class $\Omega$ of dimension $d$.
\begin{itemize}
\item Fix a K\"ahler class $\omega$, the Alexandrov-Fenchel (AF) deficit is given by
    $$A(\alpha,\beta):=A(\alpha,\beta;\omega^{d-2}\cdot \Omega)=\frac{\sqrt{(\alpha\cdot\beta\cdot\omega^{d-2}\cdot \Omega)^2-(\alpha^2\cdot\omega^{d-2}\cdot \Omega)(\beta^2\cdot\omega^{d-2}\cdot \Omega)}}{(\alpha\cdot \beta\cdot \omega^{d-2}\cdot \Omega)}.$$
    \item The Brunn-Minkowski (BM) deficit is given by
    $$B(\alpha,\beta)
    =\frac{|\alpha+\beta|_{\Omega}^{\frac{1}{d}}}
    {|\alpha|_{\Omega}^{\frac{1}{d}}+|\beta|_{\Omega}^{\frac{1}{d}}}-1.$$
    \item The Khovanskii-Teissier (KT) deficit is defined by
    $$K(\alpha,\beta)=\frac{\alpha^{d-1}\cdot\beta\cdot \Omega}{|\alpha|_{\Omega}^{\frac{d-1}{d}}|\beta|_{\Omega}^{\frac{1}{d}}}-1.$$
    \item The relative size index is given by  
    $$\sigma(\alpha,\beta)=\max\left\{\left(\frac{|\alpha|_{\Omega}}{|\beta|_{\Omega}}\right)^{\frac{1}{d}},
    \left(\frac{|\beta|_{\Omega}}{|\alpha|_{\Omega}}\right)^{\frac{1}{d}}\right\}$$
\end{itemize}
\end{defn}

When $\Omega$ is strictly Lorentzian, $A(\alpha,\beta) =0$ iff $\alpha,\beta$ are proportional. Note that the implication $(4)\Rightarrow (1)$ in Theorem \ref{thrmx proportionality} is exactly the implication $$B(\alpha,\beta)=0 \Rightarrow A(\alpha,\beta) =0.$$ 
We quantify this implication as follows.

\begin{thrmx}\label{thrmx deficitcomp}
Let $\Omega$ be a Lorentzian class of dimension $d$,
and let $\alpha,\beta$ be big nef classes on $\Omega$. 
Then there exists a constant $c(d)$ depending only on $d$ such that 
    \begin{equation*}
        A(\alpha,\beta)\leq c(d) \Big(\sigma(\alpha,\beta)B(\alpha,\beta)\Big)^{\frac{1}{2^{d-2}}}.
    \end{equation*}
    
This indeed follows from the following more refined estimates: 
\begin{enumerate}
  \item There exists a constant $c(d)$ depending only on $d$ such that 
    \begin{equation*}
        A(\alpha,\beta)\leq c(d) K(\alpha,\beta)^{\frac{1}{2^{d-2}}}.
    \end{equation*}
  \item Assume $|\alpha|_{\Omega} \geq |\beta|_{\Omega}$, then there exists a constant $c(d)$ depending only on $d$ such that: 
$$\left(\frac{|\beta|_{\Omega}}{|\alpha|_{\Omega}}\right)^{1/d}K(\alpha,\beta) \geq B(\alpha, \beta)\geq c(d) \frac{ \left(\frac{|\beta|_{\Omega}}{|\alpha|_{\Omega}}\right)^{1/d}K(\alpha,\beta) }{1+ \left(\frac{|\beta|_{\Omega}}{|\alpha|_{\Omega}}\right)^{1/d}K(\alpha,\beta)}
.$$
\end{enumerate}

\end{thrmx}

For further results on various deficits, in particular the relation with the relative radii of $\alpha, \beta$, and the domination of BM/KT deficit by AF deficit (which involves the radii), see Section \ref{sec comp defi}.

\subsubsection{FMP type estimates}
Continuing in the framework of quantitative stability estimates, our last main result is motivated 
by the breakthrough works due to Fusco-Maggi-Pratelli and Figalli-Maggi-Pratelli \cite{FMP08Annals, figalli09, Figalli10inventiones}. Let us briefly describe their result, which we call FMP estimates for short. Let $M,N\subset \mathbb{R}^n$ be convex bodies. The relative asymmetry index (differs from the usual one by a factor $\frac{1}{2}$) is given by
\begin{equation*}
    F(M,N)=\inf_{x\in\mathbb{R}}\frac{|(x+rM)\Delta N|}{2|N|}
\end{equation*}
where $|\cdot|$ denotes the Lesbesgue measure, $r=\left( \frac{|N|}{|M|} \right)^{\frac{1}{n}}$, and $$(x+rM)\Delta N=(x+rM)\backslash N \cup N\backslash (r+rM)$$ 
is the symmetric difference.
FMP estimate tells that there exists a constant $c(n)$ depending only on $n$ such that
\begin{equation}\tag{\textbf{FMP}}
    F(M,N) \leq c(n)(\sigma(M,N)B(M,N))^{\frac{1}{2}}
\end{equation}
where $\sigma (M, N)$ and $B(M, N)$ are the relative size index and BM deficit respectively. Furthermore, $c(n)$ has polynomial growth in $n$. The crucial point in (\textbf{FMP}) is that the exponent $1/2$ is dimension free and it is optimal.

Let us go back to the K\"ahler geometry world.
First, it is straightforward to check that $F(M,N)$ admits an equivalent formulation
\begin{equation*}
    F(M,N)=1-\sup_{K\leq rM,N} \frac{|K|}{|N|} 
\end{equation*}
where $K\leq rM,N$ means that $K$ ranges over all convex bodies contained in both $rM$ and $N$ up to translations.
This observation suggests a natural analogue for positive classes.
\begin{defn}
    The relative asymmetry index of two movable classes $\alpha,\beta\in \Mov^1(X)$ with $r=\left( \frac{|\beta|}{|\alpha|} \right)^{\frac{1}{n}}$ is defined by
    \begin{equation*}
        F(\alpha,\beta)=1-\sup_{\gamma\leq r\alpha,\beta} \frac{|\gamma|}{|\beta|},
    \end{equation*}
    where $\gamma \in \Mov^1(X)$ ranges over all movable classes satisfying that both $\beta-\gamma$ and $r\alpha-\gamma$ are psef.
\end{defn}

We then formulate the following stability conjecture in the absolute case.
\begin{conjx}
Let $X$ be a compact K\"ahler manifold of dimension $n$, and let $\alpha,\beta\in \Mov^1(X)$ be two big movable classes. Then there exists a constant $c(n)$ depending only on the dimension $n$ such that
    \begin{equation*}
        F(\alpha,\beta)\leq c(n)(\sigma(\alpha,\beta)B(\alpha,\beta))^{\frac{1}{2}}.
    \end{equation*}
\end{conjx}
It is sure that one can formulate the conjecture for a general (strictly) Lorentzian class, see Conjecture \ref{StabofBMOmega}.
While we do not have much evidence for the $\frac{1}{2}$-exponent in this general setting, by relating the asymmetry index with radii, we obtain a related $\frac{1}{2^{d-2}}$-exponent estimate on an arbitrary Lorentzian class of dimension $d$ (see Proposition \ref{fomegaKT}).

Over $X$, we can prove related 
$\frac{1}{2}$-exponent estimates by replacing $F$ with other asymmetry measures, that also capture the failure of proportionality.

The first one is a metric-norm asymmetry for big nef classes defined via suitable representatives in the classes on a compact K\"ahler manifold, using Calabi-Yau theorem -- complex Monge-Amp\`ere equations.
Let $X$ be a compact K\"ahler manifold of dimension $n$ and let $\alpha,\beta$ be big nef classes on $X$. Define
    \begin{equation*}
        \widehat{F}(\alpha,\beta)=\sup_{\langle\widehat{\alpha}^n \rangle=\frac{|\alpha|}{|\beta|}\langle\widehat{\beta}^n\rangle = \Phi}  \frac{1}{|\alpha|} \int_{\Amp(\alpha,\beta)} \left(\frac{\Vert \widehat{\alpha}-(\frac{|\alpha|}{|\beta|})^{1/n}
        \widehat{\beta}\Vert_{\widehat{\beta}}}{\Vert\widehat{\alpha}\Vert_{\widehat{\beta}}}\right)
        \Phi,
    \end{equation*}
    where 
    \begin{itemize}
      \item $\Amp(\alpha,\beta) = \Amp(\alpha)\cap \Amp(\beta)$ is the intersection of the ample locus of $\alpha, \beta$, which is a Zariski open subset in $X$.
      \item $\Phi$ ranges over all the smooth volume forms such that $\int_X \Phi = |\alpha|$.
      \item $\widehat{\alpha}\in\alpha$ and $\widehat{\beta}\in \beta$ range over all the pairs of positive currents such that $$\langle\widehat{\alpha}^n\rangle=\frac{|\alpha|}{|\beta|}\langle\widehat{\beta}^n\rangle =\Phi,$$ 
    whose existence is guaranteed by \cite{BEGZ10}. Here, $\langle -\rangle$ is the non-pluripolar product of positive currents.
    \item $\Vert -\Vert_{\widehat{\beta}}$ is the pointwise $L^2$-norm induced by $\widehat{\beta}$, which is a K\"ahler metric on the ample locus. 
    \end{itemize}
When $\alpha, \beta$ are K\"ahler, $\Amp(\alpha,\beta) =X$ and the existence of $\widehat{F}(\alpha,\beta)$ is then guaranteed by Yau's theorem \cite{Yau78}.
By \cite{fx19}, it is a nontrivial result that $\widehat{F}(\alpha,\beta)=0$ iff $\alpha,\beta$ are proportional.

The second one is an asymmetry measure for big movable divisor classes on a projective variety, using the theory of Newton-Okounkov bodies. Assume now that $X$ is a smooth projective variety of dimension $n$ over an algebraically closed field $k$ (one may just think $k$ as $\mathbb{C}$). 
For each valuation $\nu:k(X)^{\times}\rightarrow \mathbb{Z}^n$ with maximal rank, let $\Delta_{\nu}(\alpha)$ and $\Delta_{\nu}(\beta)$ be the Newton-Okounkov bodies of big classes $\alpha,\beta \in N^1(X)$.
Define 
    \begin{equation*}
        \widetilde{F}(\alpha,\beta)=\sup_\nu F(\Delta_{\nu}(\alpha),\Delta_{\nu}(\beta)).
    \end{equation*}
By the arguments of \cite{Jow10} (see also Lemma \ref{usingjow}), one checks that if $\alpha,\beta$ are big and movable, then $\widetilde{F}(\alpha,\beta)=0$ exactly when $\alpha,\beta$  are proportional

We prove the following FMP type estimates:
\begin{thrmx}\label{thmrx FMP}
The FMP type estimates hold in the following settings:
\begin{enumerate}
  \item Let $X$ be a compact K\"ahler manifold of dimension $n$ and let $\alpha,\beta$ be big nef classes on $X$.
Then, there is a constant $c(n)$ depending only on the dimension $n$ such that
    \begin{equation*}
        \widehat{F}(\alpha,\beta)\leq c(n)\sqrt{\sigma(\alpha,\beta)B(\alpha,\beta)}.
    \end{equation*}
  \item Let $X$ be a smooth projective variety of dimension $n$ and let $\alpha,\beta$ be big movable classes on $X$.
Then, there is a constant $c(n)$ depending only on the dimension $n$ such that
    \begin{equation*}
        \widetilde{F}(\alpha,\beta)\leq c(n)\sqrt{\sigma(\alpha,\beta)B(\alpha,\beta)}.
    \end{equation*}
\end{enumerate}
\end{thrmx}

The first part of Theorem \ref{thmrx FMP} provides a quantitative refinement of proportionality results obtained in \cite[Theorem D]{BFJ09} and \cite{fx19}, and the second part quantitatively refines some results in \cite[Section 3]{lehmxiaoPosCurveANT}.

\subsection{Organization}
This paper is organized as follows. In Section \ref{sec pre}, we introduce some basic notions and preliminary results. In Section \ref{sec Lorenclass}, we study the basic properties of Lorentzian classes and present some examples. Section \ref{sec nullface} devotes to the characterization of the kernel face on the pseudo-effective cone, including the proof of Theorem \ref{intro thrm kerf}. In Section \ref{sec propor}, we give the proof of Theorem \ref{thrmx proportionality}. Finally, in Section \ref{sec stabestim}, we strengthen the proportionality characterization by establishing quantitative stability estimates, including the proof of Theorem \ref{thrmx deficitcomp} and Theorem \ref{thmrx FMP}.

\section{Preliminaries}\label{sec pre}

In this section, we introduce some basic notions and results which will be applied in the sequel. Unless otherwise stated, throughout this section, $X$ is a compact K\"ahler manifold of dimension $n$.

\subsection{Positivity}
We first briefly introduce some positivity notions in K\"ahler geometry. The standard references are \cite{dem_analyticAG, Dem_AGbook}. For the counterparts in algebraic geometry, see \cite{lazarsfeldPosI, lazarsfeldPosII}.

Over $X$, $H^{p,p}(X, \mathbb{R})$ is the real vector space of cohomology classes of bidegree $(p,p)$. We also call the element in $H^{p,p}(X, \mathbb{R})$ has dimension $n-p$.
As we only study real classes, we sometimes just denote $H^{p,p}(X, \mathbb{R})$ by $H^{p,p}(X)$. When $X$ is projective, we denote by $N^p(X)$ the vector space of numerically equivalent cycle classes of codimension $p$ with real coefficients.

\begin{defn}
We have the following positivity notions.
\begin{itemize}
\item $\Psef^1(X) \subset H^{1,1}(X)$ is the closed convex cone of pseudo-effective (psef for short) $(1,1)$ classes, that is, $L \in \Psef^1(X)$ if and only if it can be represented by a positive $(1,1)$ current. Its interior is an open convex cone given by big $(1,1)$ class, which we denote by $\bg^1(X)$. 
    
    More generally, $\Psef^p(X) \subset H^{p,p}(X)$ is the closed convex cone generated by $d$-closed positive $(p,p)$ currents. 

\item $\Nef^1(X) \subset H^{1,1}(X) $ is the closed convex cone of nef $(1,1)$ classes. Its interior is called the K\"ahler cone, denoted by $\Amp^1 (X)$, which consists of K\"ahler classes.
\item $\Mov^1(X)\subset H^{1,1}(X)$ is the closed convex cone of movable $(1,1)$ classes, that is, $L\in \Mov^1(X)$ if and only if $$L = \lim_{k\rightarrow \infty} (\pi_k)_* \omega_k,$$ where $\pi_k: X_k \rightarrow X$ is a K\"ahler modification and $\omega_k$ is a K\"ahler class on $X_k$.
\end{itemize}
\end{defn}

We shall use the notation $\geq$ (respectively, $\leq$) to denote the partial order induced by the psef cone, that is, $\alpha \geq \beta $ (respectively, $\alpha \leq \beta $)  iff $ \alpha - \beta$  (respectively, $ \beta -\alpha $) is psef.

For a big class $L\in \bg^1 (X)$, its ample locus is the (analytic) Zariski open set 
\begin{equation*}
  \Amp(L) = X \setminus \bigcap_T \sing T
\end{equation*}
where $T$ ranges over all the K\"ahler currents in the class $L$ with analytic singularities and $\sing T$ is the singular locus of $T$. The Zariski closed subset $\bigcap_T \sing T$ is called the non-ample (or non-K\"ahler) locus, and in the algebraic case when $L$ is a big divisor class, it coincides with the augmented base locus
\begin{equation*}
\mathbb{B}_+ (L) =  \bigcap_{L\sim_{\mathbb{R}} A +E} \text{supp(E)},
\end{equation*}
where $L\sim_{\mathbb{R}} A +E$ is the Kodaira decomposition of $L$ with $A$ ample, $E$ effective.

It is clear that 
\begin{equation*}
  \Nef^1(X) \subset \Mov^1(X) \subset \Psef^1(X), 
\end{equation*}
and in general the inclusions are strict. A remarkable result relating the cones $\Psef^1(X), \Mov^1(X)$ is the divisorial Zariski decomposition due to Boucksom \cite{Bou04} (see also Nakayama \cite{Nak04} for the algebraic construction).

\begin{thrm}[Boucksom \cite{Bou04}]
Every $L\in \Eff^1(X)$ can be decomposed as follows:
\begin{equation*}
  L=P(L)+[N(L)],
\end{equation*}
where $P(L) \in \Mov^1(X)$ is called the positive part, and $N(L)$ is an effective divisor which is called the negative part. 
 Furthermore, $N(L)$ is supported by at most $\rho(X)$ prime divisors, where $\rho(X)$ is the Picard number of $X$.
\end{thrm}

This result plays an important role in Section \ref{sec nullface}. The divisor $N(L)$ has very strong rigidity in the sense that there is only one positive current in the class.

\subsection{Positive product}
The positive product is an extremely powerful tool in the study of positivity of $(1,1)$ classes.  This theory is systematically developed in \cite{BFJ09,BDPP13,BEGZ10} from diverse aspects.
We briefly recall its construction and some basic properties.

Given big classes $\alpha_1,...,\alpha_p \in H^{1,1}(X)$, the set 
$$\{\pi_*(\beta_1\cdot...\cdot\beta_p)\in \Psef^p(X) | \pi:X_{\pi}\rightarrow X \text{ K\"ahler modifications}, \beta_i \leq \pi^*\alpha \text{ K\"ahler classes}\}$$
admits a unique supremum under the partial order induced by the cone $\Psef^p (X)$. The supremum, denoted by $$\langle\alpha_1 \cdot...\cdot \alpha_p\rangle,$$ is called the positive product of $\alpha_1,...,\alpha_p$. 
Once the existence is known, we can always find a sequence of K\"ahler modifications $\pi_j:X_j\rightarrow X$ and K\"ahler classes $\beta_i^{j}\leq \pi_j^*\alpha_i$ on $X_j$ such that
$$\langle\alpha_1 \cdot...\cdot  \alpha_p\rangle=\lim{\pi_j}_*(\beta_1^j\cdot...\cdot\beta_p^j).$$

The positive product is clearly increasing in each variable. Hence we can define the positive product of psef classes $\alpha_1,...,\alpha_p \in \Psef^1(X)$ to be $$\langle\alpha_1 \cdot...\cdot \alpha_p\rangle=\lim_{\epsilon\rightarrow0^+ }\langle(\alpha_1+\epsilon\omega)\cdot...\cdot (\alpha_p+\epsilon \omega) \rangle$$ where $\omega $ is an arbitrary fixed K\"ahler class. 

\begin{prop}
The following are some basic properties of positive products:
\begin{enumerate}
    \item The positive product 
    \[
    \underbrace{\Psef^1(X)\times \cdots \times \Psef^1(X)}_{p\text{-fold}}
    \;\longrightarrow\;
    \Psef^p(X)
    \] 
    is an upper-semicontinuous symmetric function which is continuous in $\bg^1(X)\times \cdots \times \bg^1(X)$. Moreover, it is 
    1-homogeneous, increasing, superadditive in each variable.
    
    \item  The positive product depends only on the positive part of each variable.

    \item For any psef class $\alpha \in \Psef^1(X)$, $\langle\alpha\rangle$ lies in $\Mov^1(X)$ and coincides with the positive part of its divisorial Zariski decomposition, that is,
        \begin{equation*}
          \langle\alpha\rangle =P(\alpha).
        \end{equation*}
       Hence, we have $\langle\alpha \rangle =\alpha $ if $\alpha \in \Mov^1(X)$. 
%The $n$-fold positive product $\langle\alpha^n\rangle$ equals to its volume (introduced in \cite{boucksomVolume})
%    $$\vol(\alpha):=\sup_{T}\int_X T_{ac}^n$$
%    where $T$ ranges among all K\"ahler currents in $\alpha$. Of course, if there is no K\"ahler current in $\alpha$, its volume is defined to be $0$.

    \item The nef classes can be factored out in the positive products. Specifically, if $\alpha_1,...,\alpha_p\in\Nef^1(X)$ and $\beta_1,...,\beta_q\in \Psef^1(X)$, then we have $$\langle\alpha_1 \cdot ...\cdot\alpha_p\cdot \beta_1\cdot ...\cdot\beta_q\rangle=\alpha_1\cdot ...\cdot\alpha_p \cdot \langle\beta_1\cdots \beta_q\rangle.$$
    In particular, the positive product of several nef classes with a single movable class is just the usual intersection product.
            
\end{enumerate}
\end{prop}

\subsection{Numerical dimension}
A basic positivity invariant for nef classes is the numerical dimension.
\begin{defn}
    The numerical dimension of a nef class $\alpha \in \Nef^1(X)$ is 
    $$\nd(\alpha)=\max\{k\in \mathbb{ N}: \alpha^k\neq 0\}.$$
\end{defn}

Then a nef class $\alpha$ is 0 iff $\nd(\alpha)=0$.
Demailly-P\u{a}un \cite{DP04} proved that a nef class $\alpha$ is big iff $\nd(\alpha)=n$.

By the notion of positive product, we can extend the numerical dimension to arbitrary psef classes.

\begin{defn}
    The numerical dimension of a psef class $\alpha \in \Psef^1(X)$ is $$\nd(\alpha)=\max\{k\in \mathbb{ N}: \langle\alpha^k\rangle\neq 0\}.$$
\end{defn}

Then a psef class $\alpha$ with $\nd(\alpha)=0$ iff $\alpha = [N(\alpha)]$. Boucksom \cite{Bou02a} generalized \cite{DP04} to the psef setting: a psef class $\alpha$ is big iff $\nd(\alpha)=n$.

\subsection{Collections of nef classes}\label{sec collect nef}

For any positive integer $m$, we use the notation $[m]$ to denote the set $\{1,...,m\}$. For any collection $\{L_1,...,L_m\}$ of elements in a additive semigroup, we write 
$$L_I=\sum_{i\in I}L_i$$ 
for any subset $I\subset [m]$.

\begin{defn}
    Let $\mathcal{L}=\{L_1,...,L_m\}$ be a collection of nef classes where some of them are allowed to be the same. We define the numerical dimension of the collection to be $$\nd(\mathcal{L}):=\min_{\emptyset \neq I\subset [m]}\{\nd(L_I)-|I|+m\}.$$
    We denote the complete intersection class associated to the collection $\mathcal{L}$ by
    $$\mathbb{L}:=L_1\cdot ...\cdot L_m.$$
    The class $\mathbb{L}$ induces a linear operator by multiplying on cohomology classes, which we call a Lefschtez type operator.
\end{defn}

The following lemma figures out that $\nd(\mathcal{L})$ actually behaves like a numerical dimension.

\begin{lem}\label{hallrado}
We have $\mathbb{L} \neq 0$ if and only if 
    $\nd(L_I)\geq |I|, \forall I\subset [m]$.
\end{lem}

\begin{proof}
This follows from \cite[Theorem B]{hxhardlef}.
\end{proof}

\begin{lem}\label{propertynd}
    Let $\mathcal{L}=\{L_1,...,L_m\}$ be a collection of nef classes. Then 
    \begin{enumerate}
        \item $0\leq \nd(\mathcal{L})\leq \nd(L_{[m]})$ and $\nd(\mathcal{L})=0$ if and only if $L_i=0$ for each $i\in[m]$;
        \item If $L_1=\cdots =L_m=L$, then $\nd(\mathcal{L})=\nd(L)$;
        \item $\mathbb{L} \neq 0$ if and only if $\nd(\mathcal{L})\geq m$.
    \end{enumerate}
\end{lem}
\begin{proof}
    $(1)$ and $(2)$ follows directly from the definition. $(3)$ follows from Lemma \ref{hallrado} and the definition of $\nd(\mathcal{L})$.
\end{proof}

One of the basic properties of numerical dimensions of nef classes is the submodularity proved in \cite{hxinterineq}.

\begin{prop}\label{submodularity}
    Let $L,M,N\in \Nef^1(X)$ be nef classes. Then 
    $$\nd(L+M+N)+\nd(L)\leq \nd(L+M)+\nd(L+N).$$
\end{prop}

A consequence is the following:

\begin{lem}\label{maximalI}
    Let $\mathcal{L}=\{L_1,...,L_{m}\}$ be a collection of nef classes with $\nd(\mathcal{L})=m$. Then there exists a unique maximal index set $I_0\subset [m]$ such that $\nd(L_{I_0})=|I_0|$.
\end{lem}
\begin{proof}
 The condition  $\nd(\mathcal{L})=m$ is equivalent to that 
    $$\forall I, \nd(L_I)\geq |I| \text{ and } \exists J \text{ such that}\nd(L_J)=|J|.$$
    Hence, $\{ I\in  2^{[m]}\backslash\emptyset :\nd(L_{I})=|I| \}$ is a finite nonempty set which clearly contains finitely many maximal elements.
    To show the uniqueness, it suffices to show that $$\nd(L_{I\cup J})=|I\cup J|$$ for any $I,J\subset [m]$ with $\nd(L_I)=|I|$ and $\nd(L_J)=|J|$. As pointed out in the proof of Lemma \ref{propertynd}, for any $K$, $\nd(L_{K})\geq |K|$ is implied by $\nd(\mathcal{L})=m$. On the other hand, by the submodularity and the inequality $\nd(L_{I\cap J})\geq |I\cap J|$, we have 
    \begin{align*}
        \nd(L_{I\cup J})&\leq \nd(L_I)+\nd(L_J)-\nd(L_{I\cap J})\\
        &\leq |I|+|J|-|I\cap J|\\
        &=|I\cup J|
    \end{align*}
    
    This finishes the proof.
\end{proof}

\section{Lorentzian classes}\label{sec Lorenclass}
In this section, we introduce a key object -- Lorentzian classes -- in this paper. Throughout, we fix a compact K\"ahler manifold $X$ of dimension $n$.

\begin{defn}\label{deflor}
For $d\geq 2$, the class $\Omega\in H^{n-d,n-d}(X,\mathbb{R})$ is called a strictly Lorentzian class of dimension $d$ if 
\begin{enumerate}
    \item For any $\omega_1,...,\omega_{d}\in \Amp^1(X)$, we have 
            $\omega_1\cdot...\cdot\omega_{d}\cdot\Omega>0,$  and
    \item For any $\omega_1,...,\omega_{d-2}\in \Amp^1(X)$, the bilinear form 
            \begin{align*}
                &H^{1,1}(X,\mathbb{R})\times H^{1,1}(X,\mathbb{R}) \rightarrow \mathbb{R},\\ &(\alpha,\beta)\mapsto \alpha\cdot\beta\cdot\omega_1\cdot...\cdot\omega_{d-2}\cdot\Omega,
              \end{align*}
          has Lorentzian signature $(+,-,...,-)$.
\end{enumerate}
It is called a Lorentzian class if (2) is weaken by having exactly one positive eigenvalue. 
\end{defn}

The class $\Omega$ being strictly Lorentzian yields that $\omega_1\cdot...\cdot\omega_{d-2}\cdot\Omega$ has hard Lefschetz property, i.e., the map
\begin{equation*}
  \omega_1\cdot...\cdot\omega_{d-2}\cdot\Omega: H^{1,1}(X,\mathbb{R}) \rightarrow H^{n-1,n-1}(X,\mathbb{R})
\end{equation*}
is an isomorphism.

The following Lemma is well-known (see, for example, \cite[Proposition 3.4]{HuXiaoLefChar}).

\begin{lem}\label{equalAF}
    Let $\Gamma \in H^{n-2,n-2}(X,\mathbb{R})$ be a Lorentzian class. Then for any $\alpha,\beta \in \Nef^1(X)$, we have 
    $$(\alpha \cdot \beta \cdot\Gamma )^2\geq (\alpha^2 \cdot\Gamma )(\beta^2 \cdot\Gamma )$$
    with equality holds if and only if the $(n-1,n-1)$ classes $\alpha\cdot\Gamma$ and $\beta\cdot \Gamma$ are proportional.
\end{lem}

The assumption in Definition \ref{deflor} can be weaken to the following form. 

\begin{prop}\label{onetomixed}
The class  $\Omega\in H^{n-d,n-d}(X,\mathbb{R})$ is a strictly Lorentzian class if it is Lorentzian, and for any $\omega \in \Amp^1(X)$, the bilinear form 
                \begin{align*}
                &H^{1,1}(X,\mathbb{R})\times H^{1,1}(X,\mathbb{R}) \rightarrow \mathbb{R},\\ &(\alpha,\beta)\mapsto \alpha\cdot\beta\cdot\omega^{d-2}\cdot\Omega,
              \end{align*}
               is nondegenerate.
\end{prop}
\begin{proof}
It suffices to show that for any $\omega_1,...,\omega_{d-2}\in \Amp^1(X)$, the bilinear form 
    \begin{align*}
      &Q: H^{1,1}(X,\mathbb{R})\times H^{1,1}(X,\mathbb{R}) \rightarrow \mathbb{R},\\ &Q(\alpha,\beta)= \alpha\cdot\beta\cdot\omega_1\cdot...\cdot\omega_{d-2}\cdot\Omega,
    \end{align*}
is nondegenerate. 
By Poincar\'e duality, this is equivalent to show that the class $\alpha=0$ whenever
   \begin{equation}\label{321}
        \alpha \cdot \omega_1\cdot...\cdot\omega_{d-2}\cdot \Omega=0.
   \end{equation}
   By assumption, for any $\omega \in \Amp^1(X)$, $\omega \cdot \omega_2\cdot... \cdot \omega_{d-2}\cdot\Omega$ is a Lorentzian class which implies
   \begin{equation}\label{322}
    \alpha^2 \cdot \omega \cdot \omega_2\cdot... \cdot \omega_{d-2}\cdot\Omega\leq \frac{(\alpha\cdot \omega_1 \cdot \omega \cdot \omega_2\cdot... \cdot \omega_{d-2}\cdot\Omega)^2}{\omega_1^2 \cdot \omega \cdot \omega_2\cdot... \cdot \omega_{d-2}\cdot\Omega}=0.
   \end{equation}

   On the other hand, we can choose $\lambda>0$ such that $\lambda \omega_1- \omega\in \Amp^1(X)$. Replacing $\omega$ by $\lambda \omega_1- \omega$ in (\ref{322}), 
   \begin{equation}\label{323}
       0=  \lambda \omega_1\cdot \alpha^2 \cdot \omega_2\cdot... \cdot \omega_{d-2}\cdot\Omega \leq \omega \cdot \alpha^2 \cdot \omega_2\cdot... \cdot \omega_{d-2}\cdot\Omega
   \end{equation} 
   where the first equality follows from (\ref{321}). Combining (\ref{322}) and (\ref{323}) yields
   \begin{equation*}
        \omega\cdot \alpha^2 \cdot \omega_2\cdot... \cdot \omega_{d-2}\cdot\Omega=0
   \end{equation*}
   for any $\omega\in \Amp^1(X)$. This implies 
   \begin{equation*}
    \alpha^2 \cdot \omega_2\cdot... \cdot \omega_{d-2}\cdot\Omega=0,
   \end{equation*}
   since $\Amp^1(X)$ is open. A fortiori, we have
   \begin{equation*}
        \alpha^2 \cdot \omega_2^2 \cdot... \cdot \omega_{d-2}\cdot\Omega=0.
   \end{equation*}
   Together with 
   \begin{equation*}
        \alpha \cdot \omega_1 \cdot \omega_2^2 \cdot... \cdot \omega_{d-2}\cdot\Omega=0
   \end{equation*}
   which is a consequence of (\ref{321}), we obtain
   \begin{equation*}
    0= (\alpha \cdot \omega_1 \cdot \omega_2^2 \cdot... \cdot \omega_{d-2}\cdot\Omega)^2=(\alpha^2 \cdot \omega_2^2 \cdot... \cdot \omega_{d-2}\cdot\Omega)( \omega_1^2 \cdot \omega_2^2 \cdot... \cdot \omega_{d-2}\cdot\Omega).
   \end{equation*}
   Apply Lemma \ref{equalAF}, we know that $\alpha\cdot \omega_2^2 \cdot... \cdot \omega_{d-2}\cdot\Omega$ and $\omega_1\cdot \omega_2^2 \cdot... \cdot \omega_{d-2}\cdot\Omega \neq 0$ are proportional. Say 
   \begin{equation*}
    \alpha\cdot \omega_2^2 \cdot... \cdot \omega_{d-2}\cdot\Omega=\lambda \omega_1\cdot \omega_2^2 \cdot... \cdot \omega_{d-2}\cdot\Omega .
   \end{equation*}
   The constant $\lambda$ must be 0 since the LHS becomes 0 after pairing with $\omega_1$. That is,
   \begin{equation}\label{324}
    \alpha\cdot \omega_2^2 \cdot... \cdot \omega_{d-2}\cdot\Omega=0.
   \end{equation}

   What we showed above is the implication $(\ref{321})\implies (\ref{324})$. We can repeat the process and finally get 
   \begin{equation*}
    \alpha\cdot\omega_2^{d-2}\cdot\Omega=0
   \end{equation*}
   which implies $\alpha=0$ by the assumption.
\end{proof}

\subsection{Reformulations}

\subsubsection{Reformulation via Lorentzian polynomials}

We reformulate Lorentzian classes using the theory of Lorentzian polynomials \cite{branhuhlorentz,brandenleake23coneLorentzian}. This is one main reason that we use the terminology ``Lorentzian class''.

Recall the definition of cone Lorentzian polynomials introduced in \cite{brandenleake23coneLorentzian}.

\begin{defn}
Let $\mathfrak{C}$ be an open convex cone in $\mathbb{R}^s$. A homogeneous polynomial  $f$ (with $\deg f=n$) on $\mathbb{R}^s$ is called strictly $\mathfrak{C}$-Lorentzian if for all $v_1, ...,v_n \in \mathfrak{C}$,
\begin{itemize}
  \item $D_{v_1}... D_{v_n} f >0$, and
  \item the symmetric bilinear form $$(\xi, \eta) \mapsto D_\xi D_\eta D_{v_3}... D_{v_n} f$$ has Lorentzian signature $(+,-,...,-)$.
\end{itemize}
It is called $\mathfrak{C}$-Lorentzian if the second requirement is weaken by having exactly one positive eigenvalue. 
Here, $D_v$ is the directional derivative along $v$.
\end{defn}

When $\mathfrak{C}=\mathbb{R}_{>0} ^s$, these are the Lorentzian polynomials studied in \cite{branhuhlorentz}.
The above definition is equivalent to that for all positive integers $m$ and for all $v_1,...,v_m \in \mathfrak{C}$, the polynomial
\begin{equation*}
  (y_1,...,y_m)\mapsto f(y_1 v_1 +...+ y_m v_m)
\end{equation*}
is Lorentzian in the sense of \cite{branhuhlorentz} and has only positive coefficients.

Given $\Omega \in H^{n-d, n-d}(X, \mathbb{R})$, we can define the volume polynomial on $\Omega$ as 
$$f_{\Omega}: H^{1,1}(X, \mathbb{R})\rightarrow \mathbb{R}, \ \alpha \mapsto \alpha^d\cdot \Omega.$$
Then reformulating by Lorentzian polynomials on the K\"ahler cone, we have:

\begin{prop}
The class $\Omega$ is (strictly) Lorentzian iff the polynomial $f_\Omega$ is (strictly) $\Amp^1(X)$-Lorentzian.
\end{prop}

This viewpoint is helpful as we can apply Lorentzian polynomials to study Lorentzian classes, for example, as in \cite{hxinterineq}.

\subsubsection{Reformulation via restricted cohomology}
In what follows we omit the real-coefficient notation in our cohomology groups.
Let $\Omega$ be a cohomology classes of dimension $d$. It is illuminating to 
define the restricted cohomology on $\Omega$ as $$H^*(\Omega)=H^*(X)/\Ann(\Omega),$$ 
where $\Ann(\Omega)=\{\alpha\in H^*(X):\alpha\cdot\Omega=0\}$. 
It is endowed with a natural graded ring structure and a canonical perfect pairing
$$H^i(\Omega) \times H^{2d-i}(\Omega)\rightarrow \mathbb{R},\ (\alpha,\beta) \mapsto \alpha\cdot \beta\cdot \Omega$$
for each $0\leq i\leq 2d$. 
Moreover, $H^*(\Omega)$ has (not necessarily polarized) $\mathbb{R}$-Hodge structure.

We denote by $\Amp^1(\Omega)$ the image of $\Amp^1(X)$ in $H^{1,1}(\Omega)$.

\begin{defn}
    We say that $H^*(\Omega)$ satisfies the Hodge-Riemann relation of degree $\leq 1$ if the followings hold.
    \begin{enumerate}
        \item For any $\omega_1,...,\omega_{d}\in \Amp^1(\Omega)$, we have $\omega_1\cdots\omega_d\cdot \Omega>0$, and 
        \item For any $\omega_1,...,\omega_{d-2}\in \Amp^1(\Omega)$, the bilinear form 
        $$H^{1,1}(\Omega)\times H^{1,1}(\Omega) \rightarrow \mathbb{R}, \ (\alpha,\beta)\mapsto \alpha\cdot\beta\cdot\omega_1\cdots\omega_{d-2}\cdot\Omega$$
      has Lorentzian signature $(+,-,...,-)$.
    \end{enumerate}
    We say that $H^*(\Omega)$ satisfies the weak Hodge-Riemann relation if the second condition is weaken by having exactly one positive eigenvalue. 
\end{defn}

Then we have:
\begin{prop}
The class $\Omega$ is strictly Lorentzian iff $H^*(\Omega)$ satisfies the Hodge-Riemann relations of degree $\leq 1$ and $H^{1,1}(X) \cong H^{1,1}(\Omega)$. The class $\Omega$ is Lorentzian iff $H^*(\Omega)$ satisfies the weak Hodge-Riemann relations of degree $\leq 1$.
\end{prop}

Due to the discrepancy between $H^{1,1}(X)$ and $H^{1,1}(\Omega)$, we may encode the differences as follows.
\begin{defn}\label{intriLorent}
We call that $\Omega$ is \emph{extrinsically} strictly Lorentzian if it is strictly Lorentzian on $X$, and \emph{intrinsically} strictly Lorentzian if $H^*(\Omega)$ satisfies the Hodge-Riemann relation of degree $\leq 1$.
\end{defn}

The intuitive idea behind this reformulation is that we may think $\Omega$ as a combination of some cycle classes of dimension $d$ in $X$.

\subsection{Examples}\label{sec exmple lorent}
We present several examples giving (strictly) Lorentzian classes.

    \begin{enumerate}
        \item The non-zero products of nef classes $\Omega = \alpha_1\cdot...\cdot \alpha_{n-d}$ are Lorentzian classes. When all the $\alpha_i$ are K\"ahler, $\Omega$ is strictly Lorentzian. This fact follows from \cite{gromov1990convex, DN06, cattanimixedHRR}.

      \item By our previous work \cite{HuXiaoLefChar}, the products of big nef classes $\Omega = L_1\cdot...\cdot L_{n-d}$ are still strictly Lorentzian, whenever every augmented base locus (also known as non-ample/K\"ahler locus) satisfies $$\codim \mathbb{B}_+ (L_i)\geq 2.$$ 
      
        \item The fundamental class of any smooth irreducible subvariety $Z$ is Lorentzian. It is strictly Lorentzian if moreover the map 
            \begin{equation*}
                \cup [Z]: H^{1,1}(X)\rightarrow H^{n-\dim Z+1,n-\dim Z+1}(Z)
            \end{equation*}
            is injective.
        In the case when $Z$ is singular, its fundamental class is Lorentzian by applying a desingularization of $Z$.
            
        \item Schur classes $s_{\lambda}(E)$ of (ample) nef vector bundles are (strictly) Lorentzian classes by Ross-Toma \cite{ross2019hodge}.
        \item Let $\omega_1,...,\omega_e$ be K\"ahler classes and let $f$ be a nonzero dually Lorentzian polynomial (see Ross-Suss-Wannerer's work \cite{rsw23duallyLorentzian}) with $e$ variables and of degree $n-d$. Then $f(\omega_1,...,\omega_e)$ is a strictly Lorentzian class.

        This contains (1), and even (4) as special cases when $E=L_1\oplus\cdots \oplus L_{n-k}$,  since monomials and Schur polynomials are all dually Lorentzian.
        \item Let $A_{n}$ be the symmetric group of $n+1$ elements and let $X_{A_{n}}$ be the toric variety, known as permutahedron variety, associated to the permutahedron
                \begin{equation*}
                    \Xi_{n}=\conv \{ \sigma (0,1,...,n)\in \mathbb{R}^{n+1}:\sigma \in A_{n}\}
                \end{equation*}
                lying in the affine hyperplane $\sum_i x_i =\frac{n(n+1)}{2}$.
                Any loopless matroid $M$ with ground set $\{0,...,n\}$ and of rank $d+1$ determines a cohomology class $\Sigma_M\in H^{n-d,n-d}(X_{A_n})$, called the Bergman class. 
                As a consequence of the deep results by Adiprasito-Huh-Katz \cite{huhHRR}, $\Sigma_M$ is an \emph{intrinsically} strictly Lorentzian class. In fact, they showed that $H^*(\Sigma_M)$ satisfies Hodge-Riemann relations of all degrees. 
                
                This class can be represented by an effective algebraic cycle with several components. The matroid $M$ is linear (also known as representable or realizable) iff $\Sigma_M$ is the class of an irreducible subvariety \cite{huhTropMatroid}. Note that linear matroids are very sparse, as it is proved by Nelson \cite{nelsonMAtroid} that as $n$ tends to infinity, the proportion of $n$-element matroids that are linear tends to zero.
                Therefore, in an asymptotic sense, for almost all matroids $M$, $\Sigma_M$ can never be represented by an algebraic cycle with only one component.
                
                Similar construction also applies to the augmented Bergman classes for polymatroids on polystellahedral varieties \cite{polymatroidbergmanfan, polymatroidsHodge}.               
    \end{enumerate}

\subsection{Numerical dimension on a Lorentzian class}\label{sec ndlorent}

Assume that $\Omega \in H^{n-d, n-d}(X)$ is a Lorentzian class of dimension $d$.
Fix a K\"ahler class $\omega$, for $L\in \Nef^1 (X)$, we define its numerical dimension on $\Omega$ as $$\nd_{\Omega} (L) = \max \{k: L^k \cdot\omega^{d-k}\cdot \Omega \neq 0\}.$$
It is clear that $\nd_{\Omega} (L)$ is independent of the choices of $\omega$.

We have the following extension of Lemma \ref{hallrado}:

\begin{lem}\label{hallradoLorent}
Let $L_1, ..., L_m\in \Nef^1 (X)$ with $m\leq d$, then
\begin{equation*}
  L_1\cdot ...\cdot L_m \cdot \Omega\neq 0 \Leftrightarrow \nd_{\Omega} (L_I) \geq |I|,\ \forall I \subset [m].
\end{equation*}
\end{lem}

\begin{proof}
As $f_\Omega (\alpha)=\alpha^d \cdot\Omega$ is $\Amp^1 (X)$-Lorentzian, this follows exactly from \cite[Theorem 3.7]{HuXiaoLefChar}.
\end{proof}

\subsection{rKT property}\label{sec rkt}
As Lorentzian classes can be reformulated by Lorentzian polynomials, the rKT property proved in \cite{hxinterineq} transfers to Lorentzian classes. In particular, we have:

\begin{lem}\label{rkt lorent}
Let $\Omega$ be a Lorentzian class of dimension $d$ on $X$, then for any nef classes $A, B_1,...,B_d$ on $X$, we have that
\begin{equation*}
  (B_1 \cdot ...\cdot B_{d}\cdot \Omega)(A^d \cdot \Omega)\leq 2^{k(d-k)} (A^k\cdot B_{k+1} \cdot ...\cdot B_{d}\cdot \Omega)(A^{d-k}\cdot B_1 \cdots B_{k}\cdot \Omega).
\end{equation*}
\end{lem}

The rKT property plays an important role in quantitative stability estimates, see Section \ref{sec stabestim}. We refer the reader to \cite{hxinterineq} and the references therein for the improvement of the constant $2^{k(d-k)}$ when $\Omega$ satisfies additional requirement.

\section{The kernel face on the pseudoeffective cone}\label{sec nullface}

Throughout this section, we fix a compact K\"ahler manifold $X$ of dimension $n$.

Fix a collection of nef classes $\mathcal{L}=\{L_1, ...,L_m\}$, and denote the product class by $$\mathbb{L}=L_1 \cdot ...\cdot L_m.$$ 
In \cite{HuXiaoLefChar}, we systematically studied the case $m=n-2$. The kernel space
\begin{equation*}
 \mathbb{L}: H^{1,1}(X, \mathbb{R})\rightarrow H^{n-1,n-1}(X, \mathbb{R})
\end{equation*}
is closely related to the hard Lefschetz property of $\mathbb{L}$ and plays an important role in the algebro-geometric formalism of the extremals in log-concavity sequences. We gave a characterization of $\ker \mathbb{L}$
when the collection satisfies a positivity condition -- `superciticality with a rearrangement', that is,
\begin{equation*}
  \nd(L_i)\geq i+2,\ \forall i\in [n-2].
\end{equation*}
It was noted that $\ker \mathbb{L}$ is closely related to the kernel face 
on the pseudo-effective cone: $\ker \mathbb{L}\cap \Psef^1 (X)$.

The kernel face is also closely related the following extremal problem of intersection numbers.
If $\beta-\alpha \in \Psef^1(X)$ is psef, then
\begin{equation*}
    \mathbb{L}\cdot \alpha \leq \mathbb{L}\cdot \beta,
\end{equation*}
where $\leq$ means that $\mathbb{L}\cdot \beta - \mathbb{L}\cdot \alpha$ is a psef class in $H^{m+1, m+1} (X, \mathbb{R})$.
To understand when equality holds in the above monotonicity inequality, it is equivalent to the characterization of $\ker \mathbb{L}\cap \Psef^1 (X)$ -- the (\textbf{KF}) problem as mentioned in the introduction. It is easy to see that for any K\"ahler class $\omega$, $$\alpha\in \Psef^1 (X) \cap \ker \mathbb{L} \Leftrightarrow \alpha\in \Psef^1 (X) \cap \ker (\mathbb{L} \cdot \omega^{n-m-1}).$$
Therefore, we only need to consider the case when $\mathbb{L}$ is given by a product of $(n-1)$ nef classes.

In this section, we give a complete description of this kernel face.

\subsection{Numerical dimension and positivity}

We begin with the following result relating the positivity and numerical dimensions of two classes. 
The following bigness criterion due to \cite{xiao2014movable} is useful to us.

\begin{lem}
    Let $\alpha,\beta \in\Mov^1(X)$ be two movable classes. If $\langle\alpha^n\rangle-n\langle\alpha^{n-1}\cdot \beta\rangle >0$, then $\alpha-\beta$ is big.
\end{lem}

Our main result in this section is:

\begin{thrm}\label{annihilation}
    Let $L\in \Nef^1(X)$ be a nef class of numerical dimension $k$ and $\alpha \in \Mov^1(X)$ a movable class. Then the following statements are equivalent:
    \begin{enumerate}
        \item $L^k\cdot \alpha =0$;
        \item $\nd(L+\alpha)=\nd(L)$;
        \item There exists $ \lambda >0$ such that $\alpha \leq \lambda L$.
    \end{enumerate}
\end{thrm}

\begin{proof}
Fix a K\"ahler class $\omega$. 

We first show  $(2)\implies(1)$. By the assumption and the superadditivity of the positive products, we have
\begin{align*}
     0=\langle(L+\alpha)^{k+1}\rangle \cdot \omega^{n-k-1}&\geq \sum_{i=0}^{k+1}\binom{k+1}{i}\langle L^i \cdot\alpha^{k+1-i} \rangle \cdot \omega^{n-k-1}\\
     &\geq \langle L^{k}
     \cdot\alpha\rangle \cdot \omega^{n-k-1} \geq 0.
\end{align*}
    Since $L$ is nef and $\alpha$ is movable, $\langle L^{k}\cdot\alpha \rangle=L^k\cdot \alpha$, so the implication follows immediately.
    
For $(3)\implies(2)$, by the monotonicity of the numerical dimension, we obtain
    $$\nd(L)\leq \nd(L+\alpha)\leq \nd\Big((\lambda+1)L\Big)=\nd(L).$$ 
    Thus, $\nd(L+\alpha)=\nd(L)$.

Last we show $(1)\implies(3)$.
    Note that the assumption can be summarized as
    $$L^k \neq 0,L^{k+1}=0 , \text{ and }L^k \cdot \alpha=0.$$
\textbf{Claim:} there is $\lambda_0>0$ such that the intersection number
    $$(\lambda L+t \omega)^n-n(\lambda L+t \omega)^{n-1}\cdot \alpha>0$$
    for any $\lambda > \lambda_0$ and any $t>0$. Then by the bigness criterion of \cite{xiao2014movable}, the class $\lambda L-\alpha +t\omega$ is big for any $t>0$. Letting $t \rightarrow 0$ shows that $\lambda L-\alpha $ is psef. 
    
    It remains to prove the claim. This follows from the assumption as follows.
    Extending via the binomial formula, we have
    \begin{align*}
        &(\lambda L+t \omega)^n-n(\lambda L+t \omega)^{n-1}\cdot \alpha \\
        &=\sum_{i=0}^n \binom{n}{i}(L^{n-i}\cdot\omega^i)\lambda^{n-i}t^i-
        n\sum_{i=0}^{n-1}\binom{n-1}{i}(L^{n-1-i}\cdot\omega^i \cdot\alpha)\lambda^{n-1-i}t^i\\
        &=\sum_{i=n-k}^n\Big(\binom{n}{i}(L^{n-i}\cdot\omega^i)\lambda^{n-i}
        -n\binom{n-1}{i}(L^{n-1-i}\cdot\omega^i\cdot\alpha)\lambda^{n-1-i}\Big)t^i
    \end{align*}
    where we adopt the convention $\binom{n-1}{n}=0$.
    Now, choose $$\lambda >\lambda _0 :=\max\{(n-i)\frac{(L^{n-1-i}\cdot\omega^{i}\cdot\alpha)}{L^{n-i}\cdot\omega^i}:i=n-k,...,n\}.$$
    Then the polynomial
    $$(\lambda L+t \omega)^n-n(\lambda L+t \omega)^{n-1}\cdot \alpha$$
    in variable $t$ has positive coefficients, so it is positive for all $t>0$.

    This finishes the proof.
\end{proof}

Note that Theorem \ref{annihilation} implies the following: if adding a movable class $\alpha$ to a nef class $L$ does not increase the numerical dimension, then the same holds after replacing $L$ by $L+M$ for any other nef class $M$. Combining the submodularity of numerical dimensions (see Proposition \ref{submodularity}), this observation motivates the following conjecture.

\begin{conj}
    Let $L,M$ be nef classes and $\alpha$ a movable class. Then 
    \begin{equation*}
        \nd(L+M+\alpha)- \nd(L+M)\leq \nd(L+\alpha) -\nd(L).
    \end{equation*}
\end{conj}
This would generalize \cite{hxinterineq} with one class being just movable. We also conjecture the submodularity when $L, M$ are just movable.

\subsubsection{Geometric description in the semiample setting}
When $X$ is a smooth projective variety and $L$ is the first Chern class of a semiample line bundle, we may take $\Mov^1(X)\subset N^1(X)$ to be the movable cone generated by movable $\mathbb{R}$-divisors. In this case, Theorem \ref{annihilation} admits a more geometric proof. In fact, the face
$$\{\alpha \in \Mov^1(X):L^{\nd(L)}\cdot \alpha =0\}$$
has a more refined description, as shown by Lehmann \cite{lehmannpullback}.
More precisely, if $\varphi:X\rightarrow Y$ is the semiample fibration of $L$, then \cite{lehmannpullback} shows that for every $\alpha \in \Mov^1(X)$ satisfying $$L^{\nd(L)}\cdot \alpha=0,$$ 
there exists a diagram
\[
\begin{tikzcd}
X' \arrow{r}{\pi} \arrow{d}{\varphi'} & X \arrow{d}{\varphi} \\
Y' \arrow{r}{\mu}& Y
\end{tikzcd}
\]
where $\pi,\mu$ are birational morphisms, and $\beta \in \Psef^1(Y')$ satisfies that $$\alpha =\pi_*\langle\varphi'^*\beta\rangle.$$ 
That is, the face $\{\alpha \in \Mov^1(X) : L^{\nd(L)}\cdot \alpha=0\}$ essentially arises via pullback. 
It is natural to ask whether the same result holds in the transcendental setting.

\begin{conj}
    Let $\varphi:X\rightarrow Y$ be a surjective holomorphic map between compact K\"ahler manifolds of dimensions $n$ and $k$, respectively. Let $A$ be a K\"ahler class on $Y$ and set $L=\varphi^*A$. Then for any movable class $\alpha $ on $X$ with $L^k\cdot \alpha=0$, there exits a diagram
    \[
    \begin{tikzcd}
X' \arrow{r}{\pi} \arrow{d}{\varphi'} & X \arrow{d}{\varphi} \\
Y' \arrow{r}{\mu}& Y
\end{tikzcd}
\]
where $\pi,\mu$ are bimeromorphic holomorphic maps, together with a psef class $\beta$ on $Y'$ such that $$\alpha=\pi_*\langle\varphi'^*\beta\rangle.$$
\end{conj}

\subsection{Characterization of the kernel face}

In the following of this section, we set $\mathcal{L}=\{L_1,...,L_{n-1}\}$ a collection of $(n-1)$ nef classes, and recall we write $\mathbb{L}=L_1\cdots L_{n-1}$. We shall give a characterization of the kernel face $\ker \mathbb{L} \cap \Psef^1 (X)$.

We are only interested in the case $\mathbb{L}\neq0$, which by Lemma \ref{propertynd} is equivalent to $$\nd(\mathcal{L})\geq n-1\ (i.e., \nd(\mathcal{L})\in \{n-1,n\}).$$

\subsubsection{Using Zariski decomposition}
Given $\alpha\in \Psef^1 (X)$, as recalled in Section \ref{sec pre} we have the divisorial Zariski decomposition 
$$\alpha = P(\alpha)+N(\alpha)$$
such that $P(\alpha)\in \Mov^1 (X)$ and $N(\alpha)=\sum_i c_i [D_i]$ is a non-negative combination of prime divisors. Then $\alpha \in \ker \mathbb{L}$ iff $P(\alpha) \in \ker \mathbb{L}$ and $[D_i] \in \ker \mathbb{L}$ for any $i$. The condition $\mathbb{L}\cdot [D_i]=0$ can be characterized by the numerical dimensions of $\mathcal{L}$ on $D_i$. By Lemma \ref{hallrado} we have that $\mathbb{L}\cdot [D_i]=0$ iff for some $I\subset [n-1]$, $$L_I ^{|I|} \cdot [D_i] =0.$$
Note that the $D_i$ do not appear if $\alpha$ is movable. We call these hypersurfaces from the negative parts exceptional.

It remains to study the kernel face on the movable cone: $\ker \mathbb{L}\cap \Mov^1 (X)$. We divide the discussions into two cases. We first treat the easier case $\nd(\mathcal{L})=n$.

\subsubsection{Case: $\nd(\mathcal{L})=n$}

\begin{prop}\label{criticalinter}
    If $\nd(\mathcal{L})=n$, then $\ker \mathbb{L}\cap \Mov^1 (X)=\{0\}.$
\end{prop}
\begin{proof}
    Let $\omega$ be a fixed K\"ahler class.

     Assume $\alpha \in \Mov^1(X)$ is a movable class such that $$\alpha\cdot L_1\cdot L_2\cdot...\cdot L_{n-1}=0.$$
    By Hodge index theorem, this implies that
    \begin{equation}\label{hodgeeq}
    0=(\alpha\cdot L_1\cdot L_2\cdot...\cdot L_{n-1})^2 \geq (\alpha^2\cdot L_2\cdot...\cdot L_{n-1})(L_1^2\cdot L_2\cdot...\cdot L_{n-1}).
    \end{equation}
    On the other hand, since for any psef class $\beta$ and any pushforward class $\pi_* A$ via birational morphism, 
    $$\beta\cdot  \pi_* A \cdot L_2\cdot...\cdot L_{n-1} = \pi^* \beta \cdot A\cdot \pi^*L_2\cdot...\cdot \pi^*L_{n-1} \geq 0,$$
    and $\alpha$ is a limit of the classes of the form $\pi_* A$, we have
    $$\alpha^2\cdot L_2\cdot...\cdot L_{n-1}\geq 0.$$
    Therefore, (\ref{hodgeeq}) is an equality.
    
    Given that $\nd(\mathcal{L})=n$, by Lemma \ref{hallrado} we have $$L_1^2\cdot L_2\cdot...\cdot L_{n-1}>0.$$ 
    By Lemma \ref{equalAF}, there exists $c\in \mathbb{R}$ such that $$\alpha\cdot L_2\cdot...\cdot L_{n-1}=cL_1\cdot...\cdot L_{n-1}.$$
    Multiplying both sides by $L_1$ yields 
    $$0=\alpha\cdot L_1 \cdot...\cdot L_{n-1}=c L_1^2\cdot L_2 \cdot...\cdot L_{n-1},$$
    which forces $c=0$. Thus, $$\alpha\cdot L_2\cdot...\cdot L_{n-1}=0.$$ A fortiori, we have $$\alpha\cdot \omega \cdot L_2\cdot...\cdot L_{n-1}=0.$$

    Repeating the same process, we eventually obtain $\alpha\cdot\omega^{n-1}=0$, implying  $\alpha=0$.

    This finishes the proof.
\end{proof}

Next we treat the case $\nd(\mathcal{L})=n-1$.

\subsubsection{Case: $\nd(\mathcal{L})=n-1$}

The following result helps us to reduce the question of a collection of nef classes to a single nef class. 

\begin{prop}\label{abshallrado}
    Let $\alpha \in \Mov^1(X)$ be a non-zero movable classes. Then 
    $$\mathbb{L}\cdot \alpha =L_1\cdot...\cdot L_{n-1}\cdot\alpha >0 \Leftrightarrow \forall I \subset [n-1],  L_I^{|I|}\cdot\alpha \neq 0.$$
    Equivalently, $\mathbb{L}\cdot \alpha =0$ iff there is some $I \subset [n-1]$ such that $$L_I^{|I|}\cdot\alpha \neq 0.$$
\end{prop}

\begin{proof}
It is a direct consequence of \cite[Theorem 3.7]{HuXiaoLefChar}, once we know the polynomial function given by the intersection product
$$f: H^{1,1}(X,\mathbb{R})\rightarrow \mathbb{R}, \ L \mapsto L^{n-1}\cdot \alpha.$$ 
is an $\Amp^1(X)$-Lorentzian polynomial. Note that $\nd_f(L)$ therein is nothing but $$\max \{k\geq 0:L^k\cdot \alpha \neq0\},$$ 
then $L_I^{|I|}\cdot\alpha \neq 0$ means that $\nd_f(L_I) \geq |I|$.

We verify the Lorentzian property of $f$, that is, for any $A_1,..., A_{n-1} \in \Amp^1 (X)$
\begin{enumerate}
  \item The intersection number $A_1 \cdot ...\cdot  A_{n-1}\cdot \alpha >0$.
  \item The quadratic form 
  \begin{equation*}
    H^{1,1}(X,\mathbb{R}) \times H^{1,1}(X,\mathbb{R}) \rightarrow \mathbb{R},\ (C, D) \mapsto C\cdot D\cdot A_1\cdot...\cdot A_{n-3}\cdot \alpha
  \end{equation*}
  has exactly one positive eigenvalue.
\end{enumerate}
(1) is obvious. To prove (2), using the fact that $\alpha \in \Mov^1 (X)$ implies that $\alpha$ can be written as a limit of the form $\pi_* \widehat{A}$, where $\pi: Y\rightarrow X$ is a birational morphism and $\widehat{A}$ is K\"ahler on $Y$, we get for any $C\in \Nef^1 (X)$,
\begin{equation*}
  (C\cdot D\cdot A_1\cdot...\cdot A_{n-3}\cdot \alpha)^2 \geq (C^2\cdot  A_1\cdot...\cdot A_{n-3}\cdot \alpha)^2( D^2\cdot A_1\cdot...\cdot A_{n-3}\cdot \alpha)^2.
\end{equation*}
This clearly implies (2), finishing the proof that $f$ is an $\Amp^1(X)$-Lorentzian polynomial.
\end{proof}

We are now ready to settle down the case $\nd(\mathcal{L})=n-1$.
As shown in Lemma \ref{maximalI}, there exists a maximal index set $I_0\subset [n-1]$ such that $\nd(L_{I_0})=|I_0|$.
\begin{prop}\label{kern-1}
    If $\nd(\mathcal{L})=n-1$, then 
    \begin{align*}
       \ker \mathbb{L}\cap \Mov^1 (X)
        &=\sum_{I\subset[n-1]:\nd(L_{I})=|I|}\left\{\alpha\in \Mov^1(X) : L_{I}^{|I|}\cdot \alpha =0\right\}\\
        &=\left\{\alpha\in \Mov^1(X) : L_{I_0}^{|I_0|}\cdot \alpha =0\right\}.
    \end{align*}
\end{prop}
\begin{proof}
    By Proposition \ref{abshallrado}, we have
    \begin{equation}\label{ker1}
       \ker \mathbb{L}\cap \Mov^1 (X)
        =\sum_{I\subset[n-1]}\left\{\alpha\in \Mov^1(X) : L_{I}^{|I|}\cdot \alpha =0\right\}.
    \end{equation}

    We claim that for those summands $I\subset [n-1]$ such that $\nd(L_{I})\geq |I|+1$, we have 
    $$\{\alpha\in \Mov^1(X) : L_{I}^{|I|}\cdot \alpha =0\}=\{0\}.$$
    This follows from Proposition \ref{criticalinter}. In fact, $L_{I}^{|I|}\cdot \alpha=0$ is equivalent to $$\omega^{n-|I|-1}\cdot L_{I}^{|I|}\cdot \alpha=0$$ for any (or some) K\"ahler class $\omega$. The collection 
    $$\{L_I,...,L_I,\omega,...,\omega\}$$
    with $L_I$ appears for $|I|$ times and $\omega$ appears for $n-|I|-1$ times is of numerical dimension $n$.

    Hence, the only nonzero summands in (\ref{ker1}) are those given by $I\subset [n-1]$ with $\nd(L_I)=|I|$. In other words, we prove 
    $$\ker \mathbb{L}\cap \Mov^1 (X)
    =\sum_{I\subset[n-1]:\nd(L_{I})=|I|}\left\{\alpha\in \Mov^1(X) : L_{I}^{|I|}\cdot \alpha =0\right\}.$$

    Let $I\subset [n-1]$ be a subset with $\nd(L_I)=|I|$ and let $\alpha \in \Mov^1(X)$ be a movable class annihilated by $L_I^{|I|}$.
    By Theorem \ref{annihilation}, there exists $\lambda>0$ such that $\lambda \alpha \leq L_I $. By the maximality of $I_0$, we have $$\alpha \leq\lambda L_I \leq \lambda L_{I_0}.$$
    This implies that $L_{I_0}^{|I_0|}\cdot \alpha=0$ again by Theorem \ref{annihilation}. Since $I$ is an arbitrary index set with $\nd(L_I)=|I|$ and $\alpha \in \Mov^1(X)$ is an arbitrary movable class annihilated by $L_I^{|I|}$, we find
    $$\sum_{I\subset[n-1]:\nd(L_{I})=|I|}\left\{\alpha\in \Mov^1(X) : L_{I}^{|I|}\cdot \alpha =0\right\}\subset \left\{\alpha\in \Mov^1(X) : L_{I_0}^{|I_0|}\cdot \alpha =0\right\}.$$
    The converse inclusion is clear since the right hand side appears to be a summand of the left hand side.

    This finishes the proof.
\end{proof}

Combining all the above discussions, the characterization of $\ker \mathbb{L}\cap \Psef^1(X)$ can be summarized as follows:
\begin{thrm}\label{monotonicity}
        Let $\mathcal{L}=\{L_1,...,L_{n-1}\} $ be a collection of nef classes. Then we have:
    \begin{enumerate}
        \item If $\nd(\mathcal{L})\leq n-2$, then $\mathbb{L}=0$ and thus $\ker \mathbb{L}\cap \Psef^1(X) =\Psef^1 (X)$.
        \item If $\nd(\mathcal{L})=n-1$, then $\ker \mathbb{L}\cap \Psef^1(X)$ is generated by the classes of exceptional hypersurfaces $D$ such that $$\nd_D (L_{I}) < |I|$$ 
            for some $I\subset [n-1]$ and the movable classes $M$ such that 
            $$\nd(L_I+M)=\nd(L_I)$$ 
            for some $I\subset [n-1]$ with $\nd(L_I)=|I|$. 
        %$$\ker \mathbb{L}\cap \Mov^1(X)=\sum_{I\subset [n-1]:\nd(L_I)=|I|}\{M \in \Mov^1(X):\nd(L_I+M)=\nd(L_I)\}$$.
        \item If $\nd(\mathcal{L})=n$, then $\ker \mathbb{L}\cap \Psef^1(X)$ is generated by the classes of exceptional hypersurfaces $D$ such that $$\nd_D (L_{I}) < |I|$$ 
            for some $I\subset [n-1]$.
    \end{enumerate}
    Here $\nd_D (-)$ is the numerical dimension on the hypersurface $D$.
\end{thrm}

As a byproduct, we also get the following positivity criterion of some nef classes allowing one to be movable, generalizing the Hall-Rado relation in \cite{hxhardlef}.

\begin{cor}\label{hallradomov1}
    Let $L_1,...,L_{k-1} \ (1\leq k\leq n) $ be nef classes and let $L_{k}$ be a movable class. Then
    $$L_1\cdot...\cdot L_k \neq 0 \iff \nd(L_I)\geq |I|,\forall I\subset [k].$$
\end{cor}
\begin{proof}
 We first consider the implication ``$\Rightarrow$''.  
 Suppose there exists an index set $I\subset [k]$ such that $\nd(L_I)<|I|$. By the monotonicity of the positive product, we have
    \begin{equation*}
        \prod_{i\in I}L_i=\langle\prod_{i\in I}L_i\rangle \leq \langle L_I^{|I|}\rangle =0.
    \end{equation*}
    Consequently, $L_1\cdot...\cdot L_k=0$.

    We next prove the converse. That is, if $L_1\cdot...\cdot L_k=0$, then there exists an index set $I\subset [k]$ with $\nd(L_I)<|I|$. Fix a K\"ahler class $\omega$. The given condition implies
    \begin{equation*}
        \omega^{n-k}\cdot L_1\cdot...\cdot L_k=0.
    \end{equation*}
    The equality indicates that the movable class $L_k$ lies in the kernel of the operator associated with the product of classes in the collection
    $$\mathcal{L}=\{\omega,...,\omega,L_1,...,L_{k-1}\}$$
    with $\omega$ appearing for $n-k$ times. We divide the discussions into three cases based on Theorem \ref{monotonicity}.
    \begin{enumerate}
        \item If $\nd(\mathcal{L})\leq n-2$, then $$\omega^{n-k} \cdot L_1 \cdot...\cdot L_{k-1}=0.$$ This implies $L_1 \cdot...\cdot L_{k-1}=0$, and by Lemma \ref{propertynd}, there exists $I\subset [k-1]$ with $\nd(L_I)<|I|$.
        \item If $\nd(\mathcal{L})= n-1$, then there exists $I\subset [k-1]$ and $0\leq l\leq n-k$ such that $$\nd(L_I+l\omega+L_k)=|I|+l.$$
        We must have $l=0$ since $|I|+l<n$ and $\omega$ is K\"ahler. Thus, $I\cup \{k\}$ satisfies $$\nd(L_{I\cup\{k\}})=\nd(L_I+L_k)=|I|<|I|+1.$$
        \item If $\nd(\mathcal{L})= n$, then we must have $L_k=0$ in which case we just take the index set to be $\{k\}$.
    \end{enumerate}

    This finishes the proof.
\end{proof}

We conjecture that Corollary \ref{hallradomov1} extends to positive products for an arbitrary movable collection without the nefness assumption.
 
\begin{conj}\label{hallradoallmov}
    Let $L_1,...,L_{k} \ (1\leq k\leq n) $ be movable classes. Then
    $$\langle L_1\cdot...\cdot L_k \rangle \neq 0 \iff \nd(L_I)\geq |I|,\forall I\subset [k].$$
\end{conj}
The same question could also be asked for a collection of psef classes, which seems subtler due to the intractability of the divisorial Zariski decomposition for sums of psef classes.

\subsection{Extensions}
It is possible to extend Theorems \ref{annihilation}, \ref{monotonicity} to a Lorentzian class with some technical assumptions (e.g., birational invariance), which we hopefully intend to develop in a future project. As a baby example, Theorem \ref{annihilation} can be generalized as follows.

\begin{prop}\label{ndpos}
    Let $\Omega$ be a Lorentzian class of dimension $d$ and fix a K\"ahler class $\omega$. Let $L$ be a nef class of numerical dimension $\nd_{\Omega} (L) =k$ and let $\alpha$ be a nef class. Then the following statements are equivalent:
    \begin{enumerate}
        \item $L^k\cdot \alpha\cdot\omega^{d-k-1}\cdot \Omega =0$;
        \item $\nd_\Omega(L+\alpha)=\nd_\Omega(L)$;
        \item There exists $ \lambda >0$ such that for any $B_1, ..., B_{d-1}\in \Nef^1 (X)$, 
            \begin{equation*}
              (\lambda L-\alpha)\cdot  B_1 \cdot ...\cdot B_{d-1}\cdot \Omega \geq 0. 
           \end{equation*}
    \end{enumerate}
\end{prop}

Here $\nd_\Omega(L)$ is numerical dimension defined in Section \ref{sec ndlorent}:
$$\nd_{\Omega} (L) = \max \{k: L^k \cdot\omega^{d-k}\cdot \Omega \neq 0\}.$$

\begin{proof}
The implications $(2)\Rightarrow (1)$, $(3)\Rightarrow (2)$ are clear.

For the implication $(1)\Rightarrow (3)$, using the same argument as in Theorem \ref{annihilation}, there is some $\lambda_0 >0$ such that for any $\lambda > \lambda_0, t>0$,
\begin{equation*}
  \frac{(\lambda L+t \omega)^d \cdot \Omega}{2^{d-1}(\lambda L+t \omega)^{d-1}\cdot \alpha\cdot \Omega}  >1.
\end{equation*}

Therefore,
\begin{align*}
   (\lambda L + t\omega-\alpha)\cdot  B_1 \cdot ...\cdot B_{d-1}\cdot \Omega
   \geq \left(\lambda L + t\omega- \frac{(\lambda L+t \omega)^d \cdot \Omega}{2^{d-1}(\lambda L+t \omega)^{d-1}\cdot \alpha\cdot \Omega}\alpha\right)\cdot  B_1 \cdot ...\cdot B_{d-1}\cdot \Omega.
\end{align*}
By Lemma \ref{rkt lorent}, the RHS of the above inequality is non-negative.
Letting $t$ tend to 0 implies the result.
\end{proof}

\begin{comment}
***TO BE CLEARED OUT***
\begin{lem}\label{presmov}
    Let $X$ be a smooth projective variety and $A$ be a very ample divisor. 
    The restriction of a movable class to a very general hypersurface $H\in |A|$ is still movable.

\end{lem}

\begin{lem}\label{presnd}
    Let $X$ be a smooth projective variety and $A$ be a very ample divisor. Then for any psef class $\alpha \in \mathcal{E}$ and any integer $0 \leq k \leq n$,
    we have $$\langle\alpha^k \rangle \cdot A^{n-k} \leq \langle\alpha_{|V}^k\rangle$$
    where $V=\bigcap_{i=1}^{n-k}H_i$ is the complete intersection of $n-k$ very general hypersurfaces $H_1,...,H_{n-k}\in|A|$.
\end{lem}
***

\end{comment}

\section{Proportionality on a strictly Lorentzian class}\label{sec propor}

In this section, we study the extremals for ($\textbf{BM}_{\Omega}$) on a strictly Lorentzian class $\Omega$.

\subsection{Log-concavity}

In this section, on a compact K\"ahler manifold $X$ of dimension $n$, we fix a Lorentzian class $\Omega$ of dimension $d$ with $d\geq 2$. Recall that we adopt the following terminologies.

\begin{defn}
For a nef class $\alpha$ on $X$, its volume on $\Omega$ is defined by $$\vol_\Omega (\alpha):=\alpha^d \cdot \Omega.$$ 
For a nef class $\alpha$ on $X$, it is called big on $\Omega$ if $\vol_\Omega (\alpha)>0$. 
\end{defn}

For simplicity, we also denote the volume by $|\alpha|_{\Omega}$, or just $|\alpha|$ when the reference Lorentzian class $\Omega$ is clear. 
Given two nef classes $\alpha, \beta$, we set $$s_{k} = \alpha^k \cdot  \beta^{d-k}\cdot \Omega,\ 0\leq k \leq d.$$

The following lemma is well-known, see e.g. \cite{fx19}. We involve it here for completeness.

\begin{lem}\label{logconc}
    For any nef classes $\alpha $ and $ \beta$, we have the following properties:
\begin{enumerate}
  \item $s_k ^2 \geq s_{k-1} s_{k+1}, \forall\ 1\leq k\leq d-1$.
  \item $s_{k_2} \geq s_{k_1}^t s_{k_3}^{1-t}, \forall\ 0\leq k_1\leq k_2 \leq k_3\leq d,$ where $t=\frac{k_3 -k_2}{k_3 -k_1}$.
  \item $s_k ^d \geq s_d ^k s_0 ^{d-k}, \forall\ 0\leq k\leq d$. 
  \item $s_{d-1} ^d \geq s_d ^{d-1 } s_0$.
  \item $|\alpha+\beta|^{1/d}\geq  |\alpha|^{1/d} +|\beta|^{1/d}$.
\end{enumerate}

Moreover, assume further that  $\alpha $ and $ \beta$ are big on $\Omega$, then the following are equivalent:
\begin{enumerate}
  \item $s_k ^2 = s_{k-1} s_{k+1}, \forall\ 1\leq k\leq d-1$.
  \item $s_{k_2} = s_{k_1}^t s_{k_3}^{1-t}, \forall\ 0\leq k_1\leq k_2 \leq k_3\leq d,$ where $t=\frac{k_3 -k_2}{k_3 -k_1}$.
  \item $s_k ^d = s_d ^k s_0 ^{d-k}, \forall\ 0\leq k\leq d$. 
  \item $s_k ^d = s_d ^k s_0 ^{d-k}$, for some $\ 1\leq k\leq d-1$. 
  \item $s_{d-1} ^d = s_d ^{d-1 } s_0$.
  \item $|\alpha+\beta|^{1/d}=  |\alpha|^{1/d} +|\beta|^{1/d}$.
\end{enumerate}
\end{lem}

By the equivalence of (1) and (6) in the second part, we see that the extremals for ($\textbf{BM}_{\Omega}$) coincide with the extremals for the log-concave sequence $\{s_k\}$ being flat.

\begin{proof}
We first consider the properties on inequalities. By taking limits, we can assume that $\alpha, \beta$ are big on $\Omega$, then $s_k>0$ for all $k$.

First, (1) follows from Lemma \ref{equalAF}, since $\alpha^{k-1}\cdot \beta^{d-k-1}\cdot \Omega$ is a Lorentzian class of dimension 2. 

Next we prove the inequalities in (2). 
By (1), the sequence 
\begin{equation*}\label{skdecrea}
\left\{\frac{s_{k}}{s_{k-1}}\right\}_{1\leq k \leq d}
\end{equation*}
is non-increasing. For $0\leq k_1\leq k_2 \leq k_3\leq d$, write
$$\frac{s_{k_3}}{s_{k_2}} =\frac{s_{k_3}}{s_{k_3 -1}}\cdots  \frac{s_{k_2 +1}}{s_{k_2}}$$
and
$$\frac{s_{k_2}}{s_{k_1}} =\frac{s_{k_2}}{s_{k_2 -1}}\cdots  \frac{s_{k_1 +1}}{s_{k_1}}.$$
Then $\frac{s_{k_3}}{s_{k_2}}$ (respectively, $\frac{s_{k_2}}{s_{k_1}}$) is a product of $k_3 - k_2$ (respectively, $k_2 -k_1$) factors, and by (1) every factor in $\frac{s_{k_3}}{s_{k_2}}$ is bounded above by any factor in  $\frac{s_{k_2}}{s_{k_1}}$. This yields that
$$\left(\frac{s_{k_3}}{s_{k_2}}\right)^{k_2 -k_1} \leq  \left(\frac{s_{k_2}}{s_{k_1}}\right)^{k_3 -k_2},$$
which is exactly (2).

The inequalities in (3), (4) follow from (2), and (5) follows from (3). This finishes the proof of the first part.

Now we consider the second part on the equivalent descriptions on equalities. 
This follows straightforward from the proof of the inequalities. Taking the implication $(4)\Rightarrow (1)$ as an example, by the analysis for the inequality (2) in the first part, if $s_k ^d = s_d ^k s_0 ^{d-k}$ for some $\ 1\leq k\leq d-1$, then any factor in 
$$\frac{s_{d}}{s_{k}} =\frac{s_{d}}{s_{d -1}}\cdots  \frac{s_{k +1}}{s_{k}}$$
must equal to an arbitrary factor in
$$\frac{s_{k}}{s_{0}} =\frac{s_{k}}{s_{k -1}}\cdots  \frac{s_{1}}{s_{0}},$$
that is, for any $j> k$,
$$\frac{s_{j}}{s_{j-1}}=\frac{s_{k}}{s_{k-1}}=...=\frac{s_{1}}{s_{0}}.$$
This is exactly (1). The other implications are similar, we omit the details.

This finishes the proof.
\end{proof}

\subsection{Proportionality}\label{sec teisprop}

Before giving the proof of the main result. We first discuss the positivity assumptions on $\alpha,\beta$ in order to conclude proportionality for the extremals of ($\textbf{BM}_{\Omega}$). When both $\alpha$ and $\beta$ are not big on $\Omega$, in general we cannot conclude anything related to proportionality. It is easy to find examples such that $|\alpha+\beta| =0 =|\alpha|=|\beta|$ even if $\alpha, \beta$ are not proportional. When one is big and the other is not, we note that the strict inequality holds:

\begin{prop}\label{BMeqnotbig}
Let $\Omega$ be a strictly Lorentzian class of dimension $d$, and let $\alpha,  \beta$ be nef classes on $\Omega$. Assume that $|\alpha|>0, |\beta|=0$ and $\beta \neq 0$, then we will have the strict inequality:
\begin{equation*}
  |\alpha+\beta|^{1/d}> |\alpha|^{1/d}.
\end{equation*}
\end{prop}

\begin{proof}
As $\Omega$ is strictly Lorentzian, using Lemma \ref{equalAF} and the condition $\beta \neq 0$, we have $$\beta\cdot \omega^{d-1}\cdot \Omega \neq 0,$$ 
that is, $\nd_{\Omega}(\beta)\geq 1$. Then applying Lemma \ref{hallradoLorent} shows that 
\begin{equation*}
  \alpha^{d-1} \cdot \beta\cdot \Omega >0,
\end{equation*}
thus $$|\alpha+\beta|\geq |\alpha| +d  \alpha^{d-1} \cdot \beta\cdot \Omega > |\alpha|.$$

This finishes the proof.
\end{proof}

We are now ready to prove the main result of this section. It relates the proportionality of $\alpha, \beta$ to any of the equality statement (1)-(6) in the second part of Lemma \ref{logconc}.

\begin{thrm}\label{TeissierLorentzian}
Let $\Omega$ be a strictly Lorentzian class of dimension $d$, and let $\alpha,  \beta$ be big and nef classes on $\Omega$. Then the following statements are equivalent:
    \begin{enumerate}
      \item $\alpha$ and $ \beta$ are proportional;
      \item $\alpha^k$ and $ \beta^k$ are proportional for every $1\leq k\leq d-1$;
      \item $\alpha^k$ and $ \beta^k$ are proportional for some $1\leq k\leq d-1$;
      \item $\alpha^{k}\cdot \Omega$ and $\beta^{k}\cdot \Omega$ are proportional for some $1\leq k\leq d-1$;
      \item $\alpha^{d-1}$ and $ \beta^{d-1}$ are proportional;
      \item $|\alpha+\beta|^{1/d}=|\alpha|^{1/d}+|\beta|^{1/d}$.
    \end{enumerate}
\end{thrm}

\begin{proof}
It is clear that (1) implies any of the statements in (2)-(6) and (2) implies (4).

Next we show that (4) implies (6).
Assume (4) holds, i.e., for some $1\leq k\leq d-1$, there is $c>0$ such that $\alpha^{k}\cdot \Omega=c\beta^{k}\cdot \Omega$. Multiplying both sides by $\beta^{d-k}$ implies that $s_k  = c |\beta|$. Using $$s_k \geq |\alpha|^{k/d} |\beta|^{(d-k) /d},$$ 
we get
\begin{equation*}
c|\beta|^{k/d} \geq |\alpha|^{k/d}.
\end{equation*}
Similarly, multiplying the equality  $\alpha^{k}\cdot \Omega=c\beta^{k}\cdot \Omega$ by $\alpha^{d-k}$ shows that
\begin{equation*}
c|\beta|^{k/d} \leq |\alpha|^{k/d}.
\end{equation*}
Therefore, $c|\beta|^{k/d} \geq |\alpha|^{k/d}$, yielding that
\begin{equation*}
s_k=c|\beta| = |\alpha|^{k/d} |\beta|^{(d-k)/d}
\end{equation*}
which is exactly the fourth equality condition in the second part of Lemma \ref{logconc}. Combining the equivalence in Lemma \ref{logconc}, this proves $(4)\Rightarrow (6)$.

It remains to show that (6) implies (1). To this end, 
by Lemma \ref{logconc}, we have 
    \begin{equation}\label{1}
        s_l^2=s_{l-1} s_{l+1},\forall 1 \leq l \leq d-1.
    \end{equation}
Note that $s_l>0$ for any $0\leq l\leq d$. Hence, the above equality (\ref{1}) and Lemma \ref{equalAF} imply that the $(n-1,n-1)$-classes 
    $$\alpha^l\cdot\beta^{d-1-l}\cdot \Omega \text{ and } \alpha^{l-1}\cdot \beta^{d-l}\cdot \Omega$$
    are proportional for any $1\leq l\leq d-1$. As a consequence, the numbers 
    $$\frac{\gamma \cdot \alpha^l\cdot\beta^{d-1-l}\cdot \Omega}{\gamma\cdot  \alpha^{l-1}\cdot \beta^{d-l}\cdot \Omega}$$
    are independent of the choice of $\gamma \in H^{1,1}(X,\mathbb{R})$ whenever the denominator is nonzero. In particular, taking $\gamma$ to be $\alpha$ and a fixed K\"ahler class $\omega$ respectively, the equality tells that
    \begin{equation}\label{eqomeg}
    \frac{\alpha^l\cdot\beta^{d-1-l}\cdot \Omega\cdot \omega}{\alpha^{l-1}\cdot\beta^{d-l}\cdot \Omega\cdot \omega}=\frac{\alpha^{l+1}\cdot\beta^{d-1-l}\cdot \Omega}{\alpha^{l}\cdot\beta^{d-l}\cdot \Omega} =\frac{s_{l+1}}{s_l},\forall 1\leq l\leq d-1.
    \end{equation}
    Here the condition that $\alpha,\beta$ are big on $\Omega$ guarantees that the denominators are nonzero. 
    Then by (\ref{1}) and (\ref{eqomeg}), for any $ 1\leq l\leq k-2$ we have 
    $$(\alpha^l\cdot\beta^{k-1-l}\cdot \Omega\cdot \omega)^2=(\alpha^{l-1}\cdot\beta^{k-l}\cdot \Omega\cdot \omega)(\alpha^{l+1}\cdot\beta^{k-2-l}\cdot \Omega\cdot \omega).$$

    To conclude, what we proved above is that  
    $$\big(\alpha^l\cdot\beta^{k-l}\cdot\Omega\big)^2=\big(\alpha^{l-1}\cdot\beta^{k-l+1}\cdot\Omega\big)\big(\alpha^{l+1}\cdot\beta^{k-l-1}\cdot\Omega\big),\forall 1 \leq l \leq k-1$$
    imply
    $$(\alpha^l\cdot\beta^{k-1-l}\cdot \Omega\cdot \omega)^2=(\alpha^{l-1}\cdot\beta^{k-l}\cdot \Omega\cdot \omega)(\alpha^{l+1}\cdot\beta^{k-2-l}\cdot \Omega\cdot \omega),\forall 1\leq l\leq k-2.$$
    We can repeat the process and finally get
    $$(\alpha\cdot\beta\cdot \Omega\cdot \omega^{k-2})^2=(\alpha^2\cdot \Omega\cdot \omega^{k-2})(\beta^2\cdot \Omega\cdot \omega^{k-2}).$$
    This implies that $\alpha$ and $\beta$ are proportional, since $\Omega\cdot \omega^{k-2}$ has hard Lefschetz property from the assumption $\Omega$ being strictly Lorentzian.
    
    This finishes the proof.
\end{proof}

\begin{rmk}\label{intriTeisProp}
It is clear that Theorem \ref{TeissierLorentzian} extends to intrinsically strictly Lorentzian classes if we consider intrinsic proportionality, that is, proportionality in $H^*(\Omega)$.
\end{rmk}

An interesting consequence is the injectivity property on taking powers on a strictly Lorentzian class.

\begin{cor}\label{injmap}
Let $\Omega$ be a strictly Lorentzian class of dimension $d$, then for any fixed $1\leq k\leq d-1$, the map
\begin{equation*}
  \alpha \mapsto \alpha^{k}\cdot \Omega
\end{equation*}
is injective on the cone consisting of big nef classes on $\Omega$.
\end{cor} 

\begin{proof}
This follows immediately from Theorem \ref{TeissierLorentzian}.
\end{proof}

\section{Quantitative stability estimates}\label{sec stabestim}

In this section, we strengthen the proportionality characterization by establishing quantitative stability estimates.

\subsection{Quantitative deficits}\label{sec deficit}

We first introduce several deficits related to Theorem \ref{TeissierLorentzian}. Let $\alpha, \beta \in H^{1,1}(X)$ and $\Gamma$ a class of dimension 2, we introduce the notation $$\Delta(\alpha,\beta;\Gamma)=(\alpha\cdot\beta\cdot \Gamma)^2-(\alpha^2\cdot\Gamma)(\beta^2\cdot\Gamma).$$

\begin{defn}\label{defn deficit}
Let $X$ be a compact K\"ahler manifold of dimension $n$ and fix a Lorentzian class $\Omega$ of dimension $d$ ($2\leq d\leq n-2$) on $X$. Let $\alpha,\beta$ be big and nef classes on $\Omega$.
\begin{itemize}
 
  \item The Brunn-Minkowski deficit is given by
$$B(\alpha,\beta)
=\frac{|\alpha+\beta|_{\Omega}^{\frac{1}{d}}}
{|\alpha|_{\Omega}^{\frac{1}{d}}+|\beta|_{\Omega}^{\frac{1}{d}}}-1.$$
By definition, $B(\alpha,\beta) =B(\beta, \alpha)$ and $B(c\alpha,c\beta) =B(\alpha,\beta) $ for any $c>0$.

  \item The Khovanskii-Teissier deficit is defined by
$$K(\alpha,\beta)=\frac{\alpha^{d-1}\cdot\beta\cdot \Omega}{|\alpha|_{\Omega}^{\frac{d-1}{d}}|\beta|_{\Omega}^{\frac{1}{d}}}-1.$$
It is clear that $K(c_1\alpha,c_2\beta) =K(\alpha,\beta) $ for any $c_1, c_2>0$.

\item Fix a K\"ahler class $\omega$, the Alexandrov-Fenchel deficit is given by
$$A(\alpha,\beta):=A(\alpha,\beta;\omega^{d-2}\cdot \Omega)=\frac{\sqrt{\Delta(\alpha,\beta;\omega^{d-2}\cdot \Omega)}}{(\alpha\cdot \beta\cdot \omega^{d-2}\cdot \Omega)}.$$
Then we have that $A(\alpha,\beta)\leq 1$, $A(\alpha,\beta) =A(\beta, \alpha)$ and $A(c_1\alpha,c_2\beta) =A(\alpha,\beta) $ for any $c_1, c_2>0$.
\end{itemize}

\end{defn}

For $1\leq l\leq d-1$, one can also introduce the $l$-th Khovanskii-Teissier deficit as follows:
$$K_l (\alpha,\beta)=\frac{\alpha^{l}\cdot\beta^{d-l}\cdot \Omega}{|\alpha|_{\Omega}^{\frac{l}{d}}|\beta|_{\Omega}^{\frac{d-l}{d}}}-1.$$
The quantity $K(\alpha, \beta)$ above corresponds to the $(d-1)$-th Khovanskii-Teissier deficit.
Note that one can also consider $K(\beta, \alpha)$ instead. In general we have
$$K(\alpha,\beta) \neq K(\beta, \alpha).$$ 
In abbreviation, we also call the quantities introduced above BM, KT and AF deficits.

\subsubsection{Reformulation of Theorem \ref{TeissierLorentzian}}
Assume further that $\Omega$ is strictly Lorentzian, then Theorem \ref{TeissierLorentzian} can be reformulated by the deficits introduced above. 

\begin{thrm}[=Theorem \ref{TeissierLorentzian}]
Let $\Omega$ be a strictly Lorentzian class of dimension $d$ and let $\alpha,\beta$ be big and nef on $\Omega$, the following are equivalent:
\begin{itemize}
  \item $\alpha,\beta$ are proportional;
  \item $A(\alpha,\beta) = 0$;
  \item $K_l (\alpha,\beta)=0$ for some $1\leq l\leq d-1$;
  \item $K (\alpha,\beta)=0$;
  \item $B (\alpha,\beta)=0$.
\end{itemize}
\end{thrm}

Our goal is to prove quantitative stability estimates between these deficits. Taking as an example, we will show that: if $B (\alpha,\beta)$ is very small, then $A(\alpha,\beta)$ must be very small too. 

\subsection{Inradius and outradius}

\subsubsection{Bonnesen-Fenchel inequality}
The following lemma is useful to us. It is an algebro-geometric analogue of a well-known inequality in convex geometry due to Bonnesen and Fenchel \cite{BonnesenFenchel}. The proof is the same, we involve it here for readers' convenience.

\begin{lem}\label{BFineq}
    Let $\alpha, \beta, \gamma$ be nef classes and $\Omega$ a Lorentzian class of dimension 2. Then 
    \begin{align*}
        \big((\alpha\cdot\beta\cdot\Omega)(\gamma^2\cdot \Omega)-(\alpha\cdot\gamma\cdot\Omega)(\beta\cdot\gamma\cdot\Omega)\big)^2 \leq 
        \Delta(\alpha,\gamma;\Omega)\Delta(\beta,\gamma;\Omega)
    \end{align*}
\end{lem}

\begin{proof}
    Consider the polynomial $$F(\lambda,\mu)=\Delta(\alpha+\lambda \gamma,\beta +\mu \gamma;\Omega).$$
    It takes nonnegative values for all $\lambda,\mu\geq 0$ since $\Omega$ is a Lorentzian class. By definition, we can expand $F(\lambda,\mu)$ as follows:
    \begin{align*}
        F(\lambda,\mu)=\lambda^2 \Delta(\beta ,\gamma;\Omega)+\mu^2\Delta(\alpha,\gamma;\Omega) +2\lambda\mu\big((\alpha\cdot\beta\cdot\Omega)(\gamma^2\cdot \Omega)&-(\alpha\cdot\gamma\cdot\Omega)(\beta\cdot\gamma\cdot\Omega)\big)
        \\ &+\text{lower degree terms}
    \end{align*}
    Here the highest degree is 2, since terms of degree 3 and 4 are canceled.

    Then $F(\lambda ,\mu)\geq 0$ being nonnegative for all $\lambda,\mu\geq0$ implies that the quadratic term of $F(\lambda,\mu)$ is positive semi-definite. The nonnegativity of the determinant of its coefficient matrix gives the desired inequality.
    
\end{proof}

\subsubsection{Radius on the fundamental class}

We first consider the classical inradius and outradius on $X$.

\begin{defn}
    Let $\alpha,\beta$ be two psef classes. The inradius of $\alpha$ with respect to $\beta$ is 
    $$r(\alpha,\beta)=\sup\{t\geq 0:\alpha-t\beta \text{ is psef}\}.$$
    The outradius of $\alpha$ with respect to $\beta$ is defined to be 
    $$R(\alpha,\beta)=\inf\{t\geq 0:t\beta -\alpha \text{ is psef}\}.$$
\end{defn}

These two radii are related to each other via 
$$R(\alpha,\beta)=\frac{1}{r(\beta,\alpha)}.$$
The inradius or outradius is called pseudo-effective threshold in some references.

The inradius and outradius can be extracted from the intersection-number ratios by the cone duality. 

\begin{lem}\label{radiusviaduality}
    For any psef classes $\alpha$ and $\beta$, the following holds:
   \begin{align*}
     &r(\alpha,\beta)=\sup_{\delta\in \Psef^1(X)^*} \{t\geq 0:(\alpha-t\beta)\cdot \delta\geq 0 \},\\
     &R(\alpha, \beta)=\inf_{\delta\in \Psef^1(X)^*} \{t\geq 0:(t\beta -\alpha) \cdot \delta\geq 0\},
   \end{align*}
   where $\Psef^1(X)^*$ denote the dual cone of $\Psef^1(X)$.
\end{lem}

\begin{proof}
The formulas follow straightforward from the definition.
\end{proof}

As an application of Lemma \ref{BFineq}, we give a much easier proof of Teissier's result in \cite{teissier82} which is an algebro-geometric analogue of Bonnesen's inequality. 

\begin{prop}\label{Bonnesonineq}
    Let $X$ be a K\"ahler surface. Then for any big and nef classes $\alpha$ and $\beta$, the following inequality holds:
    $$\frac{\sqrt{\Delta(\alpha,\beta)}}{(\beta^2)}=\frac{\sqrt{(\alpha\cdot \beta)^2-(\alpha^2)(\beta^2)}}{(\beta^2)}\geq \frac{1}{2}(R(\alpha,\beta)-r(\alpha,\beta)).$$
\end{prop}
\begin{proof}
    By Lemma \ref{BFineq}, for any nonzero nef class $\gamma$, we have
   \begin{align*}
     \big((\beta^2)(\alpha\cdot \gamma)-(\beta\cdot\alpha)(\beta\cdot \gamma)\big)^2
     &\leq \Delta(\alpha,\beta)\Delta(\beta,\gamma)\\
     &\leq \Delta(\alpha,\beta)(\beta\cdot \gamma)^2,
   \end{align*}
   where the second inequality is trivial since by definition $\Delta(\beta,\gamma)\leq (\beta\cdot \gamma)^2$.

    Taking square root of both sides and dividing by $(\beta \cdot \gamma)(\beta^2)$ 
    yields that
    \begin{equation}\label{51}
       \frac{\sqrt{\Delta(\alpha,\beta)}}{(\beta^2)}\geq \left|\frac{(\alpha\cdot \gamma)}{(\beta \cdot \gamma)}-\frac{(\alpha\cdot\beta)}{(\beta^2)}\right|.
    \end{equation}

    Since $X$ is a surface, we have $$\Psef^1(X)^*=\Nef^1(X),$$ 
   thus
    \begin{equation}\label{52}
        r(\alpha,\beta)=\inf_{\gamma\in \Nef^1(X)\backslash 0} \frac{(\alpha\cdot \gamma)}{(\beta\cdot \gamma)},\ R(\alpha, \beta)=\sup_{\gamma\in \Nef^1(X)\backslash 0} \frac{(\alpha\cdot \gamma)}{(\beta\cdot \gamma)}.
    \end{equation}
    
    Combining (\ref{51}) and (\ref{52}) implies the desired inequality.
\end{proof}

\begin{rmk}
The inequality in Proposition \ref{Bonnesonineq} can be also deduced from the Diskant inequality proved in \cite{BFJ09}.
\end{rmk}

By the same argument, it is not hard to get a refined Khovanskii-Teissier inequality.

\begin{thrm}\label{stabKT1}
    Let $\Omega$ be a Lorentzian class of dimension 2 and $\alpha,\beta$ big nef classes on $\Omega$, then
    \begin{equation*}
        \frac{\sqrt{\Delta(\alpha,\beta;\Omega)}}{(\beta^2\cdot\Omega)}\geq \frac{1}{2}\left(\sup_{\delta\in\Nef^1(X), \beta\cdot \delta\cdot \Omega\neq 0}\frac{(\alpha\cdot\delta \cdot \Omega)}{(\beta\cdot \delta\cdot \Omega)}- \inf_{\delta\in \Nef^1(X), \beta\cdot \delta\cdot \Omega\neq 0}\frac{(\alpha\cdot\delta\cdot \Omega)}{(\beta\cdot \delta\cdot \Omega)}\right).
    \end{equation*}
   
\end{thrm}

Note that $\Omega\cdot \Nef^1(X) \subset H^{n-1,n-1}(X,\mathbb{R})$ is a full dimensional closed convex cone once $\Omega$ is strictly Lorentzian.

\begin{rmk}
When the Lorentzian property of $\Omega$ in Theorem \ref{stabKT1} has invariance via the pullback of birational morphisms, Theorem \ref{stabKT1} can be slightly improved. For example, when $\Omega=\omega^{n-2}$ for some K\"ahler class $\omega$, then
    \begin{equation*}
        \frac{\sqrt{\Delta(\alpha,\beta;\Omega)}}{(\beta^2\cdot\Omega)}\geq \frac{1}{2}\left(\sup_{\delta\in \Mov^1(X)\setminus 0}\frac{(\alpha\cdot\delta\cdot \Omega)}{(\beta\cdot \delta\cdot \Omega)}- \inf_{\delta\in \Mov^1(X)\setminus 0}\frac{(\alpha\cdot\delta\cdot \Omega)}{(\beta\cdot \delta\cdot \Omega)}\right).
    \end{equation*}
\end{rmk}

\subsubsection{Radii on a Lorentzian class}
As we are most interested in the case when $\alpha, \beta$ are nef, motivated by Lemma \ref{radiusviaduality}, on a Lorentzian class we introduce an analog of the radii for nef classes.

\begin{defn}\label{radiusonOMEGA}
Let $\Omega$ be a Lorentzian class of dimension $d$ and let $\alpha,\beta$ be nef classes on $X$. The inradius of $\alpha$ with respect to $\beta$ on $\Omega$ is defined by 
    $$r_\Omega(\alpha,\beta)=\sup\{t\geq 0:(\alpha-t\beta)\cdot B_1\cdot...\cdot B_{d-1}\cdot \Omega \geq 0, \ \forall B_i \in \Nef^1 (X)\}.$$
    The outradius of $\alpha$ with respect to $\beta$ is defined by 
   $$R_\Omega(\alpha,\beta)=\inf\{t\geq 0:(t\beta-\alpha)\cdot B_1\cdot...\cdot B_{d-1}\cdot \Omega \geq 0, \ \forall B_i \in \Nef^1 (X)\}.$$
\end{defn}

It is clear that $R_\Omega(\alpha,\beta) = \frac{1}{r_\Omega(\beta, \alpha)}$. Using the notion of radii on $\Omega$, Theorem \ref{stabKT1} can be restated as follows: for $\alpha,\beta$ big nef classes on a Lorentzian class $\Omega$ of dimension 2, 
    \begin{equation*}
        \frac{\sqrt{\Delta(\alpha,\beta;\Omega)}}{(\beta^2\cdot\Omega)}\geq \frac{1}{2}\left(R_\Omega (\alpha, \beta) - r_\Omega (\alpha, \beta)\right),
    \end{equation*}
which is analogous to Proposition \ref{Bonnesonineq}.

\begin{rmk}\label{diffradii}
For $\Omega=[X]$, we emphasis that the quantities $r_\Omega$, $R_\Omega$ are different from the radii on $X$. In general, the cone generated by the classes of the form $B_1\cdot...\cdot B_{n-1}$ with $B_i$ nef is properly contained in the dual cone $\Psef^1 (X)^*$, therefore, 
$$r_{[X]}(\alpha,\beta) \geq r(\alpha,\beta).$$  
\end{rmk}

Using Theorem \ref{stabKT1} and the definition of radii on $\Omega$, we have:

\begin{lem}\label{AFdeficitInd}
Let $\Omega$ be a Lorentzian class of dimension $d$ and let $\alpha, \beta$ be big nef classes on $\Omega$. Fix a K\"ahler class $\omega$. Then for any $1\leq l\leq d-1$, we have that
    \begin{equation*}
        \frac{\Delta(\alpha,\beta;\alpha^{l-1}\cdot \beta^{d-l-2}\cdot\omega\cdot \Omega)}{R_\Omega (\alpha,\beta)( \alpha^{l-1}\cdot \beta^{d-l}\cdot \omega \cdot \Omega )^2} \leq 2\frac{\sqrt{\Delta(\alpha,\beta; \alpha^{l-1}\cdot \beta^{d-l-1}\cdot \Omega)}}{(\alpha^{l-1}\cdot\beta^{d-l+1}\cdot\Omega)}+2\frac{\sqrt{\Delta(\alpha,\beta; \alpha^{l}\cdot \beta^{d-l-2}\cdot \Omega)}}{(\alpha^{l}\cdot\beta^{d-l}\cdot\Omega)}.
    \end{equation*}
\end{lem}

\begin{proof}

    By the very definition of $R_\Omega(\alpha,\beta)$, we have 
    \begin{equation*}
        (\alpha^l\cdot\beta^{d-l-1}\cdot \omega\cdot \Omega) \leq R_\Omega (\alpha,\beta)( \alpha^{l-1}\cdot \beta^{d-l}\cdot \omega \cdot \Omega ),
    \end{equation*}
    from which we obtain
    \begin{align*}
        \frac{\Delta(\alpha,\beta;\alpha^{l-1}\cdot \beta^{d-l-2}\cdot\omega\cdot \Omega)}{R_\Omega(\alpha,\beta)( \alpha^{l-1}\cdot \beta^{d-l}\cdot \omega \cdot \Omega )^2}
        &\leq \frac{\Delta(\alpha,\beta;\alpha^{l-1}\cdot \beta^{d-l-2}\cdot\omega\cdot \Omega)}{(\alpha^l\cdot\beta^{d-l-1}\cdot \omega\cdot \Omega)( \alpha^{l-1}\cdot \beta^{d-l}\cdot \omega \cdot \Omega )}\\
        &=
        \left|\frac{( \alpha^{l}\cdot \beta^{d-l-1}\cdot\omega\cdot \Omega)}{( \alpha^{l-1}\cdot \beta^{d-l}\cdot \omega \cdot \Omega)}-\frac{( \alpha^{l+1}\cdot \beta^{d-l-2}\cdot \omega \cdot \Omega)}{( \alpha^{l}\cdot \beta^{d-l-1}\cdot \omega\cdot \Omega)}\right| \\
        &\leq \left|\frac{( \alpha^{l}\cdot \beta^{d-l-1}\cdot\omega\cdot \Omega)}{( \alpha^{l-1}\cdot \beta^{d-l}\cdot \omega \cdot \Omega)}-\frac{( \alpha^{l+1}\cdot \beta^{d-l-1}\cdot \Omega)}{( \alpha^{l}\cdot \beta^{d-l}\cdot \Omega)}\right|\\
        &+\left|\frac{( \alpha^{l+1}\cdot \beta^{d-l-2}\cdot \omega \cdot \Omega)}{( \alpha^{l}\cdot \beta^{d-l-1}\cdot \omega\cdot \Omega)}- \frac{( \alpha^{l+1}\cdot \beta^{d-l-1}\cdot \Omega)}{( \alpha^{l}\cdot \beta^{d-l}\cdot \Omega)}\right|.
    \end{align*}

    By Theorem \ref{stabKT1}, the last two terms are bounded from above by 
    $$2\frac{\sqrt{\Delta(\alpha,\beta; \alpha^{l-1}\cdot \beta^{d-l-1}\cdot \Omega)}}{(\alpha^{l-1}\cdot\beta^{d-l+1}\cdot\Omega)}$$ and $$2\frac{\sqrt{\Delta(\alpha,\beta; \alpha^{l}\cdot \beta^{d-l-2}\cdot \Omega)}}{(\alpha^{l}\cdot\beta^{d-l}\cdot\Omega)}$$  respectively. 
    
    This finishes the proof.
\end{proof}

\subsubsection{Bounds for radii}

We will need the following estimates for the radii on a Lorentzian class. 

\begin{lem}\label{radiusbd}
Let $\Omega$ be a Lorentzian class of dimension $d$ and let $\alpha,\beta$ be big nef classes on $\Omega$. Then we have:
\begin{equation*}
  R_\Omega (\beta,\alpha) \leq 2^{d-1} \frac{ (\alpha^{d-1}\cdot \beta\cdot \Omega)}{(\alpha^d\cdot \Omega)},
\end{equation*}
and
\begin{equation*}
  r_\Omega (\alpha,\beta) \geq 2^{1-d}\frac{(\alpha^d\cdot \Omega)}{ (\alpha^{d-1}\cdot \beta\cdot \Omega)}.
\end{equation*}
\end{lem}

\begin{proof}
The estimates follow straightforward from Lemma \ref{rkt lorent}.
\end{proof}

\begin{comment}
%% TO CLEAR OUT
We want to relate $\inf_{\delta\in\Omega\cdot \Mov^1(X)}\frac{(\alpha\cdot\delta)}{(\beta\cdot \delta)}$ and $\sup_{\delta\in\Omega\cdot \Mov^1(X)}\frac{(\alpha\cdot\delta)}{(\beta\cdot \delta)}$ to some geometric quantities, for example, the inradius and outradius. This is done in the following.

\begin{lem}\label{prespsef}
    Let $H$ be a smooth hypersurface with nef fundamental class $[H]$. Then for any psef class $\alpha\in \Psef^1(X)$, its restriction $\alpha_{|H}$ to $H$ is still psef.
\end{lem}
\begin{proof}
    By the Zariski decomposition of Boucksom \cite{boucksomzariskidecomp}, we can write $\alpha=P(\alpha)+N(\alpha)$ where $P(\alpha)$ is movable and $N(\alpha)=\sum a_i [E_i]$ is a effective $\mathbb{R}$-divisor. Movable classes and prime divisors other than $H$ itself restricted to $H$ are automatically psef. Hence, it remains to show that $[H]_{|H}$ is psef. But this is straightforward by the assumption that $[H]$ is nef.
\end{proof}

\end{comment}

As a consequence of Lemma \ref{radiusbd}, we can relate the radii on $\Omega$ and the radii on $\Omega\cdot H$ as follows:

\begin{prop}\label{compradii}
Let $\Omega$ be a Lorentzian class of dimension $d$ and let $H$ be a K\"ahler class on $X$. Then for any $\alpha,\beta$ big and nef classes on $\Omega$, there is a constant $c(d)>0$ such that
    \item \begin{equation*}
        c(d) \left(1-\frac{\sqrt{\Delta(\alpha,\beta;\alpha^{d-2}\cdot \Omega)}}{(\alpha^{d-1}\cdot\beta\cdot \Omega)}\right)\leq \frac{r_\Omega(\alpha,\beta)}{r_{\Omega\cdot H}(\alpha,\beta)}\leq 1,
    \end{equation*}
    and 
    \begin{equation*}
        1\leq \frac{R_\Omega(\alpha,\beta)}{R_{\Omega\cdot H}(\alpha,\beta)} \leq \frac{1}{c(d)\left(1-\frac{\sqrt{\Delta(\alpha,\beta;\beta^{d-2}\cdot \Omega)}}{(\beta^{d-1}\cdot\alpha\cdot \Omega)}\right)}.
    \end{equation*}
    
\end{prop}
\begin{proof}
    It suffices to prove the inequality for inradii by the relation
    \begin{equation*}
        R_\bullet (\alpha,\beta)=\frac{1}{r_\bullet (\beta,\alpha)}.
    \end{equation*}

    By the definition of inradius, we always have 
    $$r_\Omega(\alpha,\beta) \leq r_{\Omega\cdot H}(\alpha,\beta).$$

For the other direction, using Lemma \ref{radiusbd}, there is some $c(d)>0$ such that
    \begin{equation}\label{571}
        r_\Omega (\alpha,\beta)\geq c(d) \frac{(\alpha^{d}\cdot \Omega)}{(\alpha^{d-1}\cdot \beta\cdot \Omega)}.
    \end{equation}
On the other hand, by the definition of inradius, we have
    \begin{equation}\label{572}
        r_{\Omega\cdot H}(\alpha,\beta)\leq \frac{(\alpha^{d-1}\cdot \Omega\cdot H)}{(\alpha^{d-2}\cdot \beta\cdot \Omega\cdot H)}.
    \end{equation}
    Together with (\ref{571}) and (\ref{572}), 
    $$\frac{r_\Omega(\alpha,\beta)}{r_{\Omega\cdot H}(\alpha,\beta)}\geq c(d)\frac{(\alpha^{d}\cdot \Omega)(\alpha^{d-2}\cdot \beta\cdot {\Omega\cdot H})}{(\alpha^{d-1}\cdot \beta \cdot \Omega)(\alpha^{d-1}\cdot {\Omega\cdot H})}.$$

    Hence, it remains to show 
    \begin{equation}\label{573}
       \frac{(\alpha^{d}\cdot \Omega)(\alpha^{d-2}\cdot \beta\cdot {\Omega\cdot H})}{(\alpha^{d-1}\cdot \beta \cdot \Omega)(\alpha^{d-1}\cdot {\Omega\cdot H})}\geq 1-\frac{\sqrt{\Delta(\alpha,\beta;\alpha^{d-2}\cdot \Omega)}}{(\alpha^{d-1}\cdot\beta\cdot \Omega)}.
    \end{equation}
    To this end, we make use of Lemma \ref{BFineq}, which tells that
    \begin{align*}
        &((\alpha^{d}\cdot \Omega)(\alpha^{d-2}\cdot \beta\cdot {\Omega\cdot H})-
        (\alpha^{d-1}\cdot \beta \cdot \Omega)(\alpha^{d-1}\cdot {\Omega\cdot H}))^2\\
        &\leq \Delta(\alpha,\beta;\alpha^{d-2}\cdot \Omega)\Delta(\alpha,H;\alpha^{d-2}\cdot \Omega).
    \end{align*}
    Taking square root and dividing both sides by $(\alpha^{d-1}\cdot \beta\cdot \Omega)(\alpha^{d-1}\cdot H\cdot \Omega)$ yield (\ref{573}).
%    \begin{equation*}
%        |\frac{(\alpha^{n})(\alpha^{n-2}\cdot \beta\cdot h)}{(\alpha^{n-1}\cdot \beta)(\alpha^{n-1}\cdot h)}-1|\leq \frac{\sqrt{\Delta(\alpha,\beta;\alpha^{n-2})}}{(\alpha^{n-1}\cdot\beta)} \frac{\sqrt{\Delta(\alpha,h;\alpha^{n-2})}}{(\alpha^{n-1}\cdot h)} \leq \frac{\sqrt{\Delta(\alpha,\beta;\alpha^{n-2})}}{(\alpha^{n-1}\cdot\beta)}
%    \end{equation*}
 
 This finishes the proof.
\end{proof}

\subsection{Comparison of deficits}\label{sec comp defi}

In this section, we study the relation between various deficits on a Lorentzian class $\Omega$ as introduced in Section \ref{sec deficit}.

For a Lorentzian class $\Omega$ of dimension $d$, recall that we adopt the notations: $$s_k = \alpha^k \cdot \beta^{d-k}\cdot\Omega,\ |\alpha|=\alpha^d \cdot \Omega. $$

We will use the convention: the notation $c(d)$ is a constant depending only on $d$, and it may represent different values line by line.

\subsubsection{KT dominates BM}
\begin{prop}\label{ktDbm}
Fix a Lorentzian class $\Omega$ of dimension $d$.
For any big and nef classes $\alpha,\beta$ on $\Omega$ with $|\alpha|\geq |\beta|$, 
$$B(\alpha, \beta)\leq \left(\frac{|\beta|}{|\alpha|}\right)^{1/d}K(\alpha,\beta).$$
\end{prop}

\begin{proof}
By scaling invariance for BM and KT deficits, after replacing $\alpha$ by $\frac{1}{|\alpha|^{1/d}}\alpha$ and $\beta$ by $\frac{1}{|\alpha|^{1/d}}\beta$, we may assume that
$$1=|\alpha|\geq |\beta|,$$
and then we need to prove $B(\alpha, \beta)\leq |\beta|^{1/d}K(\alpha,\beta)$ in this setting.

We write 
\begin{equation}\label{bmkt1}
B(\alpha, \beta)=\frac{1}{1+|\beta|^{1/d}}(|\alpha+\beta|^{1/d}-|\beta|^{1/d}-1).
\end{equation}
We estimate the upper bound of $|\alpha+\beta|$ as follows. Note that by Lemma \ref{logconc}, for any $k\geq 1$, 
$$s_{d-1} \geq s_{d-k}^{1/k}s_d ^{(k-1) /k} =s_{d-k}^{1/k} ,$$
which implies that
$$|\alpha+\beta| = 1 +\sum_{k=1}^d \binom{d}{k} s_{d-k} \leq 1 +\sum_{k=1}^d \binom{d}{k} s_{d-1}^k. $$

Therefore, 
$$|\alpha+\beta|^{1/d} \leq 1+ s_{d-1}.$$
Plugging this into the formula (\ref{bmkt1}) shows that
$$B(\alpha, \beta)\leq \frac{|\beta|^{1/d}}{1+|\beta|^{1/d}} K(\alpha, \beta)\leq  |\beta|^{1/d}K(\alpha, \beta).$$

This finishes the proof.

 \end{proof}

\subsubsection{BM dominates KT}

\begin{prop}\label{BMcontrKT}
Fix a Lorentzian class $\Omega$ of dimension $d$.
There is a constant $c(d)>0$ such that for any big and nef classes $\alpha,\beta$ on $\Omega$ with $|\alpha|\geq |\beta|$, 
$$B(\alpha, \beta)\geq c(d) \frac{ \left(\frac{|\beta|}{|\alpha|}\right)^{1/d}K(\alpha,\beta) }{1+ \left(\frac{|\beta|}{|\alpha|}\right)^{1/d}K(\alpha,\beta)}.$$
\end{prop}

\begin{proof}
By scaling invariance for BM, KT deficits, we may assume that $1=|\alpha|\geq |\beta|$. In this setting, $K(\alpha, \beta)=s_{d-1} -1$.

For $k\geq 1$, using Lemma \ref{logconc} implies $$s_{d-k} \geq s_{d-1} ^{\frac{d-k}{d-1}} s_0 ^{\frac{k-1}{d-1}},$$ 
which yields that
\begin{align*}
  |\alpha+\beta|^{1/d} & = \left(1 + \sum_{l=0} ^{d-1} \binom{d}{l+1} s_{d-l-1}\right)^{1/d} \\
 & \geq \left(1 + \sum_{l=0} ^{d-1} \binom{d}{l+1}s_{d-1} ^{\frac{d-1-l}{d-1}} s_0 ^{\frac{l}{d-1}}\right)^{1/d}\\
 &= \left(1 + \sum_{l=0} ^{d-1} \binom{d}{l+1}(1+K(\alpha, \beta)) ^{\frac{d-1-l}{d-1}} |\beta| ^{\frac{l+1}{d}}\right)^{1/d}\\
 &:=P^{1/d},
\end{align*}
where in the third line we use the definition of $K(\alpha, \beta)$ and $s_0 = |\beta|$.

Note that 
\begin{align*}
  1+|\beta|^{1/d} & = \left( (1+|\beta|^{1/d})^d\right)^{1/d} \\
 & =\left(1 + \sum_{l=0} ^{d-1} \binom{d}{l+1}|\beta| ^{\frac{l+1}{d}}\right)^{1/d}\\
 &:= Q^{1/d}.
\end{align*}

We get:
\begin{align*}
  |\alpha+\beta|^{1/d} - (1+|\beta|^{1/d})
 &\geq P^{1/d} - Q^{1/d}\\
 & = \frac{P-Q}{\sum_{j=0} ^{d-1} P^{j/d} Q^{(d-1-j) /d}}\\
 &\geq \frac{1}{2^d} \frac{P-Q}{\sum_{j=0} ^{d-1} P^{j/d}},
\end{align*}
where in the equality we use the identity $$a^m-b^m=(a-b)\sum_{j=1}^{m}a^{j-1}b^{m-j},$$ 
and in the last line we use the bound $Q=(1+|\beta|^{1/d})^d \leq 2^d$, following our assumption $|\beta|\leq 1$. For the term $P-Q$, we have:
\begin{align*}
P-Q &= \sum_{l=0} ^{d-1} \binom{d}{l+1}\left( (1+K(\alpha, \beta)) ^{\frac{d-1-l}{d-1}}-1 \right) |\beta| ^{\frac{l+1}{d}}\\
&\geq d |\beta|^{1/d}K(\alpha, \beta), \ \textrm{by taking the first summand with}\ l=0.
\end{align*}
For the upper bound for $\sum_{j=0} ^{d-1} P^{j/d} $, we estimate as follows:
\begin{align*}
\sum_{j=0} ^{d-1} P^{j/d} & \leq d P^{(d-1) /d} \\
&\leq d( 1+ \frac{d-1}{d}\sum_{l=0} ^{d-1} \binom{d}{l+1}(1+K(\alpha, \beta)) ^{\frac{d-1-l}{d-1}} |\beta| ^{\frac{l+1}{d}}\\
&\leq c(d)(1+|\beta|^{1/d}K(\alpha, \beta)), \ \textrm{by comparing with the first summand $l=0$},
\end{align*}
where in the first inequality we use $P\geq 1$, in the second inequality we use the fact that for $0<s<1$ and $x>0$, $$(1+x)^s \leq 1+sx,$$  
and in the third inequality we also use the monotonicity of the summand with respect to $l$ and the assumption $|\beta|\leq 1$.

In summary, there is some $c(d)>0$ such that
$$B(\alpha, \beta)\geq c(d) \frac{|\beta|^{1/d}K(\alpha,\beta) }{1+|\beta|^{1/d}K(\alpha,\beta)}.$$

This finishes the proof.
\end{proof}

\subsubsection{KT dominates AF}

\begin{prop}\label{KTcontrAF}
Fix a Lorentzian class $\Omega$ of dimension $d$.
There is a constant $c(d)>0$ such that for any big and nef classes $\alpha,\beta$ on $\Omega$, 
$$K(\alpha, \beta)\geq c(d) A(\alpha,\beta) ^{2^{d-2}}.$$
\end{prop}

\begin{proof}
As KT and AF deficits are 0-homogeneous with respect to all the variables, we may assume that $|\alpha|=|\beta|=1$.  Then 
\begin{equation*}
K(\alpha, \beta)= s_{d-1}-1.
\end{equation*}
In the sequel, we will frequently use $|\alpha|=|\beta|=s_0 =s_d=1$.

To prove the desired estimate, we may further assume that $K(\alpha,\beta)\leq 1$, since by definition $A(\alpha,\beta)\leq 1$. Hence, we have 
$$1\leq s_{d-1} \leq 2$$
where the lower bound follows from Lemma \ref{logconc}: $$s_{d-1}\geq |\alpha|^{d-1/d}|\beta|^{1/d}=1.$$

Using Lemma \ref{logconc} again, for any $k\leq d-1$,
\begin{equation}\label{sd1}
 s_{d-1} \geq s_k ^{1/(d-k)} s_d ^{(d-k-1) /(d-k)} =    s_k ^{1/(d-k)}.
\end{equation}
Using the identity $$a^m-b^m=(a-b)\sum_{j=1}^{m}a^{j-1}b^{m-j},$$ 
we have
\begin{equation}\label{sk}
s_k ^{1/(d-k)}-1 =\frac{s_k -1}{\sum_{j=0} ^{d-k-1} s_k ^{j/(d-k)}}.
\end{equation}
Regarding the denominator above, we first apply Lemma \ref{radiusbd} to obtain
  \begin{equation}\label{rbd1}
    R_\Omega(\beta, \alpha)\leq 2^{d-1} \frac{\alpha^{d-1}\cdot \beta\cdot \Omega}{\alpha^d\cdot \Omega}=2^{d-1} s_{d-1}\leq 2^d,
  \end{equation}
from which we get:
\begin{equation*}
 1 = s_d ^{k/d} s_0 ^{(d-k) /d}\leq s_{k} \leq  (2^d)^{d-k}.
\end{equation*}
Thus, there is a constant $c(d)$ such that
\begin{equation*}
s_k ^{1/(d-k)}-1 \geq c(d) (s_k -1),
\end{equation*}
which combining (\ref{sd1}) implies that for any $k\leq d-1$,
\begin{equation}\label{kk}
K(\alpha,\beta) \geq c(d) (s_k -1).
\end{equation}

As $s_{l} \geq 1$ for any $l\geq 0$, we have
    \begin{equation}\label{583}
       s_k -1 = \alpha^k\cdot\beta^{d-k}\cdot \Omega-1\geq\sqrt{s_{k}^2-s_{k-1}s_{k+1}} \geq\frac{\sqrt{s_{k}^2-s_{k-1}s_{k+1}}}{s_{k-1}},
    \end{equation}
Combining (\ref{kk}), this yields that
    \begin{equation}\label{lowerbd1}
     K(\alpha,\beta) \geq c(d)   \frac{\sqrt{s_{k}^2-s_{k-1}s_{k+1}}}{s_{k-1}}, \forall \ 1\leq k\leq d-1.
    \end{equation}

Applying (\ref{lowerbd1}) to $k, k+1$ and adding the inequalities, by Lemma \ref{AFdeficitInd} we obtain: for any $1\leq k\leq d-2$,
\begin{equation}\label{step1}
    K(\alpha,\beta) \geq c(d) \frac{(\alpha^k\cdot\beta^{d-k-1}\cdot \omega\cdot \Omega)^2-(\alpha^{k+1}\cdot\beta^{d-k-2}\cdot \omega\cdot \Omega)(\alpha^{k-1}\cdot\beta^{d-k}\cdot \omega\cdot \Omega)}{R_\Omega(\alpha, \beta)(\alpha^{k-1}\cdot\beta^{d-k}\cdot\omega\cdot \Omega)^2}.
    \end{equation}

Note that
\begin{align}\label{Rbd}
&R_\Omega(\alpha, \beta)\leq  2^{d-1} \frac{\beta^{d-1}\cdot \alpha\cdot \Omega}{\beta^d\cdot \Omega}
\leq 2^{d-1} R_\Omega(\beta, \alpha)^{d-1}\leq 2^{d-1} 2^{d(d-1)},
\end{align}
where in the first inequality we apply Lemma \ref{radiusbd}, in the second inequality we use the definition of the radius and $|\alpha|=|\beta|=1$, and in the last inequality we apply (\ref{rbd1}). 
Using the definition of $R_\Omega(\beta, \alpha)$ and (\ref{rbd1}),
\begin{equation}\label{rbd2}
(\alpha^{k-1}\cdot\beta^{d-k}\cdot\omega\cdot \Omega)\leq 2^d (\alpha^{k}\cdot\beta^{d-k-1}\cdot\omega\cdot \Omega) .
\end{equation}

Applying (\ref{step1}), (\ref{Rbd}) and (\ref{rbd2}) shows that for any $ 1\leq k\leq d-2$,
\begin{equation}\label{step2}
    K(\alpha,\beta)^{1/2} \geq c(d) \frac{\sqrt{(\alpha^k\cdot\beta^{d-k-1}\cdot \omega\cdot \Omega)^2-(\alpha^{k+1}\cdot\beta^{d-k-2}\cdot \omega\cdot \Omega)(\alpha^{k-1}\cdot\beta^{d-k}\cdot \omega\cdot \Omega)}}{\alpha^{k}\cdot\beta^{d-k-1}\cdot\omega\cdot \Omega}.
    \end{equation}
    
In conclusion, the process from (\ref{lowerbd1}) to (\ref{step2}) gives the lower bound of $K(\alpha,\beta)^{1/2}$ when we replace $\Omega$ by $\omega\cdot \Omega$.   
Repeating the process, we finally have
  $$K(\alpha, \beta)\geq c(d) A(\alpha,\beta) ^{2^{d-2}}.$$

    This finishes the proof. 
\end{proof}

%%%%%%%%%%%%%%%%%
%%%replacing ntod%%%

More generally, we have:

\begin{prop}
Fix a Lorentzian class $\Omega$ of dimension $d$.
For any $1\leq l\leq d-1$, there is a constant $c(d)>0$ such that for any big and nef classes $\alpha,\beta$ on $\Omega$, 
$$K_l(\alpha, \beta)\geq c(d) A(\alpha,\beta) ^{2^{d-2}}.$$
\end{prop}

\begin{proof}
By scaling invariance, Proposition \ref{KTcontrAF} and $A(\alpha,\beta)\leq 1$, it suffices to show that there is some $c(d)>0$ such that
\begin{equation*}
  K_l(\alpha, \beta)\geq c(d) K(\alpha, \beta)
\end{equation*}
under the assumption that $|\alpha|=|\beta|=1$ and $s_l \leq 2$.

Using Lemma \ref{logconc}, we have $s_l \geq s_{d-1} ^{l/(d-1)}$, implying that $s_{d-1} \leq 2^{(d-1)/l}$. We also have $s_{d-1} \geq 1$. Then we get:
\begin{align*}
  K_l (\alpha, \beta) &\geq s_{d-1} ^{l/(d-1)} -1\\
  &=(s_{d-1} ^{1/(d-1)} -1) \sum_{j=0} ^{l-1} s_{d-1} ^{j/(d-1)}\\
  &\geq l (s_{d-1} ^{1/(d-1)} -1). 
\end{align*}
Then applying the identity 
$$ s_{d-1} ^{1/(d-1)} -1 =\frac{s_{d-1}-1}{\sum_{j=0}^{d-2} s_{d-1} ^{j/(d-1)}}$$
and the upper bound for $s_{d-1}$ yield a constant $c(d)>0$ such that
\begin{equation*}
  K_l(\alpha, \beta)\geq c(d) K(\alpha, \beta),
\end{equation*}
finishing the proof.
\end{proof}

In summary, combining Propositions \ref{ktDbm}, \ref{BMcontrKT} and \ref{KTcontrAF}, we have:
\begin{thrm}\label{deficitscompa}
Let $\Omega$ be a Lorentzian class of dimension $d$.
For any big and nef classes $\alpha,\beta$ on $\Omega$, with $|\alpha|_{\Omega} \geq |\beta|_{\Omega}$, we have that: 
$$\left(\frac{|\beta|_{\Omega}}{|\alpha|_{\Omega}}\right)^{1/d}K(\alpha,\beta) \geq B(\alpha, \beta)\geq c(d) \frac{ \left(\frac{|\beta|_{\Omega}}{|\alpha|_{\Omega}}\right)^{1/d}K(\alpha,\beta) }{1+ \left(\frac{|\beta|_{\Omega}}{|\alpha|_{\Omega}}\right)^{1/d}K(\alpha,\beta)}
.$$
For any big and nef classes $\alpha,\beta$ on $\Omega$,
$$K(\alpha, \beta)\geq c(d) A(\alpha,\beta) ^{2^{d-2}}. $$
\end{thrm}

As a corollary, we show that the KT deficit dominates the ratio of the radii.

\begin{cor}\label{KTradii}
Let $\Omega$ be a Lorentzian class of dimension $d$, and let $\alpha,\beta$ be big and nef classes on $\Omega$. Denote $\Gamma = \omega^{d-2}\cdot \Omega$. Then there is a constant $c(d)>0$ such that 
    \begin{equation*}
       K(\alpha,\beta)^{\frac{1}{2^{d-2}}} \geq c(d) \left( 1-\frac{r_\Gamma(\alpha,\beta)}{R_\Gamma(\alpha,\beta)}  \right).
    \end{equation*}
\end{cor}

\begin{proof}
By Theorem \ref{stabKT1} and Theorem \ref{deficitscompa},
\begin{align*}
K(\alpha,\beta)^{\frac{1}{2^{d-2}}} &\geq c(d) A(\alpha,\beta)\\
  & = c(d) \frac{\Delta(\alpha,\beta; \omega^{d-2}\cdot \Omega)^{1/2}}{\alpha \cdot \beta\cdot  \omega^{d-2}\cdot \Omega}\\
  & \geq c(d) \frac{\Delta(\alpha,\beta; \omega^{d-2}\cdot \Omega)^{1/2}}{R_\Gamma(\alpha,\beta)( \beta^2\cdot  \omega^{d-2}\cdot \Omega)}\\
  & \geq c(d) \frac{1}{R_\Gamma(\alpha,\beta)} \cdot \frac{1}{2}(R_\Gamma(\alpha,\beta) -r_\Gamma(\alpha,\beta)).
\end{align*}

The finishes the proof.
\end{proof}

\subsubsection{AF dominates BM/KT}
The above discussions show that the AF deficit is the smallest one.
We next show that, in some sense, the converse of Theorem \ref{deficitscompa} is also true. That is, the AF deficit can also control the BM/KT deficit but the estimates involve the radii of the classes.

Such an estimate depends on an algebro-geometric analogue of \cite{Schneider90stabAF}.

\begin{lem}\label{schneider}
Let $\Omega$ be a Lorentzian class of dimension $d$. Let $\alpha, \beta ,\gamma_3,\overline{\gamma_3},\gamma_4,...,\gamma_d$ be big and nef classes on $\Omega$
 such that $\gamma_3':=\gamma_3-\overline{\gamma_3}$ is nef. Denote 
 $$\Gamma=\gamma_3 \cdot \gamma _4\cdots\gamma_d \cdot \Omega \text{ and } \overline{\Gamma}=\overline{\gamma_3}\cdot \gamma_4\cdots\gamma_d \cdot \Omega,$$
    then 
    \begin{equation*}
        \frac{\Delta(\alpha,\beta;\overline{\Gamma})}{(\alpha\cdot\beta\cdot \overline{\Gamma})^2} \leq 
        4 \min\left\{\frac{R_\Omega(\alpha,\gamma_3)}{r_\Omega(\alpha,\gamma_3)},\frac{R_\Omega(\beta,\gamma_3)}
        {r_\Omega(\beta,\gamma_3)}\right\}R_\Omega(\gamma_3,\overline{\gamma_3})
        \frac{\sqrt{\Delta(\alpha,\beta;\Gamma)}}{(\alpha\cdot\beta\cdot \Gamma)}.
    \end{equation*}
    
\end{lem}

Roughly speaking, this result tells that decreasing the positivity of the class $\Gamma$ will decrease the quantity $\Delta(\alpha,\beta;\Gamma)/(\alpha\cdot\beta\cdot \Gamma)^2$. 

%question: 
%\red{Can we replace $\Omega=\gamma_3 \cdot \gamma _4\cdots\gamma_n$ by $\Omega=\langle\gamma_3 \cdot \gamma _4\cdots\gamma_n\rangle$?}

\begin{proof}
Rewriting the estimate in Lemma \ref{BFineq} for the classes $\alpha,\beta,\gamma_3,\overline
    {\Gamma}$, we get 
    \begin{equation}\label{521}
        \Delta(\alpha,\beta;\overline{\Gamma})\leq \frac{(\alpha\cdot\gamma_3\cdot \overline{\Gamma})(\beta\cdot\gamma_3\cdot \overline{\Gamma})}{(\gamma_3^2\cdot \overline{\Gamma})}
        \big(2(\alpha\cdot\beta\cdot\overline{\Gamma})-(\alpha^2\cdot \overline{\Gamma})\frac{(\beta\cdot\gamma_3\cdot \overline{\Gamma})}{(\alpha\cdot\gamma_3\cdot \overline{\Gamma})}-(\beta^2 \cdot \overline{\Gamma})\frac{(\alpha\cdot\gamma_3\cdot \overline{\Gamma})}{(\beta\cdot\gamma_3\cdot \overline{\Gamma})}\big).
    \end{equation}
    Denote
    \begin{align*}
      \Lambda(\xi,\eta)=2(\alpha\cdot\beta\cdot \xi\cdot \gamma_4\cdots\gamma_d \cdot \Omega)
      &-(\alpha^2\cdot \xi\cdot \gamma_4\cdots\gamma_d \cdot \Omega)
    \frac{(\beta\cdot\eta\cdot \gamma_3 \cdot \gamma_4\cdots\gamma_d \cdot \Omega)}{(\alpha\cdot\eta\cdot \gamma_3\cdot \gamma_4\cdots\gamma_d \cdot \Omega)}\\
    &-
    (\beta^2 \cdot \xi \cdot\gamma_4\cdots\gamma_d \cdot \Omega)
    \frac{(\alpha\cdot\eta\cdot \gamma_3 \cdot \gamma_4\cdots\gamma_d \cdot \Omega)}{(\beta\cdot\eta\cdot \gamma_3 \cdot \gamma_4\cdots\gamma_d \cdot \Omega)}.
    \end{align*}
    Then (\ref{521}) can be written as
    \begin{equation}\label{523}
        \frac{\Delta(\alpha,\beta;\overline{\Gamma})}{(\alpha\cdot\beta\cdot \overline{\Gamma})^2}\leq 
        \frac{(\alpha\cdot\gamma_3\cdot \overline{\Gamma})(\beta\cdot\gamma_3\cdot \overline{\Gamma})}{(\gamma_3^2\cdot \overline{\Gamma})(\alpha\cdot\beta\cdot \overline{\Gamma})} \frac{1}{(\alpha\cdot\beta\cdot \overline{\Gamma})} \Lambda(\overline{\gamma_3},\overline{\gamma_3}).
    \end{equation}

    \textbf{Claim:}
    \begin{equation}\label{522}
        \Lambda (\overline{\gamma_3},\overline{\gamma_3})\leq 
        4 \sqrt{\Delta(\alpha,\beta;\Gamma)}.
    \end{equation}

    On the other hand, by definition of the radius (on $\Omega$), we have
    \begin{align*}
     ( \alpha\cdot\gamma_3\cdot \overline{\Gamma}) &\leq R_\Omega(\alpha, \gamma_3) (\gamma_3 ^2 \cdot \overline{\Gamma}),\\
      (\alpha\cdot\beta \cdot \overline{\Gamma} )& \geq r_\Omega(\alpha, \gamma_3) (\gamma_3\cdot\beta \cdot \overline{\Gamma}).
    \end{align*}
Using similar inequalities for $\beta, \gamma_3$ and the inequality $$(\alpha\cdot \beta\cdot \Gamma)\leq R_\Omega(\gamma_3, \overline{\gamma_3}) (\alpha\cdot \beta\cdot \overline{\Gamma})$$
implies:
    \begin{equation}\label{525}
    \frac{(\alpha\cdot\gamma_3\cdot \overline{\Gamma})(\beta\cdot\gamma_3\cdot \overline{\Gamma})}{(\gamma_3^2\cdot \overline{\Gamma})(\alpha\cdot\beta\cdot \overline{\Gamma})} \frac{1}{(\alpha\cdot\beta\cdot \overline{\Gamma})} \leq 
    \min\left\{\frac{R_\Omega(\alpha,\gamma_3)}{r_\Omega(\alpha,\gamma_3)},\frac{R_\Omega(\beta,\gamma_3)}
    {r_\Omega(\beta,\gamma_3)}\right\}R_\Omega(\gamma_3,\overline{\gamma_3})\frac{1}{(\alpha\cdot\beta\cdot \Gamma)}.
    \end{equation}
    Combining (\ref{523}), (\ref{525}) and the claimed (\ref{522}), we get
    \begin{equation*}
        \frac{\Delta(\alpha,\beta;\overline{\Gamma})}{(\alpha\cdot\beta\cdot \overline{\Gamma})^2} \leq 
        4 \min\left\{\frac{R_\Omega(\alpha,\gamma_3)}{r_\Omega(\alpha,\gamma_3)},\frac{R_\Omega(\beta,\gamma_3)}
        {r_\Omega(\beta,\gamma_3)}\right\}R_\Omega(\gamma_3,\overline{\gamma_3})\frac{\sqrt{\Delta(\alpha,\beta;\Gamma)}}{(\alpha\cdot\beta\cdot \Gamma)}.
    \end{equation*}

Next we prove the claim.

Denote $\Gamma' =\gamma'_3 \cdot \gamma _4\cdots\gamma_d \cdot \Omega$, then $\Gamma=\Gamma'+\overline{\Gamma}$.

Note that, for any nef class $\eta$, we have
\begin{align}\label{5c1}
    \Lambda(\gamma_3,\eta)&= 2(\alpha\cdot\beta\cdot \Gamma)-
    (\alpha^2\cdot \Gamma)\frac{(\eta \cdot\beta\cdot \Gamma)}{(\eta \cdot\alpha \cdot \Gamma)}-
    (\beta^2 \cdot \Gamma)\frac{(\eta \cdot\alpha\cdot \Gamma)}{(\eta \cdot\beta\cdot \Gamma)} \nonumber \\
    &= \frac{(\beta\cdot \eta\cdot \Gamma)}{(\beta^2\cdot\Gamma)(\alpha\cdot \eta\cdot \Gamma)}\Delta(\alpha,\beta;\Gamma)-
    \frac{(\beta^2\cdot\Gamma)(\beta\cdot \eta\cdot \Gamma)}{(\alpha\cdot \eta\cdot \Gamma)}
    (\frac{(\alpha\cdot \beta\cdot\Gamma)}{(\beta^2\cdot \Gamma)}-\frac{(\alpha\cdot \eta \cdot\Gamma)}{(\beta\cdot \eta \cdot \Gamma)})^2
\end{align}
where the second equality can be checked directly. Taking $\eta= \beta$ yields
\begin{equation}\label{5c2}
    \Lambda(\gamma_3,\beta )=\frac{\Delta(\alpha,\beta;\Gamma)}{(\alpha\cdot\beta\cdot \Gamma)}\leq \sqrt{\Delta(\alpha,\beta;\Gamma)}.
\end{equation}

Using (\ref{521}) (with $\overline{\Gamma}$ replaced by $\Gamma$), we have 
\begin{equation}\label{5c3}
    \Lambda(\gamma_3,\eta)\geq 0.   %    ,\Lambda(\eta,\eta)
\end{equation}
This and (\ref{5c1}) show that 
\begin{equation}\label{5c3}
    \left|\frac{(\alpha\cdot \beta\cdot\Gamma)}{(\beta^2\cdot \Gamma)}-\frac{(\alpha\cdot \eta \cdot\Gamma)}{(\beta\cdot \eta \cdot \Gamma)}\right|\leq
    \frac{\sqrt{\Delta(\alpha,\beta;\Gamma)}}{(\beta^2\cdot\Gamma)}
\end{equation}
and as a consequence,
\begin{align}\label{5c4}
    \left|\frac{(\beta^2\cdot \Gamma)}{(\alpha\cdot \beta\cdot\Gamma)}-\frac{(\beta\cdot \eta \cdot \Gamma)}{(\alpha\cdot \eta \cdot\Gamma)}\right|
    &=\frac{(\beta^2\cdot \Gamma)(\beta\cdot \eta \cdot \Gamma)}{(\alpha\cdot \beta\cdot\Gamma)(\alpha\cdot \eta \cdot\Gamma)}\left|\frac{(\alpha\cdot \beta\cdot\Gamma)}{(\beta^2\cdot \Gamma)}-\frac{(\alpha\cdot \eta \cdot\Gamma)}{(\beta\cdot \eta \cdot \Gamma)}\right|\\
    &\leq\frac{(\beta\cdot \eta \cdot\Gamma)\sqrt{\Delta(\alpha,\beta;\Gamma)}}{(\alpha \cdot \beta\cdot\Gamma)(\alpha \cdot \eta\cdot\Gamma)}. \nonumber
\end{align}
Using the definition of $\Lambda$, (\ref{5c3}) and (\ref{5c4}), we have
\begin{align*}
    |\Lambda(\gamma_3',\gamma_3')-\Lambda(\gamma_3',\beta)|&\leq 
    (\alpha^2\cdot \Gamma')\left|\frac{(\beta\cdot\gamma_3'\cdot \Gamma)}{(\alpha\cdot \gamma_3'\cdot\Gamma)}-\frac{(\beta^2\cdot \Gamma)}{(\alpha \cdot \beta\cdot \Gamma)}\right|+ (\beta^2\cdot \Gamma')\left|\frac{(\alpha\cdot \gamma_3'\cdot\Gamma)}{(\beta\cdot\gamma_3'\cdot \Gamma)}-\frac{(\alpha \cdot \beta\cdot \Gamma)}{(\beta^2\cdot \Gamma)}\right|
    \\
    &\leq (\frac{(\alpha^2\cdot \Gamma')(\beta\cdot\gamma_3'\cdot \Gamma)}{(\alpha\cdot \beta\cdot \Gamma)(\alpha\cdot \gamma_3'\cdot\Gamma)}+\frac{(\beta^2\cdot \Gamma')}{(\beta^2\cdot \Gamma)})\sqrt{\Delta(\alpha,\beta;\Gamma)}
\end{align*}

The rKT property for $\Gamma$ tells that 
$$2(\alpha\cdot \beta\cdot \Gamma)(\alpha\cdot \gamma_3'\cdot\Gamma) \geq  (\alpha^2\cdot \Gamma)(\beta\cdot\gamma_3'\cdot \Gamma),$$
which yields 
$$\frac{(\alpha^2\cdot \Gamma')(\beta\cdot\gamma_3'\cdot \Gamma)}{(\alpha\cdot \beta\cdot \Gamma)(\alpha\cdot \gamma_3'\cdot\Gamma)}\leq 2 \frac{(\alpha^2\cdot\Gamma')}{(\alpha^2\cdot \Gamma)}.$$
Hence, we get 
\begin{equation}\label{5c5}
    |\Lambda(\gamma_3',\gamma_3')-\Lambda(\gamma_3',\beta)|\leq \left(2 \frac{(\alpha^2\cdot\Gamma')}{(\alpha^2\cdot \Gamma)}+\frac{(\beta^2\cdot \Gamma')}{(\beta^2\cdot \Gamma)}\right)\sqrt{\Delta(\alpha,\beta;\Gamma)}.
\end{equation}
Similarly, we have 
\begin{equation}\label{5c6}
    |\Lambda(\overline{\gamma_3},\overline{\gamma_3})-\Lambda(\overline{\gamma_3},\beta)|\leq \left(2 \frac{(\alpha^2\cdot\overline{\Gamma})}{(\alpha^2\cdot \Gamma)}+\frac{(\beta^2\cdot \overline{\Gamma})}{(\beta^2\cdot \Gamma)}\right)\sqrt{\Delta(\alpha,\beta;\Gamma)}.
\end{equation}
With the inequality
\begin{equation*}
    \Lambda(\overline{\gamma_3},\overline{\gamma_3})\leq
    \Lambda(\gamma_3,\beta)+
    |\Lambda(\gamma_3',\gamma_3')-\Lambda(\gamma_3',\beta)|+
    |\Lambda(\overline{\gamma_3},\overline{\gamma_3})-\Lambda(\overline{\gamma_3},\beta)|,
\end{equation*}
and (\ref{5c2}), (\ref{5c5}), (\ref{5c6}) applying to the three terms on the right side respectively and $\Gamma = \Gamma' + \overline{\Gamma}$, we finally get
\begin{equation*}
    \Lambda(\overline{\gamma_3},\overline{\gamma_3})\leq 4\sqrt{\Delta(\alpha,\beta;\Gamma)}.
\end{equation*}

This finishes the proof.
\end{proof}

%\hfill\qedsymbol{}

Now we can prove that the AF deficit dominates the BM/KT deficit.

\begin{thrm}
Let $\Omega$ be a Lorentzian class of dimension $d$, and let $\alpha, \beta$ be big and nef classes on $\Omega$. Fix a K\"ahler class $\omega$. We write $$\delta(\alpha,\omega)=\frac{R_\Omega(\alpha,\omega)}{r_\Omega(\alpha,\omega)^2}\inf\{\lambda>0:\lambda\omega-\alpha \text{ is nef }\},$$
and set 
\begin{equation*}
 C= \max\{R_\Omega(\alpha,\beta),R_\Omega(\beta,\alpha)\}^{\frac{d}{2}+2} \max\{\delta (\alpha,\omega),\delta(\beta,\omega)\}^2
\end{equation*}
    then we have:
 \begin{equation*}
        B(\alpha,\beta)
        \leq 4Cd^3 A(\alpha,\beta;\omega^{d-2}\cdot \Omega)^{\frac{1}{2^{d-3}}},
    \end{equation*}
and for any $1\leq l\leq d-1$,
    \begin{equation*}
        K_l(\alpha,\beta)
        \leq 4Cd^2 A(\alpha,\beta;\omega^{d-2} \cdot\Omega)^{\frac{1}{2^{d-3}}}.
    \end{equation*}

\end{thrm}

\begin{proof}

We first prove that the AF deficit dominates the BM deficit.

Recall that we adopt the notation $s_k = \alpha^k \cdot \beta^{d-k}\cdot \Omega$.

    By the identity $$a^n-b^n=(a-b)\sum_{k=1}^na^{k-1}b^{n-k},$$ we have
    \begin{align*}
        B(\alpha,\beta)
        &=\frac{|\alpha+\beta|^{\frac{1}{d}}}{|\alpha|^{\frac{1}{d}}+|\beta|^{\frac{1}{d}}}-1
        =\frac{|\alpha+\beta|^{\frac{1}{d}}-(s_d^{\frac{1}{d}}+s_0^{\frac{1}{d}})}
        {s_d^{\frac{1}{d}}+s_0^{\frac{1}{d}}}\\
        &=\frac{|\alpha+\beta|-(s_d^{\frac{1}{d}}+s_0^{\frac{1}{d}})^{d}}
        {(s_d^{\frac{1}{d}}+s_0^{\frac{1}{d}})
        \sum_{k=1}^d|\alpha+\beta|^{\frac{k-1}{d}}(s_d^{\frac{1}{d}}+s_0^{\frac{1}{d}})^{d-k}}\\
        &\leq \frac{|\alpha+\beta|-(s_d^{\frac{1}{d}}+s_0^{\frac{1}{d}})^{d}}
        {(s_d^{\frac{1}{d}}+s_0^{\frac{1}{d}})^d}\\
        &=\sum_{k=0}^d\binom{d}{k}\frac{s_k-s_d^{\frac{k}{d}}s_0^{\frac{d-k}{d}}}
        {(s_d^{\frac{1}{d}}+s_0^{\frac{1}{d}})^d}\\
         &\leq \sum_{k=0}^d \left(\frac{s_k}{s_d^{^{\frac{k}{d}}}s_0^{\frac{d-k}{d}}}-1\right),
    \end{align*}
    where in the estimate from the second line to the third line we only keep the summand given by $k=1$ in the denominator, and in the last line we simply use 
    $$(s_d^{\frac{1}{d}}+s_0^{\frac{1}{d}})^d\geq \binom{d}{k} s_d^{\frac{k}{d}}s_0^{\frac{d-k}{d}}.$$
    Using the inequality 
    \begin{equation*}
        x-1 \leq x \log x \ \text{for}\ x\geq 1,
    \end{equation*}
    we obtain
    \begin{equation*}
        B(\alpha,\beta)\leq \sum_{k=0}^d \frac{s_k}{s_d^{^{\frac{k}{d}}}s_0^{\frac{d-k}{d}}} \log\left(\frac{s_k}{s_d^{^{\frac{k}{d}}}s_0^{\frac{d-k}{d}}}\right).
    \end{equation*}

    By the definition of outradii, we have 
    \begin{align*}
        \frac{s_k}{s_d^{^{\frac{k}{d}}}s_0^{\frac{d-k}{d}}}
        &\leq R_\Omega(\beta,\alpha)^{\frac{k(d-k)}{d}}R_\Omega(\alpha,\beta)^{\frac{k(d-k)}{d}}\\
        &\leq \max\{R_\Omega(\alpha,\beta),R_\Omega(\beta,\alpha)\}^{\frac{d}{2}},
    \end{align*}
    where in the last line we apply $4k(d-k)\leq d^2$.

    Hence, we have 
    \begin{equation*}
        B(\alpha,\beta)\leq \max\{R_\Omega(\alpha,\beta),R_\Omega(\beta,\alpha)\}^{\frac{d}{2}} \sum_{k=1}^{d-1}(\log s_k-\frac{k}{d}\log s_d-\frac{d-k}{d}\log s_0).
    \end{equation*}
    The summands can be estimated as follows:
    \begin{align*}
        &\log s_k-\frac{k}{d}\log s_d-\frac{d-k}{d}\log s_0\\
        & =\frac{k(d-k)}{d}\left(\frac{\log s_k - \log s_0}{k} -\frac{\log s_d -\log s_k }{d-k}\right)\\
        &\leq \frac{k(d-k)}{d}\left((\log s_1 -\log s_0) -(\log s_d -\log s_{d-1})\right)\\
        &= \frac{k(d-k)}{d}\sum_{k=1}^{d-1}((\log s_k -\log s_{k-1}) -(\log s_{k+1} -\log s_{k}))\\
        &\leq \frac{d^2}{4}\max\{(\log s_k -\log s_{k-1}) -(\log s_{k+1} -\log s_{k}) :k=1,...,d-1\}\\
        &=\frac{d^2}{4}\max\left\{\log\left(\frac{s_k^2}{s_{k-1}s_{k+1}}\right) :k=1,...,d-1\right\}
    \end{align*}
    where for the first inequality, we use the log-concavity of the sequence $\{s_k\}$.
    For each $k$, applying the inequality $$\log x\leq x-1\ \text{for}\ x\geq 1,$$ 
    we have 
    \begin{align*}
        \log\left(\frac{s_k^2}{s_{k-1}s_{k+1}}\right)&\leq \frac{s_k^2}{s_{k-1}s_{k+1}}-1 \\ &=\frac{s_k^2}{s_{k-1}s_{k+1}}\left(1- \frac{s_{k-1}s_{k+1}}{s_k^2}\right) \\
        &\leq \max\{R_\Omega(\alpha,\beta),R_\Omega(\beta,\alpha)\}^{2} A(\alpha,\beta;\alpha^{k-1}\cdot\beta^{d-k-1}\cdot \Omega)^2,
    \end{align*}
    where the last line follows from the definitions of $R_\Omega$ and $A(-,-;\alpha^{k-1}\cdot\beta^{d-k-1}\cdot \Omega)$.

    To conclude, we obtain 
    \begin{equation*}
        B(\alpha,\beta)\leq \frac{d^3}{4}\max\{R_\Omega(\alpha,\beta),R_\Omega(\beta,\alpha)\}^{\frac{d}{2}+2} \max\{A(\alpha,\beta;\alpha^{k-1}\cdot\beta^{d-k-1}\cdot \Omega)^2:k=1,...,d-1\}.
    \end{equation*}

    Hence, it remains to show that, for each $1\leq k\leq d-1$, we have 
    \begin{equation*}
        A(\alpha,\beta;\alpha^{k-1}\cdot\beta^{d-1-k}\cdot \Omega)\leq 4\max\{\delta(\alpha,\omega),\delta(\beta,\omega)\} A(\alpha,\beta;\omega^{d-2}\cdot \Omega)^\frac{1}{2^{d-2}}.
    \end{equation*}
    To this end, we apply Lemma \ref{schneider} as follows: for any $k\geq 2$ and for any $\lambda >0$ such that $\lambda \omega-\alpha $ is nef, we have
    \begin{equation*}
        A(\alpha,\beta;\alpha^{k-1}\cdot\beta^{d-k-1}\cdot \Omega)^2\leq 4 \frac{R_\Omega(\alpha,\lambda \omega )}{r_\Omega(\alpha,\lambda\omega)} R_\Omega(\lambda \omega, \alpha) A(\alpha,\beta;\lambda \omega \cdot \alpha^{k-2}\cdot \beta^{d-k-1}\cdot \Omega).
    \end{equation*}
    Note that the inradius and outradius are both 1-homogeneous with respect to the first variable and $(-1)$-homogeneous with respect to the second variable, and the AF deficit $A(\alpha,\beta;\Gamma)$ is $0$-homogeneous with respect to the reference class $\Gamma$. Thus, we have 
    \begin{equation*}
    A(\alpha,\beta;\alpha^{k-1}\cdot\beta^{d-k-1}\cdot \Omega)^2\leq 4\lambda \frac{R_\Omega(\alpha, \omega )}{r_\Omega(\alpha,\omega)} R_\Omega( \omega, \alpha) A(\alpha,\beta; \omega \cdot \alpha^{k-2}\cdot \beta^{d-k-1}\cdot \Omega).
    \end{equation*}
    Since the inequality holds for any $\lambda$ satisfying $\lambda \omega -\alpha $ is nef, using $r_\Omega(\alpha,\omega) =R_\Omega(\omega,\alpha)^{-1}$ we obtain:
    \begin{equation*}
        A(\alpha,\beta;\alpha^{k-1}\cdot\beta^{d-k-1}\cdot \Omega)^2\leq 4 \delta(\alpha,\omega) A(\alpha,\beta;\omega\cdot \alpha^{k-2}\cdot\beta^{d-k-1}\cdot \Omega).
    \end{equation*}
    Similar inequality holds if we replace one $\beta$ in the reference $\alpha^{k-1}\cdot\beta^{d-k-1}\cdot \Omega$ by $\omega$:
    \begin{equation*}
        A(\alpha,\beta;\alpha^{k-1}\cdot\beta^{d-k-1}\cdot \Omega)^2\leq 4 \delta(\beta,\omega) A(\alpha,\beta;\omega\cdot \alpha^{k-1}\cdot\beta^{d-k-2}\cdot \Omega).
    \end{equation*}
    We can repeat this process, and finally get
    \begin{align*}
        A(\alpha,\beta;\alpha^{k-1}\cdot\beta^{d-k-1}\cdot \Omega) &\leq (4\max\{\delta(\alpha,\omega),\delta (\beta,\omega)\})^{\frac{2^{d-2}-1}{2^{d-2}}}A(\alpha,\beta;\omega^{d-2}\cdot \Omega)^{\frac{1}{2^{d-2}}}\\
        &\leq 4\max\{\delta (\alpha,\omega),\delta(\beta,\omega)\}A(\alpha,\beta;\omega^{d-2}\cdot \Omega)^{\frac{1}{2^{d-2}}}
    \end{align*}
    where for the second inequality, we use $\delta (\alpha,\omega),\delta(\beta,\omega)\geq 1$.

    Checking the above argument also shows that the AF deficit dominates the KT deficit in a similar way.

    This finishes the proof.
\end{proof}

%\begin{thrm}
%    Let $\Omega$ be a Lorentzian class of dimension $d$, and let $\alpha, \beta$ be big and nef classes on $\Omega$. Fix a K\"ahler class $\omega$. Then for any $1\leq l\leq n-1$,
%    \begin{equation*}
%        K_l(\alpha,\beta)
%        \leq 4d^2\max\{R_\Omega(\alpha,\beta),R_\Omega(\beta,\alpha)\}^{\frac{d}{2}+2} \max\{\delta (\alpha,\omega),\delta(\beta,\omega)\}^2 A(\alpha,\beta;\omega^{d-2} \cdot\Omega)^{\frac{1}{2^{d-3}}}.
%    \end{equation*}
%\end{thrm}

\subsection{FMP type estimates}

This section is motivated by the sharp quantitative isoperimetric and the refined Brunn-Minkowski inequality in convex geometry, established by Fusco-Maggi-Pratelli and Figalli-Maggi-Pratelli (see \cite{FMP08Annals, figalli09, Figalli10inventiones}). We call the stability estimates in this section FMP type estimates.

\subsubsection{Relative asymmetry index}
The key quantity is the relative asymmetry index, defined as follows. 
For any two convex bodies $A$ and $B$ of dimension $n$, their relative asymmetry index (differs from the usual definition by a multiple $\frac{1}{2}$) is defined as
\begin{equation*}
    F(A,B)=\inf_{x\in\mathbb{R}^n}\frac{|(x+rA)\Delta B|}{2|B|}.
\end{equation*}
Here $$(x+rA)\Delta B= (x+rA)\backslash B \cup B\backslash (x+rA)$$ 
is the symmetric difference and $r=\left(\frac{|B|}{|A|}\right)^{\frac{1}{n}}$. Note that $|rA|=|B|$ and thus, $F(A,B)$ can be rewritten as
\begin{equation*}
    F(A,B)=\inf_{x\in \mathbb{R}^n}\left(1-\frac{|(x+rA)\cap B|}{|B|}\right)=1-\sup_{x\in \mathbb{R}^n}\frac{|(x+rA)\cap B|}{|B|}.
\end{equation*}

To motivate the algebro-geometric construction, 
we claim that $$\sup_{x\in \mathbb{R}^n}\frac{|(x+rA)\cap B|}{|B|}= \sup_{K \leq rA, K\leq B}\frac{|K|}{|B|}$$
where the notation $K\leq L$ means that $K$ is contained in $L$ up to a translation. 
First, $(x+rA)\cap B$ is clearly contained in both $rA$ and $B$ up to translations which implies that 
$$\sup_{x\in \mathbb{R}^n}\frac{|(x+rA)\cap B|}{|B|} \leq \sup_{K \leq rA, K \leq B}\frac{|K|}{|B|}.$$
On the other hand, for any convex body $K$ satisfying $x+K\subset rA$ and $y+K \subset B$, we have  that
$$y+K \subset (y-x+rA) \cap B$$
which gives the converse inequality. To conclude, we obtain
\begin{equation*}
    F(A,B)=1- \sup_{K \leq rA,K \leq B}\frac{|K|}{|B|}.
\end{equation*}

It is clear that $F(A,B)=0$ iff $A$ and $B$ are homothetic, that is, there exists $\lambda>0$ and $x\in \mathbb{R}^n$ such that 
$x+\lambda A=B$. Hence, $F(A,B)$ is a natural function reflecting how far of $A$ and $B$ from being homothetic.

Figalli, Maggi and Pratelli proved the following remarkable result: let $K, L \subset \mathbb{R}^n$ be convex bodies, then there exists a constant $c(n)$ depending only on the dimension $n$ (and having polynomial growth) such that
\begin{equation*}
        F(K, L)\leq  c(n) \sqrt{\sigma(K,L)B(K, L)} 
\end{equation*}
where $$\sigma(K,L)=\max\left\{\left(\frac{|K|}{|L|}\right)^{\frac{1}{n}},
\left(\frac{|L|}{|K|}\right)^{\frac{1}{n}}\right\} $$ 
is the relative size of $K, L$. Furthermore, the exponent $1/2$ on the BM deficit is optimal.

\subsubsection{Algebro-geometric analog and conjecture}
The above discussion motivates the following definition:

\begin{defn}\label{asymIndex}
Let $X$ be a compact K\"ahler manifold of dimension $n$. Let $\alpha$ and $\beta$ be two big movable classes on $X$. Denote $r=\left(\frac{|\beta|}{|\alpha|}\right)^{\frac{1}{n}}$. Their relative asymmetry index is defined as
    \begin{equation*}
        F(\alpha,\beta)=1-\sup_{\gamma\leq r\alpha,\gamma\leq\beta}\frac{|\gamma|}{|\beta|}
    \end{equation*}
    where $\gamma$ ranges over all big movable classes such that $\gamma\leq r\alpha, \gamma\leq\beta$. Here $A\leq B$ means that $B-A$ is psef.
\end{defn}

It is easy to see that $F(\alpha,\beta) = F(\beta, \alpha)$ and $F(a\alpha,b\beta)=F(\alpha,\beta)$ for any $a, b>0$.

Similarly to the convexity setting, $F(\alpha,\beta)$ measures the failure of $\alpha$ and $\beta$ being proportional.

\begin{lem}\label{movableusingLX}
Let $X$ be a projective manifold of dimension $n$.
For $\alpha,\beta$ big movable on $X$, $F(\alpha,\beta)=0$ iff $\alpha$ and $\beta$ are proportional.
\end{lem}

\begin{proof}
It is clear that $F(\alpha,\beta)=0$ if $\alpha$ and $\beta$ are proportional. 

In the other direction, take a sequence $\gamma_k$ satisfying $\gamma_k \leq r\alpha, \gamma_k\leq\beta$ such that 
$$|\gamma_k|\rightarrow \sup_{\gamma\leq r\alpha,\gamma\leq\beta} |\gamma|.$$
Note that such a sequence lies in a compact subset of $\Mov^1 (X)$, by taking a subsequence, we may assume that $\gamma_k \rightarrow M$ for some big movable class $M$. The class $M$ satisfies that $M\leq r\alpha,M\leq\beta$ and by $F(\alpha,\beta)=0$, $|M|=|\beta|$. Then by the monotonicity of positive product and Khovanskii-Teissier inequality, we must have
\begin{equation*}
  \langle M^{n-1}\rangle\cdot \beta = |M|^{\frac{n-1}{n}}|\beta|^{\frac{1}{n}}.
\end{equation*}
Using \cite[Section 3]{lehmxiaoPosCurveANT} implies that $M, \beta$ are proportional. Similar argument also shows that  $M, \alpha$ are proportional. Therefore, $\alpha, \beta$ are proportional.
\end{proof}

\begin{lem} \label{propusingfx}
Let $X$ be a compact K\"ahler manifold of dimension $n$.
For $\alpha,\beta$ big nef on $X$, $F(\alpha,\beta)=0$ iff $\alpha$ and $\beta$ are proportional.
\end{lem}

We restrict on nef classes since we do not have the K\"ahler version of \cite{lehmxiaoPosCurveANT} yet.

\begin{proof}
The argument is the same as Lemma \ref{movableusingLX}, while in the last step we conclude proportionality by \cite{fx19} or Theorem \ref{TeissierLorentzian} instead.
\end{proof}

 Motivated by FMP estimate in the convexity setting, we propose the following conjecture.

\begin{conj}\label{StabofBM}
Let $X$ be a compact K\"ahler manifold of dimension $n$.
Let $\alpha,\beta$ be big movable classes on $X$.
Then there exists a constant $c(n)$ depending only on the dimension $n$ such that 
    \begin{equation}\tag{\textbf{FMP}}
        F(\alpha,\beta)\leq  c(n) \sqrt{\sigma(\alpha,\beta)B(\alpha,\beta)} 
    \end{equation}
    where $\sigma(\alpha,\beta)
    =\max\left\{\left(\frac{|\alpha|}{|\beta|}\right)^{\frac{1}{n}},
    \left(\frac{|\beta|}{|\alpha|}\right)^{\frac{1}{n}}\right\}$ is the relative size.
\end{conj}

%**This should be optimal via toric variety**

It is sure that one can formulate the conjecture for big nef classes on a general Lorentzian class.
\begin{conj}\label{StabofBMOmega}
Let $X$ be a compact K\"ahler manifold, and let $\Omega$ be a Lorentzian class of dimension $d$ on $X$.
Assume that $\alpha,\beta$ be big nef classes on $\Omega$,
then there exists a constant $c(d)$ depending only on $d$ such that 
    \begin{equation}\tag{$\textbf{FMP}_\Omega$}
        F_\Omega(\alpha,\beta)\leq  c(d) \sqrt{\sigma(\alpha,\beta)B(\alpha,\beta)} 
    \end{equation}
    where $\sigma(\alpha,\beta)$ is the relative size on $\Omega$ and 
   \begin{equation*}
        F_\Omega(\alpha,\beta)=1-\sup_{\gamma\leq r\alpha,\gamma\leq\beta}\frac{|\gamma|_\Omega}{|\beta|_\Omega}
    \end{equation*}
    where $\gamma$ ranges over all big nef classes such that $\gamma\leq r\alpha, \gamma\leq\beta$. Here $A\leq B$ means that for any $B_i$ nef, 
    $$(B-A)\cdot B_1 \cdot...\cdot B_{d-1}\cdot \Omega \geq 0,$$ 
    hence we emphasis the class $\Omega$ in the notation $F_\Omega$.
\end{conj}

It should be noted that $F_{[X]}$ can be different with $F$ due the difference of the partial orders (see also Remark \ref{diffradii}).

\begin{lem}
Assume that $\Omega$ is strictly Lorentzian and $\alpha, \beta$ are big nef on $\Omega$, then $F_\Omega(\alpha,\beta)=0$ iff $\alpha, \beta$ are proportional.
\end{lem}

\begin{proof}
This follows from Theorem \ref{TeissierLorentzian} (and the equality characterizations in Lemma \ref{logconc}) and the same arguments as in Lemma \ref{movableusingLX} or \ref{propusingfx}.
\end{proof}

\begin{rmk}\label{rmkKTBM}
Without loss of generality, we assume $|\alpha|\geq |\beta|$. By Theorem \ref{deficitscompa}, Conjecture \ref{StabofBMOmega} is equivalent to that: there exists $c(d)>0$ such that
    \begin{equation*}\label{eqfk}
        F_\Omega (\alpha,\beta) \leq c(d) \sqrt{K(\alpha,\beta)}.
    \end{equation*}
Indeed, fix $t>0$, the following two statements are equivalent:
\begin{enumerate}
  \item there exists $c(d)>0$ such that
    $F_\Omega (\alpha,\beta) \leq c(d) K(\alpha,\beta)^t$.
  \item there exists $c(d)>0$ such that
    $F_\Omega (\alpha,\beta) \leq c(d) (\sigma(\alpha,\beta)B(\alpha,\beta))^t$.
\end{enumerate}
The above constants $c(d)$ may be different.
It is obvious from Theorem \ref{deficitscompa} that (2) implies (1). For the converse,
by Theorem \ref{deficitscompa} again, we have some constant $c_1 (d)>0$ such that
\begin{equation}\label{eqBK}
c_1 (d) \sigma (\alpha, \beta)B(\alpha, \beta)\geq  \frac{ K(\alpha,\beta) }{1+ K(\alpha,\beta)}.
\end{equation}
As $ F_\Omega(\alpha,\beta)\leq 1$, to prove (1) we only need to consider the case when $c_1 (d) \sigma (\alpha, \beta)B(\alpha, \beta)<1/2$, yielding that $K(\alpha,\beta)<1$. Hence, using (\ref{eqBK}),
\begin{equation*}
 K(\alpha,\beta)\leq 2 c_1 (n) \sigma (\alpha, \beta)B(\alpha, \beta).
\end{equation*}
This shows that (1) implies (2). These discussions also apply to Conjecture \ref{StabofBM}.
\end{rmk}

While do not have much evidence for the $\frac{1}{2}$-exponent in Conjecture \ref{StabofBMOmega}, from Corollary \ref{KTradii} we obtain the $\frac{1}{2^{d-2}}$-exponent by ``cutting down'' a Lorentzian class of dimension $d$.

\begin{prop}\label{fomegaKT}
Let $X$ be a compact K\"ahler manifold of dimension $n$ and let $\Omega$ be a Lorentzian class of dimension $d$. Fix a K\"ahler class $\omega$, and denote $\Gamma = \omega^{d-2}\cdot \Omega$. Assume that $\alpha,\beta$ are big nef classes on $\Omega$, then there exists a constant $c(d)>0$ such that
 $$F_\Gamma(\alpha,\beta)\leq c(d) K(\alpha,\beta)^{\frac{1}{2^{d-2}}}.$$
    In particular, there exists a constant $c(d)$ such that 
    \begin{equation*}
        F_\Gamma(\alpha,\beta)\leq  c(d) \left(\sigma(\alpha,\beta)B(\alpha,\beta)\right)^{\frac{1}{2^{d-2}}}.
    \end{equation*}
\end{prop}

\begin{proof}
By similar discussion as in Remark \ref{rmkKTBM}, it suffices to compare $F_\Gamma$ and the KT deficit.
Since $K(\alpha,\beta)$ and $F_\Gamma(\alpha,\beta)$ are both $0$-homogeneous with respect to all variables, we may assume $$|\alpha|_\Gamma=|\beta|_\Gamma=1.$$

By the definition of radii on $\Gamma$ (see Definition \ref{radiusonOMEGA}), we have $r_\Gamma (\alpha, \beta)\beta \leq \alpha$ and due to  $|\alpha|_\Gamma=|\beta|_\Gamma=1$, we get $$r_\Gamma (\alpha, \beta) \leq 1, R_\Gamma (\alpha,\beta)\geq 1.$$ Thus, 
\begin{equation}\label{Fomega}
  F_\Gamma(\alpha,\beta) \leq 1 -r_\Gamma (\alpha, \beta)^2 \leq 2(1 - r_\Gamma (\alpha, \beta)).
\end{equation}
Here, in the second inequality we apply the elementary inequality 
$$1 -x^2 \leq 2(1 - x),\ \text{if}\ x\in [0,1].$$ 
On the other hand, by Corollary \ref{KTradii}, 
there is a constant $c(d)>0$ such that 
    \begin{equation}\label{Fomega2}
       K(\alpha,\beta)^{\frac{1}{2^{d-2}}} \geq c(d) \left( 1-\frac{r_\Gamma(\alpha,\beta)}{R_\Gamma(\alpha,\beta)}  \right).
    \end{equation}

Combining (\ref{Fomega}), (\ref{Fomega2}) and $R_\Gamma (\alpha,\beta)\geq 1$ implies the result.
\end{proof}

In the following sections, we focus on FMP type estimates on $X$.

\subsubsection{The $\frac{1}{n}$-exponent via Diskant inequality}

We show that the Diskant inequality gives the $\frac{1}{n}$-exponent.
 
\begin{prop}
    Assume that $X$ is a projective manifold of dimension $n$ and $\alpha,\beta$ are big nef classes on $X$, then $$F(\alpha,\beta)\leq 2n K(\alpha,\beta)^{\frac{1}{n}}.$$
    In particular, there exists a constant $c(n)$ depending only on the dimension $n$ such that 
    \begin{equation*}
        F(\alpha,\beta)\leq  c(n) (\sigma(\alpha,\beta)B(\alpha,\beta))^{1/n}
    \end{equation*}
\end{prop}

\begin{proof}
As noted in Remark \ref{rmkKTBM}, using Theorem \ref{deficitscompa}, we only need to establish the estimate relating $F(\alpha,\beta)$ and $K(\alpha, \beta)$.

    Since $K(\alpha,\beta)$ and $F(\alpha,\beta)$ are both $0$-homogeneous with respect to all variables, we may assume $$|\alpha|=|\beta|=1.$$
    Moreover, we may assume $K(\alpha,\beta)\leq 1$ since $F(\alpha,\beta)\leq 1$. Therefore under the assumptions, $$1\leq \alpha^{n-1}\cdot\beta\leq 2.$$

    The Diskant inequality by \cite{BFJ09} tells
    \begin{equation*}
        (\alpha^{n-1}\cdot\beta)^{\frac{n}{n-1}}-1\geq \Big((\alpha^{n-1}\cdot\beta)^{\frac{1}{n-1}}-r(\alpha,\beta)\Big)^n\geq (1-r(\alpha,\beta))^n.
    \end{equation*}
    Using the elementary inequality 
    \begin{equation*}
        t^{\frac{n}{n-1}}-1\leq 4(t-1), \forall 1\leq t \leq 2
    \end{equation*}
    yields
    $$(1-r(\alpha,\beta))^n\leq 4 K(\alpha,\beta),$$
    which further implies
    $$1-r(\alpha,\beta)\leq 2 K(\alpha,\beta)^{\frac{1}{n}}.$$

    Since $|\alpha|=|\beta|$, we have $r(\alpha,\beta)\leq 1$ which, together with the definition of $r(\alpha,\beta)$, implies that $r(\alpha,\beta)\beta\leq \alpha,\beta$. Therefore,
    $$F(\alpha,\beta)=1-\sup_{\gamma\leq \alpha,\beta}\frac{|\gamma|}{|\beta|}\leq 1-r(\alpha,\beta)^n .$$
    All these yield
    \begin{equation*}
        F(\alpha,\beta)\leq 1-r(\alpha,\beta)^n\leq n(1-r(\alpha,\beta))\leq 2nK(\alpha,\beta)^{\frac{1}{n}}.
    \end{equation*}

    This finishes the proof.
\end{proof}

\subsubsection{Deficit via Calabi-Yau theorem}

While we are not able to prove Conjecture \ref{StabofBM} in full generality yet, for big nef classes, we can prove a metric-norm version by replacing $F(\alpha,\beta)$ with $\widehat{F}(\alpha,\beta)$ defined below.

\begin{defn}
Let $X$ be a compact K\"ahler manifold of dimension $n$ and let $\alpha,\beta$ be big nef classes on $X$. Define
    \begin{equation*}
        \widehat{F}(\alpha,\beta)=\sup_{\langle\widehat{\alpha}^n \rangle=\frac{|\alpha|}{|\beta|}\langle\widehat{\beta}^n\rangle = \Phi}  \frac{1}{|\alpha|} \int_{\Amp(\alpha,\beta)} \left(\frac{\Vert \widehat{\alpha}-(\frac{|\alpha|}{|\beta|})^{1/n}
        \widehat{\beta}\Vert_{\widehat{\beta}}}{\Vert\widehat{\alpha}\Vert_{\widehat{\beta}}}\right)
        \Phi
    \end{equation*}
    where 
    \begin{itemize}
      \item $\Amp(\alpha,\beta) = \Amp(\alpha)\cap \Amp(\beta)$ is the intersection of the ample locus of $\alpha, \beta$, which is a Zariski open subset in $X$.
      \item $\Phi$ ranges over all the smooth volume forms on $X$ such that $\int_X \Phi = |\alpha|$.
      \item $\widehat{\alpha}\in\alpha$ and $\widehat{\beta}\in \beta$ range over all the pairs of positive currents such that $$\langle\widehat{\alpha}^n\rangle=\frac{|\alpha|}{|\beta|}\langle\widehat{\beta}^n\rangle =\Phi,$$ 
    whose existence is guaranteed by \cite{BEGZ10}. Here, $\langle -\rangle$ is the non-pluripolar product of positive currents.
    \item $\Vert -\Vert_{\widehat{\beta}}$ is the pointwise $L^2$-norm on $\Amp(\alpha,\beta)$ induced by $\widehat{\beta}$, which is a K\"ahler metric on the ample locus. 
    \end{itemize}
\end{defn}

It is clear that $\widehat{F}(\alpha, \beta)$ is 0-homogeneous with respect to the variables $\alpha, \beta$.

\begin{lem}
For $\alpha,\beta$ big nef, $\widehat{F}(\alpha, \beta)=0$ iff $\alpha, \beta$ are proportional.
\end{lem}

\begin{proof}
The equality $\widehat{F}(\alpha, \beta)=0$ shows that $$\widehat{\alpha}=\left(\frac{|\alpha|}{|\beta|}\right)^{1/n}
        \widehat{\beta}\ \text{on}\ \Amp(\alpha, \beta).$$ 
        Then using \cite[Lemmas 2.2, 2.3]{fx19} implies that they are actually equal on $X$, finishing the proof. 
\end{proof}

A few words are needed for the last bullet in the above definition, which will be used later. Regarding the equations $$\langle\widehat{\alpha}^n\rangle=\frac{|\alpha|}{|\beta|}\langle\widehat{\beta}^n\rangle =\Phi,$$ 
by \cite{BEGZ10},  the current $\widehat{\alpha}$ (respectively, $\widehat{\beta}$) has minimal singularities and is smooth on the Zariski open set $\Amp(\alpha)$ (respectively, $\Amp(\beta)$). Here, the smoothness applies the seminal work of Yau \cite{Yau78}. Hence, $\widehat{\alpha}$ (respectively, $\widehat{\beta}$) is a K\"ahler metric on $\Amp(\alpha)$ (respectively, $\Amp(\beta)$).

In the special case when $\alpha, \beta$ are K\"ahler, 
the deficit function $\widehat{F}(\alpha,\beta)$ is given by 
    \begin{equation*}
        \widehat{F}(\alpha,\beta)=\sup_{\widehat{\alpha}^n=\frac{|\alpha|}{|\beta|}\widehat{\beta}^n} \frac{1}{|\alpha|} \int_X \left(\frac{\Vert \widehat{\alpha}-(\frac{|\alpha|}{|\beta|})^{1/n}
        \widehat{\beta}\Vert_{\widehat{\beta}}}{\Vert\widehat{\alpha}\Vert_{\widehat{\beta}}}\right)
        \widehat{\alpha}^n
    \end{equation*}
    where $\widehat{\alpha}\in\alpha$ and $\widehat{\beta}\in \beta$ range over all the pairs of K\"ahler metrics such that $$\widehat{\alpha}^n=\frac{|\alpha|}{|\beta|}\widehat{\beta}^n$$ 
    whose existence is guaranteed by Yau's theorem \cite{Yau78}.

The following lemma is taken from \cite[Lemma 2.5]{Figalli10inventiones}.

\begin{lem}\label{figallilem}
Let $0 <\lambda_1 \leq \cdots \leq\lambda_n$ be positive real numbers, and denote 
$$\lambda_A = \frac{1}{n} \sum_{i=1}^n \lambda_i,\ \lambda_G = \left(\prod_{i=1}^n \lambda_i\right)^{1/n},$$ then
$$7n^2 (\lambda_A - \lambda_G)  \geq \frac{1}{\lambda_n}\sum_{i=1} ^n (\lambda_i - \lambda_G)^2.$$
\end{lem}

The main result in this subsection is:

\begin{thrm}\label{KTcontrF}
    If $\alpha,\beta$ are two big nef classes, then there exists a constant $c(n)>0$ such that
    \begin{equation}\label{eqfhatk}
        \widehat{F}(\beta,\alpha) \leq c(n) K(\alpha,\beta)^\frac{1}{2}.
    \end{equation}
    As a consequence,  assume further that $|\alpha|\geq |\beta|$, then there exists a constant $c(n)>0$ such that
    \begin{equation}\label{eqfhatb}
        \widehat{F}(\beta, \alpha)^2 \leq c(n) \left(\frac{|\alpha|}{|\beta|}\right)^{1/n} B(\beta, \alpha).
    \end{equation}
\end{thrm}

The proof of this result is almost the same as \cite{figalli09}, while we apply complex Monge-Amp\`ere equations instead of mass transport method. 

\begin{proof}
We first prove the estimate (\ref{eqfhatk}).
Since both $\widehat{F}(\beta,\alpha)$ and $K(\alpha,\beta)$ are 0-homogeneous with respect to the variables $\alpha, \beta$, without loss of generalities, we may assume that
$$|\alpha|=|\beta|=1.$$

Let $\Phi$ be a smooth volume form with volume 1, and let $\widehat{\alpha}\in\alpha$ and $\widehat{\beta}\in \beta$ be a pair of positive currents such that $$\langle\widehat{\alpha}^n\rangle=\langle\widehat{\beta}^n\rangle=\Phi.$$

At any point $x\in \Amp(\alpha, \beta)$, we can choose a local coordinate system $(z_1,...,z_n)$ such that 
$$\widehat{\beta}=\sqrt{-1}\sum_{i=1}^n \lambda_i dz_i\wedge {d\bar z_i}\text{ and } \widehat{\alpha}=\sqrt{-1}\sum_{i=1}^ndz_i\wedge {d\bar z_i}$$
    with 
    $$0<\lambda_1\leq ...\leq \lambda_n. $$    
    Then $\widehat{\alpha}^n=\widehat{\beta}^n$ implies that $\prod_{i=1} ^n \lambda_i =1$. Thus,
    $$\lambda_n \geq 1, \lambda_1 \leq 1. $$    
    By $\sum_{i=1} ^n \lambda_i ^2 \geq 1$, it is easy to see that
    \begin{equation*}
     \frac{(\sum_{i=1} ^n (\lambda_i -1)^2)^{1/2}}{(\sum_{i=1} ^n \lambda_i ^2)^{1/2}} \leq 2\sqrt{n},
    \end{equation*}
    implying that 
      \begin{equation}\label{upperbdF}
        \widehat{F}(\beta,\alpha)=\sup_{\widehat{\alpha}^n=\widehat{\beta}^n=\Phi} \int_{\Amp(\alpha, \beta)} \frac{(\sum_{i=1} ^n (\lambda_i -1)^2)^{1/2}}{(\sum_{i=1} ^n \lambda_i ^2)^{1/2}}\Phi \leq 2\sqrt{n}.
    \end{equation}

Therefore, to prove the existence of $c(n)$, we can reduce to the case when $K(\alpha, \beta)<1$.

Regarding the KT deficit, we have:    
\begin{align*}
  K(\alpha, \beta)&= \alpha^{n-1} \cdot \beta -1\\
  &=\int_X \langle\widehat{\alpha}^{n-1}\wedge \widehat{\beta}\rangle -1\\
  &=  \int_{\Amp(\alpha, \beta)} \widehat{\alpha}^{n-1}\wedge \widehat{\beta} -1\\
  &=\int_{\Amp(\alpha, \beta)} \left(\frac{\widehat{\alpha}^{n-1}\wedge \widehat{\beta}}{\widehat{\alpha}^{n}} -1 \right) \widehat{\alpha}^{n}\\
  &= \int_{\Amp(\alpha, \beta)} (\lambda_A -\lambda_G) \widehat{\alpha}^{n},
\end{align*}
where in the second equality we use the fact that $\widehat{\alpha}, \widehat{\beta}$ have minimal singularities, and in the third equality we use that non-pluripolar products put no mass on proper analytic subsets. 
Using Lemma \ref{figallilem} and the inequality 
$$(\lambda_n-\lambda_1)^2\leq 2((\lambda_n-1)^2+(1-\lambda_1)^2),$$
we obtain: 
\begin{align*}
  K(\alpha, \beta)&\geq\int_{\Amp(\alpha, \beta)} \frac{1}{7n^2 \lambda_n} \sum_{i=1} ^n (\lambda_i -1)^2 \widehat{\alpha}^{n}\\
  &\geq \frac{1}{14n^2 } \int_{\Amp(\alpha, \beta)}  \frac{ (\lambda_n -\lambda_1)^2}{\lambda_n} \widehat{\alpha}^{n}\\
  &\geq \frac{1}{14n^2 } \left(\int_{\Amp(\alpha, \beta)}  (\lambda_n -\lambda_1) \widehat{\alpha}^{n}\right)^2 \left(\int_{\Amp(\alpha, \beta)} \lambda_n \widehat{\alpha}^{n}\right)^{-1},
\end{align*}
where the last line follows from H\"older inequality, therefore,
\begin{equation}\label{lambda}
 \int_{\Amp(\alpha, \beta)}  (\lambda_n -\lambda_1) \widehat{\alpha}^{n} \leq \sqrt{14}n \left(\int_{\Amp(\alpha, \beta)} \lambda_n \widehat{\alpha}^{n}\right)^{1/2} K(\alpha, \beta)^{1/2}.
\end{equation}

 Using $\lambda_1 \leq 1$, we get that
 \begin{equation*}
 \int_{\Amp(\alpha, \beta)}  \lambda_n  \widehat{\alpha}^{n} \leq 1+ \sqrt{14}n \left(\int_{\Amp(\alpha, \beta)} \lambda_n \widehat{\alpha}^{n}\right)^{1/2} K(\alpha, \beta)^{1/2}
 \end{equation*}
 which, combining with $\lambda_n \geq 1$ and the reduction $K(\alpha, \beta)<1$, implies that 
 \begin{equation}\label{lambdan}
  \int_{\Amp(\alpha, \beta)} \lambda_n  \widehat{\alpha}^{n} \leq c(n) (1 +  K(\alpha, \beta)^{1/2})
 \end{equation}
for some constant $c(n)>0$.

By (\ref{lambda}), (\ref{lambdan}) and again the reduction $K(\alpha, \beta)<1$, we get
\begin{equation}\label{lambda1}
 \int_{\Amp(\alpha, \beta)} (\lambda_n -\lambda_1) \widehat{\alpha}^{n} \leq c(n) K(\alpha, \beta)^{1/2},
\end{equation}
yielding 
\begin{equation}\label{lambda1}
 \int_{\Amp(\alpha, \beta)} \left(\frac{\lambda_n -\lambda_1}{\lambda_n}\right) \widehat{\alpha}^{n} \leq c(n) K(\alpha, \beta)^{1/2},
\end{equation}
since $\lambda_n \geq 1$.

Finally, note that
at each point $x\in {\Amp(\alpha, \beta)}$, we have
    \begin{equation*}
        \frac{1}{\sqrt{2}}(\lambda_n-\lambda_1)\leq \Vert\widehat{\beta}-\widehat{\alpha}\Vert_{\widehat{\alpha}}
        =(\sum_{i=1}^n(\lambda_i-1)^2)^{1/2}\leq \sqrt{n}(\lambda_n-\lambda_1)
    \end{equation*}
and 
\begin{equation*}
        \lambda_n\leq  \Vert\widehat{\beta}\Vert_{\widehat{\alpha}} =(\sum_{i=1}^n\lambda_i^2)^{1/2}\leq \sqrt{n}\lambda_n,
    \end{equation*}
implying that
\begin{equation*}
    \int_{\Amp(\alpha, \beta)}  \left(\frac{\Vert\widehat{\beta}-\widehat{\alpha}\Vert_{\widehat{\alpha}}}
   {\Vert\widehat{\beta}\Vert_{\widehat{\alpha}}}\right) \Phi=\int_{\Amp(\alpha, \beta)}  \left(\frac{\Vert\widehat{\beta}-\widehat{\alpha}\Vert_{\widehat{\alpha}}}
   {\Vert\widehat{\beta}\Vert_{\widehat{\alpha}}}\right) \widehat{\alpha}^{n} \leq c(n) K(\alpha, \beta)^{1/2}.
\end{equation*}

This finishes the proof of (\ref{eqfhatk}).

Next we prove (\ref{eqfhatb}). The argument is similar to the discussion in Remark \ref{rmkKTBM}, which we repeat below. 
Note that $B(\beta, \alpha)=B(\alpha,\beta)$. For $|\alpha|\geq |\beta|$, Theorem \ref{deficitscompa} shows that for some constant $c_1(n)>0$,
\begin{equation*}
  c_1(n)\left(\frac{|\alpha|}{|\beta|}\right)^{1/n} B(\beta, \alpha) \geq \frac{K(\alpha,\beta)}{1+K(\alpha,\beta)}. 
\end{equation*}

Since $\widehat{F}(\alpha, \beta)$ is 0-homogeneous with respect to the variables $\alpha, \beta$, by (\ref{upperbdF}), we have noted that $\widehat{F}(\beta, \alpha)$ has an upper bound depending only on $n$, therefore to bound it from above as required, we may reduce to the case when 
\begin{equation*}
 c_1(n)\left(\frac{|\alpha|}{|\beta|}\right)^{1/n} B(\beta, \alpha)<1/2,
\end{equation*}
which in turn shows that
$$K(\alpha,\beta)\leq 1. $$

Hence,
\begin{equation*}
 2 c_1(n)\left(\frac{|\alpha|}{|\beta|}\right)^{1/n} B(\beta, \alpha) \geq K(\alpha,\beta).
\end{equation*}
Applying (\ref{eqfhatk}) finishes the proof of (\ref{eqfhatb}).

\end{proof}

Theorem \ref{KTcontrF} gives a quantitative refinement of \cite{fx19}.

\subsubsection{Deficit via Newton-Okounkov bodies}

In this section, we prove a variant of Conjecture \ref{StabofBM} via the construction of Newton-Okounkov bodies.

We recall briefly the construction of Newton-Okounkov bodies (see \cite{lazamustNewtonOkou, kk12} for more details). Let $X$ be a normal projective variety of dimension $n$ over an algebraically closed field $k$ (one may just take $k= \mathbb{C}$), and let $L$ be a big line bundle on $X$. Fix an admissible full flag 
\begin{equation*}
  Y_\bullet : X=Y_0 \supsetneq Y_1 \supsetneq...\supsetneq Y_1 \supsetneq Y_n =\{pt\},
\end{equation*}
that is, $Y_k$ is a sequence of normal irreducible subvarieties such that $\codim Y_k = k$ and every $Y_k$ is smooth at the point $Y_n$. Using the flag, one can define a valuation like function $\nu_{Y_\bullet}$ on $H^0 (X, L)$:
\begin{equation*}
  \nu_{Y_\bullet}: H^0 (X, L)\setminus \{0\} \rightarrow \mathbb{Z}^n,\ \nu_{Y_\bullet}(s)=(\nu_{1}(s),...,\nu_{n}(s)).
\end{equation*}
Here, $\nu_{k}(s)$ is inductively defined by $\nu_{k}(s)=\ord_{Y_k} (s_{k-1})$ as follows: we begin with $s_0 =s$, then $\nu_{1}(s)=\ord_{Y_1} (s_{0})$ 
and set 
$$s_1 =\frac{s_0}{s_{Y_1} ^ {\nu_1 (s)}}|_{Y_1} \in H^0 (Y_1, (L - \nu_1 (s)[Y_1])_{|Y_1})$$ 
where $s_{Y_1}$ is the canonical section of the divisor $Y_1$,  
from $s_1$ we can construct $s_2$ in a similar way, $\nu_{2}(s)=\ord_{Y_2} (s_{1})$ and 
$$s_2 =\frac{s_1}{s_{Y_2} ^ {\nu_2 (s)}}|_{Y_2} \in H^0 (Y_2, (L - \nu_1 (s)[Y_1]- \nu_2 (s)[Y_2])_{|Y_2}),$$
continue the process,  
$$ s_k  \in H^0 (Y_k, (L - \nu_1 (s)[Y_1]-...-\nu_k (s)[Y_k])_{|Y_k}),$$
and we obtain $\nu_{Y_\bullet} (s)$. More generally, one can consider the function $\nu_{Y_\bullet}$ associated to the admissible full flag $Y_\bullet$ as the restriction of a valuation on the rational function field
\begin{equation*}
  \nu: k(X)^{\times} \rightarrow \mathbb{Z}^n.
\end{equation*}
Given a rational function $f$ around the point $Y_n =pt$, then $\nu_{k}(f)$ is inductively defined by letting $f_0 =f$, $\nu_{1}(f) =\ord_{Y_1} (f_0)$ and 
\begin{equation*}
 \nu_{k} (f) =\ord_{Y_k} (f_{k-1}),\ f_k = \frac{f_{k-1}}{g_k ^{\nu_{k} (f)}}|_{Y_k}
\end{equation*}
where $g_k =0$ is the local equation of $Y_k$ in $Y_{k-1}$ around $pt$. Note that $H^0 (X, L)$ can be identified as a subspace of $k(X)$. The valuation $\nu$ induced by $\nu_{Y_\bullet}$ has maximal rank $n$. Conversely, by \cite{kuronyaNOBODYvaluation}
 any valuation on $k(X)$ of maximal rank is equivalent to a valuation arising from an admissible flag on a suitable proper birational model of $X$. All known results for Newton-Okounkov bodies defined in
terms of flag valuations generalize to Newton-Okounkov bodies in terms of valuations of
maximal rank, modulo passing to some different birational model.

Let $\nu: k(X)^{\times} \rightarrow \mathbb{Z}^n$ be a valuation of maximal rank. 
Set 
$$\Gamma(L) _m  = \Image (\nu: H^0 (X, mL)\setminus \{0\} \rightarrow \mathbb{Z}^n), $$
then the Newton-Okounkov body of $L$ with respect to the valuation $\nu$ is the compact convex set
\begin{equation*}
  \Delta_{\nu} (L) = \text{closed convex hull} \left\{\bigcup_{m\geq 1} \frac{1}{m} \Gamma(L) _m \right\} \subset \mathbb{R}_{\geq 0}^n.
\end{equation*}
The basic fact is that the construction for $\Delta_{\nu} (L)$ extends to real divisor classes and 
$$|\Delta_{\nu} (L)| =\frac{1}{n!} |L|,$$
and 
$$\Delta_{\nu} (L+M) \supseteq \Delta_{\nu} (L) +\Delta_{\nu} (M). $$

Using Newton-Okounkov bodies and the relative asymmetry index of convex bodies, we introduce a variant of $F(\alpha,\beta)$:

\begin{defn}
    Let $X$ be a normal projective variety and let $\alpha,\beta \in N^1(X)$ be big movable classes, we define 
        \begin{equation*}
            \widetilde{F}(\alpha,\beta)=\sup_{\nu} F(\Delta_{\nu}(\alpha),\Delta_{\nu}(\beta))
        \end{equation*}
        where $\nu$ ranges over all valuations $\nu: k(X)^{\times} \rightarrow \mathbb{Z}^n$ with maximal rank.
\end{defn}

\begin{lem}\label{usingjow}
    For big movable $\alpha,\beta\in N^1(X)$, $\widetilde{F}(\alpha,\beta)=0$ iff $\alpha,\beta$ are proportional.
\end{lem}

\begin{proof}
As $\widetilde{F}$ is 0-homogeneous with respect to all the variables, without loss of generality, we may assume that $|\alpha|=|\beta|$.

It is clear that $\widetilde{F}(\alpha,\beta)=0$, if $\alpha,\beta$ are proportional. 

It remains to prove the other direction.
If $\widetilde{F}(\alpha,\beta)=0$, then by definition, for any admissible flag $Y_\bullet$ on $X$, we have that $\Delta_{Y_\bullet} (\alpha), \Delta_{Y_\bullet} (\beta)$ are homothetic. The assumption $|\alpha|=|\beta|$ shows that 
\begin{equation*}
 \Delta_{Y_\bullet} (\alpha)+ x = \Delta_{Y_\bullet} (\beta)
\end{equation*}
for some translation $x \in \mathbb{R}^n$. 

By taking the flag $Y_\bullet$ sufficiently general such that the point $Y_n$ lies in the intersection of the ample locus of $\alpha, \beta$, we have that 
$$0\in  \Delta_{Y_\bullet} (\alpha) \cap \Delta_{Y_\bullet} (\beta),$$
from which and the fact that 
$$\Delta_{Y_\bullet} (\alpha), \Delta_{Y_\bullet} (\beta) \subset \mathbb{R}_{\geq 0}^n$$
we must have $x=0$. Hence, 
\begin{equation*}
 \Delta_{Y_\bullet} (\alpha)= \Delta_{Y_\bullet} (\beta)
\end{equation*}
By the proof of \cite[Theorem A]{Jow10}, this kind of flags (or valuations) is enough to conclude that $\alpha, \beta$ are proportional since $\alpha,\beta$ are big and movable.

This finishes the proof.
\end{proof}

Now we prove:

\begin{thrm}\label{algFMP}
    Let $X$ be a normal projective variety of dimension $n$ and let $\alpha,\beta \in N^1(X)$ be big classes. Then there exits a constant $c(n)>0$ such that
    \begin{equation*}
        \widetilde{F}(\langle\alpha\rangle,\langle\beta\rangle)\leq c(n)\sqrt{\sigma(\alpha,\beta)B(\alpha,\beta)}.
    \end{equation*}
    In particular, $B(\alpha,\beta)=0$ iff $\langle\alpha\rangle, \langle\beta\rangle$ are proportional.
\end{thrm}

\begin{proof}
Recall that for any valuation $\nu$ with maximal rank, the volume of the Newton-Okounkov body equals the volume of the class:
    \begin{equation}\label{5273}
        n! |\Delta_{\nu}(L)|=|L|.
    \end{equation}
Together with the superadditivity
    \begin{equation*}
        \Delta_{\nu}(\alpha+\beta) \supset \Delta_{\nu}(\alpha)+\Delta_{\nu}(\beta),
    \end{equation*}
    we obtain
    \begin{equation}\label{5271}
        B(\alpha,\beta)=\frac{|\alpha+\beta|^{\frac{1}{n}}}{|\alpha|^{\frac{1}{n}}+|\beta|^{\frac{1}{n}}}-1=\frac{|\Delta_{\nu}(\alpha+\beta)|^{\frac{1}{n}}}{|\Delta_{\nu}(\alpha)|^{\frac{1}{n}}+|\Delta_{\nu}(\beta)|^{\frac{1}{n}}}-1 \geq B(\Delta_{\nu}(\alpha),\Delta_{\nu}(\beta)).
    \end{equation}

By FMP estimate, there exists a polynomial function $c(n)$ depending only on the dimension $n$ such that
    \begin{equation}\label{5272}
        F(\Delta_{\nu}(\alpha),\Delta_{\nu}(\beta))\leq c(n)\sqrt{\sigma(\Delta_{\nu}(\alpha),\Delta_{\nu}(\beta))B(\Delta_{\nu}(\alpha), \Delta_{\nu}(\beta))}.
    \end{equation}
    By (\ref{5273}), we have 
    \begin{equation}\label{5274}
        \sigma(\Delta_{\nu}(\alpha),\Delta_{\nu}(\beta))=\sigma(\alpha,\beta).
    \end{equation}
    Combining (\ref{5271}), (\ref{5272}) and (\ref{5274}), and since $\nu$ is arbitrary, we obtain
    \begin{equation*}
        \sup_{\nu} F(\Delta_{\nu}(\alpha),\Delta_{\nu}(\beta)) \leq c(n)\sqrt{\sigma(\alpha,\beta)B(\alpha,\beta)}.
    \end{equation*}

    Finally, since for any big class $\alpha$, the Newton-Okounkov bodies of $\alpha$ and $\langle \alpha \rangle$ are equal up to a translation, and the relative asymmetry index of convex bodies is translation-invariant, we obtain
    \begin{equation*}
        \widetilde{F}(\langle\alpha\rangle,\langle\beta\rangle) \leq c(n)\sqrt{\sigma(\alpha,\beta)B(\alpha,\beta)}.
    \end{equation*}

    This finishes the proof.
\end{proof}

The argument above shows that the constant $c(n)$ is the same as the constant in FMP estimate, in particular, it has polynomial growth with respect to $n$.

\begin{rmk}
The proof of Theorem \ref{algFMP} also indicates an alternative approach to \cite[Theorem 3.9]{lehmxiaoPosCurveANT} as follows. For $\alpha, \beta$ big divisor classes,  
$ B(\alpha,\beta)=0$ yields that $B(\Delta_{\nu}(\alpha),\Delta_{\nu}(\beta))=0$.
Thus, $\Delta_{\nu}(\alpha),\Delta_{\nu}(\beta)$ are homothetic. Using divisorial Zariski decomposition and \cite{Jow10} implies that $\langle\alpha\rangle, \langle\beta\rangle$ are proportional.
\end{rmk}

The comparison of $\widetilde{F}$ and $F$ motivates the following interesting question on the relation between the partial order of positive classes given by pseudo-effectivity and the the partial order of convex bodies given by the inclusion of sets.

\begin{ques}
Let $X$ be a normal projective variety of dimension $n$ and let $A,B\in N^1(X)$ be big movable divisor classes. Denote $r=\left(\frac{|A|}{|B|}\right)^{1/n}$. Do we have that
\begin{equation*}
  \sup_{C \leq A, C \leq rB} |C|   =\sup_\nu \{|K|: K \subset \Delta_{\nu}(A)\ \&\ K \subset \Delta_{\nu}(rB)\} ?
\end{equation*}
\end{ques}

It is clear that LHS $\leq$ RHS. It is unclear whether they are equal.

\begin{rmk}
By \cite{Jow10}, for $B_1, B_2$ nef, it is easy to see that the inclusion $ \Delta_\nu (B_1) \subset  \Delta_\nu (B_2)$ for a general valuation implies $B_1\leq B_2$. This can be proved as follows. 
The nefness assumption ensures that the asymptotic order of vanishing (sometimes also called minimal vanishing order) $\ord_E (\| B_i\|) =0$ for any prime divisor $E$.
By \cite[Corollary 3.3 and Theorem 3.4]{Jow10}, for an admissible flag $Y_\bullet$ given by a general complete intersection of hypersurfaces from a very ample divisor $A$, 
\begin{equation*}
  [Y_{n-1}] \cdot B_i = \vol (\Delta_{Y_\bullet} (B_i)_{|\textbf{0}^{n-1}}),
\end{equation*}
where $\Delta_{Y_\bullet} (B_i)_{|\textbf{0}^{n-1}}=\{x\in \mathbb{R}: (0,...,0, x)\in \Delta_{Y_\bullet} (B_i)\}$. Then the inclusion of Newton-Okounkov bodies yields that 
\begin{equation*}
  [Y_{n-1}] \cdot B_1 \leq [Y_{n-1}] \cdot B_2.
\end{equation*}
Note that $[Y_{n-1}]=A^{n-1}$ and the same argument also works on any birational model $\pi: \widehat{X}\rightarrow X$, therefore, for any very ample divisor $\widehat{A}$ on $\widehat{X}$,
\begin{equation*}
  (B_2 -B_1)\cdot \pi_* (\widehat{A} ^{n-1})\geq 0.
\end{equation*}
Applying \cite{BDPP13} shows that $B_2-B_1$ is psef.

%For the other direction, assume then $B_2-B_1$ is psef, then there is some $x\in \mathbb{R}_{\geq 0} ^n$ such that
%\begin{equation*}
%  \Delta_\nu (B_1) +x \subset  \Delta_\nu (B_2).
%\end{equation*}
%For $\nu$ general, we have that $0\in \Delta_\nu (B_1) \cap \Delta_\nu (B_2)$, yielding that $x=0$. Therefore, for $\nu$ general, 
%\begin{equation*}
%  \Delta_\nu (B_1) \subset  \Delta_\nu (B_2).
%\end{equation*}

%This finishes the proof.
\end{rmk}

\section*{Acknowledgements}
This work is supported by the National Key Research and Development Program of China (No. 2021YFA1002300) and National Natural Science Foundation of China (No. 11901336). We would like to thank Jixiang Fu and R\'emi Reboulet for some comments on this work. At the finalization of the manuscript, R\'emi asked us questions related to Examples \ref{sec exmple lorent} (2).
The question of whether FMP type estimates hold for positive classes had been asked by Brian Lehmann and some arguments using complex Monge-Amp\`ere equations appeared in a potential collaboration project with Brian several years ago, the second named author would like to thank Brian for many insightful discussions.

\bibliography{reference}
\bibliographystyle{amsalpha}

\bigskip

\bigskip

\noindent
\textsc{Tsinghua University, Beijing 100084, China}\\
\noindent
\verb"Email: hujj22@mails.tsinghua.edu.cn"\\
\verb"Email: jianxiao@tsinghua.edu.cn"\\
\end{document}